\newif\ifpdf
\ifx\pdfoutput\undefined
\pdffalse % we are not running PDFLaTeX
\else
\pdfoutput=1 % we are running PDFLaTeX
\pdftrue \fi

\newif\iffinal
\finalfalse	% Not final version
\finaltrue % Final version

\documentclass[reqno,twoside,11pt]{amsart}
\iffinal\else\usepackage[notref,notcite]{showkeys}\fi
\iffinal\else\IfFileExists{pdfsync.sty}{\usepackage{pdfsync}}{}\fi
\usepackage{cite}
\usepackage{color}
\usepackage{amsmath}
\usepackage{amsfonts}
\usepackage{amssymb}
\usepackage{epsfig}
\usepackage{verbatim}

\usepackage{mathrsfs}
\usepackage{mathpazo}
\usepackage{mathabx}

\DeclareFontFamily{OT1}{eusb}{} \DeclareFontShape{OT1}{eusb}{m}{n} {<5> <6> <7> <8> <9> <10> <11> <12> <14.4> eusb10}{}
\DeclareMathAlphabet{\eusb}{OT1}{eusb}{m}{n}

\DeclareFontFamily{OT1}{eusm}{} \DeclareFontShape{OT1}{eusm}{m}{n} {<5> <6> <7> <8> <9> <10> <11> <12> <14.4> eusm10}{}
\DeclareMathAlphabet{\eusm}{OT1}{eusm}{m}{n}

\DeclareFontFamily{OT1}{eufm}{} \DeclareFontShape{OT1}{eufm}{m}{n} {<5> <6> <7> <8> <9> <10> <11> <12> <14.4> eufm10}{}
\DeclareMathAlphabet{\mathfrak}{OT1}{eufm}{m}{n}

\DeclareFontFamily{OT1}{fraktura}{}
\DeclareFontShape{OT1}{fraktura}{m}{n} {<5> <6> <7> <8> <9> <10> <11> <12> <13> <14.4> [1.1] eufm10}{}
\DeclareMathAlphabet{\fraktura}{OT1}{fraktura}{m}{n}

\DeclareFontFamily{OT1}{cmfi}{} \DeclareFontShape{OT1}{cmfi}{m}{n} {<5> <6> <7> <8> <9> <10> <11> <12> <13> <14.4> [0.9] cmfi10}{}
\DeclareMathAlphabet{\cmfi}{OT1}{cmfi}{b}{n}

\DeclareFontFamily{OT1}{cmss}{} \DeclareFontShape{OT1}{cmss}{m}{n} {<5> <6> <7> <8> <9> <10> <11> <12> <13> <14.4> cmss10}{}
\DeclareMathAlphabet{\cmss}{OT1}{cmss}{m}{n}
\newtheoremstyle{thm}{1.5ex}{1.5ex}{\itshape\rmfamily}{} {\bfseries\rmfamily}{}{2ex}{}

\newtheoremstyle{def}{1.5ex}{1.5ex}{\rmfamily\sl}{} {\bfseries\rmfamily}{}{2ex}{}

\newtheoremstyle{rem}{1.3ex}{1.3ex}{\rmfamily}{} {\itshape}
{} {1.5ex}{}

\newenvironment{proofsect}[1] {\vskip0.1cm\noindent{\rmfamily\itshape#1.}}{\qed\vspace{0.15cm}}%{\newline\vspace{0.15cm}}

\theoremstyle{thm}
\newtheorem{theorem}{Theorem}[section]
\newtheorem{lemma}[theorem]{Lemma}
\newtheorem{proposition}[theorem]{Proposition}

\newtheorem*{Main Theorem}{Main Theorem.}
\newtheorem{corollary}[theorem]{Corollary}

\newtheorem{conjecture}[theorem]{Conjecture}
\newtheorem{problem}[theorem]{Problem}

\theoremstyle{def}
\newtheorem{definition}[theorem]{Definition}

\theoremstyle{rem}
\newtheorem{remark}[theorem]{{Remark}}

\numberwithin{equation}{section}

\renewcommand{\section}{\secdef\sct\sect}
\newcommand{\sct}[2][default]{\refstepcounter{section}
\addcontentsline{toc}{section}
{{\tocsection {}{\thesection}{\!\!\!\!#1\dotfill}}{}}
\vspace{0.7cm}
\centerline{ %\large
\scshape\arabic{section}.\ #1} \nopagebreak \vspace{0.2cm}}
\newcommand{\sect}[1]{
\vspace{0.4cm} \centerline{\large\scshape\rmfamily #1}
\vspace{0.2cm}}

\renewcommand{\subsection}{\secdef\subsct\sbsect}
\newcommand{\subsct}[2][default]{\refstepcounter{subsection}
\addcontentsline{toc}{subsection}
{{\tocsection{\!\!}{\hspace{1.2em}\thesubsection}{\!\!\!\!#1\dotfill}}{}}
\vspace{0.45\baselineskip} {\flushleft\bf
\arabic{section}.\arabic{subsection}~\bf #1.~}
\\*[3mm]\noindent
\nopagebreak}
\newcommand{\sbsect}[1]{\vspace{0.1cm}\noindent
\textbf{#1.~}\vspace{0.1cm}}

\renewcommand{\subsubsection}{%
\secdef \subsubsect\sbsbsect}
\newcommand{\subsubsect}[2][default]{%
\refstepcounter{subsubsection} 
\addcontentsline{toc}{subsubsection}{{\tocsection{\!\!}
{\hspace{3.05em}\thesubsubsection}{\!\!\!\!#1\dotfill}}{}}
\nopagebreak
\vspace{0.15\baselineskip} \nopagebreak {\flushleft\rmfamily
\itshape\arabic{section}.\arabic{subsection}.\arabic{subsubsection}
\ \rmfamily #1\/.}\ }
\newcommand{\sbsbsect}[1]{\vspace{0.1cm}\noindent
\rmfamily \itshape
\arabic{section}.\arabic{subsection}.\arabic{subsubsection} \
\sffamily #1\/.\ }

\iffinal

\else

\fi

\renewcommand{\caption}[1]{%
\vglue0.5cm
\refstepcounter{figure}
\begin{minipage}{0.9\textwidth}\small {\sc Figure~\thefigure. }#1\end{minipage}}

\newcommand{\textd}{\text{\rm d}\mkern0.5mu}

\newcommand{\texte}{\text{\rm  e}\mkern0.7mu}

\newcommand{\Cov}{\text{\rm Cov}}

\newcommand{\BB}{\mathcal B}

\newcommand{\EE}{\mathcal E}
\newcommand{\FF}{\mathcal F}
\newcommand{\GG}{\mathcal G}
\newcommand{\HH}{\mathcal H}

\newcommand{\KK}{\mathcal K}
\newcommand{\LL}{\mathcal L}

\newcommand{\NN}{\mathcal N}

\newcommand{\CalS}{\mathcal S}

\newcommand{\VV}{\mathcal V}

\newcommand{\D}{\mathbb D}
\newcommand{\E}{\mathbb E}

\newcommand{\N}{\mathbb N}

\newcommand{\BbbP}{\mathbb P}
\newcommand{\Q}{\mathbb Q}
\newcommand{\R}{\mathbb R}

\newcommand{\V}{\mathbb V}

\newcommand{\Z}{\mathbb Z}

\newcommand{\twoeqref}[2]{(\ref{#1}--\ref{#2})}

\def\myffrac#1#2 in #3{\raise 2.6pt\hbox{$#3 #1$}\mkern-1.5mu\raise 0.8pt\hbox{$#3/$}\mkern-1.1mu\lower 1.5pt\hbox{$#3 #2$}}

\newcommand{\cc}{\text{\rm c}}

\newcommand{\wt}{\widetilde}

\newcommand{\Cloc}{C_{\text{\rm loc}}}
\newcommand{\Cb}{C_{\text{\rm b}}}
\newcommand{\frakm}{\mathfrak m}

\newcommand{\ttau}{H}

\begin{document}

\vglue-2mm

\title[Deterministic conductance model\hfill]{
\large Homogenization theory of random walks\\among deterministic conductances}
\author[\hfill M.~Biskup]
{Marek~Biskup}
\thanks{\hglue-4.5mm\fontsize{9.6}{9.6}\selectfont\copyright\,\textrm{2025}\ \ \textrm{M.~Biskup. Reproduction for non-commercial use is permitted.}}
\maketitle

\vspace{-5mm}
\centerline{\it Department of Mathematics, UCLA, Los  Angeles, California, USA}
\smallskip

%\smallskip
%\centerline{\version}

\vskip0.5cm
\begin{quote}
\footnotesize \textbf{Abstract:}
%We study asymptotic laws of random walks on~$\mathbb Z^d$ ($d\ge1$) in deterministic reversible environments defined by an assignment of a positive and finite conductance to each edge of~$\mathbb Z^d$. The walk jumps across an edge with probability proportional to its conductance. We identify a deterministic set of conductance configurations for which an Invariance Principle (i.e., convergence in law to Brownian motion under diffusive scaling of space and time) provably holds. This set is closed under translations and zero-density perturbations and carries all ergodic conductance laws subject to certain moment conditions. The proofs are based on martingale approximations whose control relies on the conversion of averages in time and physical space under the deterministic environment to those in a suitable stochastic counterpart. Our study sets up a framework for proofs of ``deterministic homogenization'' in other motions in disordered media.
We study  the  asymptotic  distribution  of random walks on~$\mathbb Z^d$ ($d\ge1$) in deterministic reversible environments defined by an assignment of a positive conductance to each edge of~$\mathbb Z^d$.
We identify a deterministic set of conductance configurations for which  the walk obeys an  Invariance Principle; i.e.,  converges  in law to  a non-degenerate  Brownian motion under diffusive scaling of space and time. This set is closed under translations and zero-density perturbations and carries all ergodic conductance laws subject to certain moment conditions. The proofs  rely  on martingale approximations whose  main step is  the conversion of averages in time and physical space under the deterministic environment to those in a suitable stochastic counterpart. Our study sets up a framework for  ``de-randomized homogenization'' of  other motions in disordered media.

\end{quote}

%

%\tableofcontents\vglue-1cm

\section{Introduction and results}
\noindent
\vglue-1.5mm
\subsection{Main question and outlook}
Let~$d\ge1$ be an integer and let~$E(\Z^d)$ denote the set of ordered nearest-neighbor pairs (a.k.a.\ edges) of~$\Z^d$. Assume that each edge~$e\in E(\Z^d)$ is assigned a positive quantity called the \emph{conductance} of~$e$, which is symmetric upon reversal of the orientation of the edge. We wish to study the asymptotic distribution of random walks on~$\Z^d$ that jump only between nearest-neighbor vertices and that so with probability proportional to the conductance of the traversed edge. Our objective is to find general criteria that provably guarantee a diffusive scaling limit of the walk.

Let $\Omega$ denote the set of admissible conductance configurations. Writing~$\omega$ for a generic element of~$\Omega$, the said random walk is  technically a discrete-time Markov chain with state space $\Z^d$ and transition probability~$\cmss P_\omega$ defined by
\begin{equation}
%\label{}
\cmss P_\omega(x,y):=\frac{\cc_\omega(x,y)}{\pi_\omega(x)}\quad\text{for}\quad\pi_\omega(x):=\sum_{z\colon (x,z)\in E(\Z^d)}\cc_\omega(x,z),
\end{equation}
where $\cc_\omega(e)$, representing the conductance of~$e$ in~$\omega$, is the projection of~$\omega$ on the coordinate indexed by~$e$. The edge-reversal symmetry of the conductances,
\begin{equation}
\label{E:1.2w}
\forall (x,y)\in E(\Z^d)\colon\quad\cc_\omega(x,y)=\cc_\omega(y,x),
\end{equation}
 is then equivalent to~$\pi_\omega$ being a reversible measure. We will use~$X$ to denote a sample path of the random walk and write~$P^x_\omega$ for the law of~$X$ subject to $P^x_\omega(X_0=x)=1$. 
 
%The large-time asymptotic behavior of above random walks has a long history in probability and, via closely related elliptic and parabolic problems in disordered environments, in a subarea of PDE analysis called homogenization theory. This starts from the observation that, for $e\mapsto\cc_\omega(e)$ constant, the chain~$X$ reduces to the simple symmetric random walk. Donsker's Invariance Principle (Donsker~\cite{Donsker}) then asserts convergence in law to an isotropic Brownian motion under the diffusive scaling of space and time. A question is then for what other assignments of the conductances we get (qualitatively) the same result.

 Writing  $C([0,\infty),\R^d)$ for the set of continuous  maps  $[0,\infty)\to\R^d$ endowed with the topology of locally uniform convergence and thus the Borel $\sigma$-algebra $\BB(C([0,\infty),\R^d))$ of measurable sets, the above ``diffusive scaling limit'' is then formalized as:

\begin{definition}[Individual Invariance Principle]
Given~$\omega\in\Omega$ and a sample~$X$ from~$P_\omega^0$, for each natural~$n\ge1$ and real $t\ge0$ set  
\begin{equation}
\label{E:1.3}
B^{(n)}_t:= \frac1{\sqrt n}\Bigl(\,X_{\lfloor tn\rfloor}+(tn-\lfloor tn\rfloor)(X_{\lfloor tn\rfloor+1}-X_{\lfloor tn\rfloor})\Bigr).
\end{equation}
An \emph{Individual Invariance Principle (IIP)} is then said to hold for~$\omega$ if there exists a non-degenerate $d$-dimensional covariance matrix~$\Sigma$ such that, as $n\to\infty$, the law of the process $t\mapsto B_t^{(n)}$ induced by~$P^0_\omega$ on~$(C([0,\infty),\R^d),\BB(C([0,\infty),\R^d)))$ 
tends to that of the Brownian motion $t\mapsto B_t$ with $E B_t=0$ and $\Cov(B_t)=t\Sigma$ for all~$t\ge0$.
\end{definition}

The definition is motivated by Donsker's Invariance Principle (Don\-sker~\cite{Donsker}) which yields an IIP whenever $e\mapsto\cc_\omega(e)$ is a positive constant. The adjective ``Individual'' has been added to draw contrast with similar, but technically distinct, notions of Invariance Principle in this subject area; see Section~\ref{sec-2.1}  for their overview.  

 The  main question  we wish to answer is  for what non-constant assignments $e\mapsto\cc_\omega(e)$ an IIP provably  holds. This has been studied in great detail in the past using the framework of \textit{stochastic} homogenization whose starting point is a conductance configuration drawn at {random} from a law that is invariant and ergodic under translates of~$\Z^d$.  
Our aim here is to depart from, and go beyond, the stochastic homogenization approach and work solely with \emph{deterministic} environments. 

Somewhat to our surprise we find out that, modulo certain restrictions on the conductances that are too large or too small, all that is needed for an IIP to hold is that~$\omega$ obeys the \emph{conclusion} of the Spatial Ergodic Theorem --- namely, that block averages of translates of every continuous local function converge and the limit can be expressed as expectation  under  an ergodic measure. As a consequence, our approach subsumes most of the  previous derivations  of an IIP in  the  stochastic setting.

 In our proofs  we still rely on  standard tools from stochastic homogenization; namely,  martingale approximations and the corrector method so, naturally,  many of   our  arguments  do run parallel to the stochastic approach. However, since we do not start from a probability space, all the needed averaging  has to be derived from spatial averages. This forces us to find ways around certain approximations whose control is soft in the stochastic setting but hard for individual environments.

A take-away message of our work is that, while stochastic homogenization has been extremely useful for the development of the theory of random motions in disordered media, it may in fact be redundant. We believe that our  ``de-randomized  homogenization'' take on the problem is worthy of developing in other models and contexts. An annotated list of some of these is given at the end of the paper.

\subsection{ Framework for de-randomized homogenization}
\label{sec-1.2}\noindent
The statement of our results  requires that we first introduce concepts that underlie  all derivations in this work. Recall that, in stochastic homogenization, we start with a law~$\BbbP$ on the space of admissible conductances and sample conductance configurations from it. In our approach we instead start with an individual conductance configuration~$\omega$ and, imposing suitable mixing assumptions, associate with it a law~$\BbbP_\omega$ that makes the  conclusion of the  Spatial Ergodic Theorem true for~$\omega$.

Recall that
\begin{equation}
%\label{}
\Omega:=\bigl\{\omega\in(0,\infty)^{E(\Z^d)}\colon \text{\eqref{E:1.2w} holds}\bigr\}
\end{equation}
denotes the set of admissible conductance configurations.  Let~$\tau_x\colon\Omega\to\Omega$ denote the ``shift by~$x$'' acting as
\begin{equation}
%\label{}
\tau_x(\omega)(y,z):=\omega(y+x,y+z).
\end{equation}
 We endow~$\Omega$ with the product Euclidean topology and write $C(\Omega)$ for the set of continuous functions~$\Omega\to\R$. We   use~$\Cloc(\Omega)$ to denote  the subset of \emph{local functions}; i.e., those $f\in C(\Omega)$ for which there is a finite~$B\subseteq E(\Z^d)$ such that~$f$ depends only on $\omega$ restricted to~$B$ and, if regarded as a function on~$(0,\infty)^B$, has compact support. For each $r\ge0$, let
\begin{equation}
\label{E:2.2i}
\Lambda_r:=[-r,r]^d\cap\Z^d.
\end{equation}
The aforementioned mixing assumption then comes in:

\begin{definition}
\label{def-1.2}
We say that $\omega\in\Omega$ is \emph{averaging} if the limit
\begin{equation}
\label{E:1.1}
\ell_\omega(f):=\lim_{n\to\infty}\,\frac1{|\Lambda_n|}\!\sum_{x\in\Lambda_n}
f\circ\tau_x(\omega)
\end{equation}
exists for each $f\in \Cloc(\Omega)$.
\end{definition}

Alternative definitions may be considered that yield a similar (or even the same) notion. For instance, we could change the class of generating~$f$'s ---   e.g.,  by requiring \eqref{E:1.1} for all bounded~$f\in C(\Omega)$ or restricting it to just indicators of ``finite patterns'' --- or work with other sequences of domains. The above satisfies our needs and, at the same time, does not formally rule out  issues potentially associated with this concept.

For each averaging~$\omega$, the map $f\mapsto\ell_\omega(f)$ is a continuous positive linear functional on~$\Cloc(\Omega)$ and so one is naturally tempted to realize~$\ell_\omega(f)$ as an integral of~$f$ against a  finite measure. This is possible but getting a non-zero measure requires that we curb the growth and decay of the conductances. Let~$E(\Lambda)$ denote the set of edges incident with at least one vertex in~$\Lambda\subseteq\Z^d$. We then have:

\begin{proposition}
\label{prop-2.2}
Write~$\Cb(\Omega)$ for the set of bounded $f\in C(\Omega)$ and let $\FF$ be the product $\sigma$-algebra on~$\Omega$.
Define the set of  ``tight''  conductance configurations by
\begin{equation}
%\label{E:1.11}
\Omega':=\biggl\{\omega\in\Omega\colon\lim_{\epsilon\downarrow0}\limsup_{n\to\infty}\frac1{n^d}\sum_{e\in E(\Lambda_n)}1_{\R\smallsetminus [\epsilon,1/\epsilon]}(\cc_\omega(e))=0\biggr\}.
\end{equation}
Then for each averaging~$\omega\in\Omega'$, the limit \eqref{E:1.1} exists for all $f\in \Cb(\Omega)$ and there exists a unique probability measure~$\BbbP_\omega$ on $(\Omega,\FF)$ such that
\begin{equation}
\label{E:1.2}
\forall f\in \Cb(\Omega)\colon\, E_{\BbbP_\omega}(f)=\ell_\omega(f).
\end{equation}
Moreover, $\BbbP_\omega$ is translation invariant, i.e., $\forall x\in\Z^d\colon\BbbP_\omega\circ\tau_x^{-1}=\BbbP_\omega$.
\end{proposition}

Unfortunately, the existence of the limits \eqref{E:1.1}  alone  is not sufficient for our needs because it does not imply that samples from~$\BbbP_\omega$ look locally like those from~$\omega$. This is fixed by imposing the following restriction:

\begin{definition}
\label{def-1.4}
An averaging  and tight  $\omega\in\Omega$ is said to be \emph{ergodic} if the associated measure~$\BbbP_\omega$ is ergodic with respect to the translations of~$\Z^d$, i.e., $\BbbP_\omega(A)\in\{0,1\}$ holds for each~$A\in\FF$ satisfying $\tau_x^{-1}(A)=A$ for all~$x\in\Z^d$.
\end{definition}

\noindent
To check that this is not a superfluous concept, assume $d=1$ and let  $\omega$ be such that $\cc_\omega(x,x+1)=1$ if~$x<0$ and $\cc_\omega(x,x+1)=2$ if~$x\ge0$.  Then $\BbbP_\omega$ exists but is a mixture of point-masses on all-one's and all-two's, and so is not ergodic.  A less extreme example (with the same~$\BbbP_\omega$) arises from~$\omega$ with  $\cc_\omega(x,x+1):=1$ when $|x|\in[n^{3/2},(n+1)^{3/2})$ for~$n\in\N$ even and $\cc_\omega(x,x+1):=2$ otherwise.

As it turns out, the requirement in Definition~\ref{def-1.4} can equivalently be cast as a stronger version of the ``averaging'' in Definition~\ref{def-1.2}. Given~$z\in\R^d$, let~$\lfloor z\rfloor$ denote the unique $y\in\Z^d$ such that~$z-y\in[0,1)^d$. We then have:

\begin{proposition}
\label{prop-1.5}
For all~$\omega\in\Omega'$, the property that for each~$f\in\Cloc(\Omega)$ there exists $\ell_\omega(f)\in\R$ such that
\begin{equation} 
\label{E:1.9i}
\lim_{r\to\infty}\limsup_{n\to\infty}\frac1{|\Lambda_n|}\sum_{x\in\Lambda_n}\biggl|\frac1{|\Lambda_r|}\!\sum_{y\in \Lambda_r}f\circ\tau_{x+y}(\omega)-\ell_\omega(f)\biggr|=0
\end{equation}
is equivalent to~$\omega$ being averaging and ergodic.
\end{proposition}

As is easy to check, the condition \eqref{E:1.9i} can equivalently be stated as that in which~$x$ is replaced by~$(2r+1)x$ or~$rx$ in the expression $f\circ\tau_{x+y}(\omega)$. The main point of the condition is to ensure the uniformity of ``averaging'' in most of the translates of~$\Lambda_r$.

\subsection{ Main result}
With  all  the important notions in place, we are ready to state our main result:

\begin{theorem}
\label{thm-2.4}
The set
\begin{equation}
%\label{}
\Omega^\star:=\bigl\{\omega\in\Omega\colon \text{\rm\ averaging,  tight  and ergodic}\bigr\}
\end{equation}
is $\FF$-measurable, translation invariant and obeys $\nu(\Omega^\star)=1$ for each translation-invariant, ergodic probability measure $\nu$ on $(\Omega,\FF)$. In addition, denoting
\begin{equation}
\label{E:2.6w}
\Omega_{p,q}:=\biggl\{\omega\in\Omega\colon\sup_{n\ge1}\,\frac1{n^d}\!\sum_{e\in E(\Lambda_n)}\cc_\omega(e)^p<\infty\,\wedge\,\sup_{n\ge1}\,\frac1{n^d}\!\sum_{e\in E(\Lambda_n)}\cc_\omega(e)^{-q}<\infty\biggr\},
\end{equation}
an IIP  holds for all $\omega\in\Omega^\star\cap\Omega_{p,q}$ with $p,q\in(1,\infty)$ such that
\begin{equation}
\label{E:2.6}
\frac1p+\frac1q<\frac2d \quad\text{\rm if}\quad d\ge2.
\end{equation}
The limiting covariance~$\Sigma$ obeys, and is determined by,
\begin{equation}
\label{E:1.5}
a\cdot\Sigma a =\frac1{E_{\BbbP_\omega}\pi(0)}\,\inf_{\varphi\in \Cloc(\Omega)}E_{\BbbP_\omega}\Biggl(\,\sum_{\begin{subarray}{c}
x\in\Z^d\\(0,x)\in E(\Z^d)
\end{subarray}}
 \cc(0,x) |a\cdot x+\varphi\circ\tau_x-\varphi|^2\Biggr),
\end{equation}
for all~$a\in\R^d$.
\end{theorem}

The upshot of Theorem~\ref{thm-2.4} is that, barring configurations that grow or decay too densely too quickly, being averaging and ergodic is sufficient for an IIP to hold. In particular, our result gives an IIP for a.e.-sample from any translation-invariant ergodic measure~$\nu$ on~$\Omega$ for which 
\begin{equation}
\label{E:2.7i}
\cc(e)\in L^p(\nu)\quad\text{and}\quad \cc(e)^{-1}\in L^q(\nu),\quad e\in E(\Z^d),
\end{equation}
for some~$p,q>1$ satisfying \eqref{E:2.6}. This reproduces the conclusion of Andres, Deuschel and Slowik~\cite{ADS15}  though it  comes short of that in Bella and Sch\"affner~\cite{BSch20} in~$d\ge3$; see Section~\ref{sec-2.1} for  more discussion. 

The set~$\Omega^\star$ of course contains many configurations that are not  typical  samples from any ergodic conductance law. This is because any  (elliptic)  perturbation of a configuration in~$\Omega^\star$ on a zero-density subset of~$E(\Z^d)$ still lies in~$\Omega^\star$. For instance, if $\omega$ is  a non-periodic, elliptic configuration arising as a zero-density  perturbation of an $R$-periodic configuration~$\omega'$, then writing $\Lambda_R':=[0,R)^d\cap\Z^d$,  we have 
\begin{equation}
%\label{}
\BbbP_\omega=R^{-d}\sum_{x\in \Lambda_R'}\delta_{\tau_x(\omega')}
\end{equation}
which is translation invariant and ergodic. The stochastic homogenization approach under the law~$\BbbP_\omega$ gives us an IIP for every translate for~$\omega'$, but since $\BbbP_\omega(\{\omega\})=0$, it has nothing to say about~$\omega$ itself. 

\newcommand{\hate}{\hat{\text{e}}}

A similar situation occurs for aperiodic configurations whose explicit examples in~$\Omega^\star$ can be produced by, e.g., iterating local dynamical rules (such as irrational rotations of the unit circle) or imposing local constraints (as in constructions of aperiodic tilings). To give a concrete example, write $\hate_1,\dots,\hate_d$ for the coordinate vectors in~$\Z^d$, pick $\alpha_1,\dots,\alpha_d\in(0,1)\smallsetminus\Q$ and define~$\omega$ by imposing
\begin{equation}
%\label{}
\cc_\omega\bigl(x+k\hate_i,x+(k+1)\hate_i\bigr):=1+1_{\{\lfloor \alpha_i k\rfloor\text{ even}\}},\quad k\in\Z,
\end{equation}
for each~$i=1,\dots,d$ and and each~$x\in\Z^d$ with~$x\cdot\hate_i=0$. Then~$\omega\in\Omega^\star$  with  the resulting~$\BbbP_\omega$ supported on the  (topological)  closure of $\{\tau_x(\omega)\colon x\in\Z^d\}$ in $\Omega$, which  is a set of  cardinality of the continuum containing only aperiodic configurations.  The stochastic homogenization approach yields an IIP for~$\BbbP_\omega$-a.e.\ sample of the conductances but, since~$\BbbP_\omega$ does not charge singletons, it has nothing to say about~$\omega$ or its translates. Our  ``de-randomized  homogenization'' approach gives an IIP for~$\omega$  itself  and, in fact, \emph{any}  elliptic  zero-density perturbation of \emph{any} element in the closure of $\{\tau_x(\omega)\colon x\in\Z^d\}$.

We remark that the stochastic-homogenization treatment of random walk in aperiodic environments  (which the author learned of from Bartha and Telcs~\cite{BT16}) has been a source of motivation for the setup presented in Section~\ref{sec-1.2}. Attempts for similar axiomatization of convergence assumptions have appeared in different albeit related contexts; e.g., in Blanc, Le Bris and Lions~\cite{BBL}.

\section{Connections and proof ideas}
\noindent
Having stated the results, we now discuss the broader context of this work and give relevant references. Then we will outline the main steps of the proof while highlighting the key differences to the stochastic homogenization approach.

 %In order to state our main result, we need some definitions and notation. To motivate these we note that the lack of explicitness of the set~$\Omega_\nu$ above arises from basically just two sources: the use of the Spatial Ergodic Theorem for space and time averages and reliance on $L^2$-limits in the definition and control of the corrector. Both of these inherently produce an exceptional set which is null but otherwise impossible to control.
%We thus need a way to bypass both problematic steps. 

\subsection{ Earlier work}
\label{sec-2.1}\noindent
As noted  above,  the stochastic homogenization approach to conductance models is generally based on the assumption that~$\omega$ is drawn from a stationary and ergodic law~$\BbbP$ on~$(\Omega,\FF)$. The results along these lines  appeared, roughly, in three stages.  

\medskip\noindent
\textsl{(1) Annealed Invariance Principle}: An early approach based on stochastic homogenization was due to Kipnis and Varadhan~\cite{KV86}, drawing on earlier contributions of Kozlov~\cite{K77,K78,K85}, Papanicolau and Varadhan~\cite{PV79} and K\"unnerman~\cite{Ku83}. Assuming  (besides translation invariance and ergodicity)  the moment conditions
\begin{equation}
\label{E:1.8}
\cc(e),\cc(e)^{-1}\in L^1(\BbbP), \quad e\in E(\Z^d) 
\end{equation}
they gave a proof of weak convergence of~$B^{(n)}$ from \eqref{E:1.3} to Brownian motion under the so called \emph{annealed}, or \emph{averaged}, law on $\Omega\times(\Z^d)^\N$ defined by
\begin{equation}
%\label{}
\textbf{P}(E\times F):=\frac1{\E\pi(0)}\int_E P_\omega^0(F)\pi_\omega(0)\BbbP(\textd \omega).
\end{equation}
The conclusion is then  generally  referred to as the \emph{Annealed Invariance Principle} (AIP).

The homogenization theory enters the argument through the so called \emph{corrector}, which is a random function $\chi\colon\Omega\times\Z^d\to\R^d$ such that $\psi(\omega,x):=x+\chi(\omega,x)$ is $\cmss P_\omega$-harmonic at all~$x\in\Z^d$ and $\BbbP$-a.e.~$\omega$.
The harmonicity is useful as it implies that $n\mapsto\psi(\omega,X_n)$ is a martingale under~$P_\omega^0$. Since the construction of~$\psi$ also gives
\begin{equation}
%\label{}
\E\Bigl(\cc(x,y)\bigl|\psi(\cdot,y)-\psi(\cdot,x)\bigr|^2\Bigr)<\infty,\quad (x,y)\in\E(\Z^d),
\end{equation}
and shows that it satisfies the cocycle condition
\begin{equation}
%\label{}
\psi(\cdot,x+z)-\psi(\cdot,x) = \psi\bigl(\tau_x(\cdot),z\bigr),\quad\BbbP\text{-a.s.}
\end{equation}
 for all~$x,z\in\Z^d$,  one can employ the  technique of  the ``point of view of the particle'' to verify the conditions of the Martingale Functional Central Limit Theorem of Brown~\cite{Brown} and get convergence of the process $t\mapsto n^{-1/2}\psi(\omega,X_{\lfloor tn\rfloor})$ to Brownian motion with covariance given by the formula of the kind \eqref{E:1.5}.

In order to extract the corresponding statement for the random walk $n\mapsto X_n$,  one then needs to show   that the corrector is subdiffusive along the path of the walk, i.e.,
\begin{equation}
\label{E:1.7a}
\max_{1\le k\le n}\frac{\chi(\omega,X_k)}{\sqrt n}\,\,\underset{n\to\infty}{\overset{\textbf{P}}{\longrightarrow}}\,0.
\end{equation}
The annealed law enters here because the proof of this is based on reversibility of the Markov chain on environments $n\mapsto\tau_{X_n}(\omega)$ under the law where~$X$ is sampled from~$P_\omega^0$ for~$\omega$ sampled from
\begin{equation}
\label{E:2.6u}
\Q(\textd\omega):=\frac{\pi_\omega(0)}{\E\pi(0)}\BbbP(\textd\omega).
\end{equation}
(This chain on~$\Omega$ is what defines the term ``point of view of the particle.'') The conclusion \eqref{E:1.7a} gives us actually a bit more than just AIP; indeed, we get that for any bounded continuous function $F\colon C([0,\infty),\R^d)\to\R$,
\begin{equation}
%\label{}
E_\omega^0\bigl(F(B^{(n)})\bigr)\,\,\underset{n\to\infty}{\overset{\Q}\longrightarrow}\,\,E\bigl(F(B)\bigr).
\end{equation}
This is sometimes called ``Invariance Principle in probability''  or ``Semi-quenched Invariance Principle'' (see, e.g., T\'oth~\cite{Toth}).

\medskip\noindent
\textsl{(2) Quenched Invariance Principle}: While the Invariance Principle in probability suffices when averaging is taken with respect to the starting point (as often happens in applications) it has little to say about an IIP for individual samples from~$\BbbP$ and/or the walk started from a fixed point. Attempts to remedy this deficiency gave rise to a sequence of works whose objective is to prove an IIP for~$\BbbP$-a.e.~$\omega$ and~$X$ started from the origin, which  is referred  to as a \emph{Quenched Invariance Principle} (QIP).  The following referencess address  progressively more challenging classes of random environments:
\settowidth{\leftmargini}{(11)}
\begin{enumerate}
\item[(1)] uniformly elliptic environments covered by Sidoravicius and Sznitman~\cite{SS04} (with Quenched CLT settled by Osada~\cite{Osada}, Boivin~\cite{Boivin}, Boivin and Depauw~\cite{Boivin-Depauw}),
\item[(2)] random walk on supercritical percolation cluster (Sidoravicius and Sznitman~\cite{SS04}, Berger and Biskup~\cite{BB07}, Mathieu and Piatnitski~\cite{MP07}),
\item[(3)] i.i.d.\ nearest-neighbor conductance models subject to the restriction that edges with positive conductances percolate (Mathieu~\cite{M08}, Biskup and Prescott~\cite{BP07}, Barlow and Deuschel~\cite{BD12}, Andres, Barlow, Deuschel and Hambly~\cite{ABDH}),
\item[(4)] nearest-neighbor conductances subject to moment conditions \eqref{E:2.7i}
with $p,q=1$ in~$d=1,2$ (Biskup~\cite{B11}) and for $p,q>1$ such that $1/p+1/q<2/d$ (Andres, Deuschel and Slowik~\cite{ADS15}) in~$d\ge3$; these were improved to work under
\begin{equation}
\label{E:1.4w}
\frac1p+\frac1q<\frac2{d-1}\quad\text{if}\quad d\ge3
\end{equation}
by Bella and Sch\"affner~\cite{BSch20},
\item[(5)] long-range conductance models with the first condition in \eqref{E:2.7i} replaced by
\begin{equation}
%\label{}
\sum_{x\in\Z^d} \cc(0,x)|x|^2\in L^p(\BbbP)
\end{equation}
and $p,q>1$ such that $1/p+1/q<2/d$ (Biskup, Chen, Kumagai and Wang~\cite{BCKW}),
\item[(6)] conductance models in subdomains of~$\Z^d$ (Chen, Croydon and Kumagai~\cite{CCK15}) and/or  for the walk  started from a point that ``slides-off'' on the diffusive scale (Rhodes~\cite{Rhodes}).
\end{enumerate}

A key additional ingredient of all proofs in this program (except those in~$d=1,2$) is \emph{elliptic regularity}, which comes either in the form of heat-kernel estimates, exit time estimates or iterative techniques such as Moser iteration. 
 The role of these is  to upgrade the $L^1$-type control
\begin{equation}
%\label{}
\frac1{n^d}\sum_{x\in\Lambda_n}\frac{\bigl|\chi(\omega,x)\bigr|}{n}\,\,\underset{n\to\infty}\longrightarrow\,0, \qquad \BbbP\text{-a.s.},
\end{equation}
which under the second condition in \eqref{E:2.7i} with~$q=1$ comes for free from the construction of~$\chi$, to an $L^\infty$-type control,
\begin{equation}
\label{E:1.9}
\max_{x\in\Lambda_n}\frac{\bigl|\chi(\omega,x)\bigr|}{n}\,\,\underset{n\to\infty}\longrightarrow\,0, \qquad \BbbP\text{-a.s.}
\end{equation}
This is then sufficient to extract an IIP for~$n\to X_n$ from an IIP for~$n\mapsto\psi(\omega,X_n)$, for $\BbbP$-a.e.~$\omega$. The improved conditions in~(4) are basically optimal as, except perhaps for the boundary values in \eqref{E:1.4w}, the convergence \eqref{E:1.9} \emph{fails} in general unless \eqref{E:2.7i} and \eqref{E:1.4w} are in force (Biskup, Chen, Kumagai and Wang~\cite{BCKW}). 

%Note that the above moment conditions coincide with \eqref{E:1.8} in $d=1,2$ but are strictly weaker otherwise. This lends some credence to the conjecture that \eqref{E:1.8} is actually sharp for a QIP in all~$d\ge1$.

\medskip\noindent
\textsl{(3) Quantitative homogenization}:
A separate track of research, running parallel to the above quest for a sharp QIP, focused on quantitative versions of stochastic homogenization. The main objective here is to derive rates of convergence  for the random walk as well as the associated PDEs.  In order to get these rates, one needs to assume additional regularity of the  stochastic environment; e.g., that the conductances are i.i.d.\ or with a good decay of correlations or  obey  a functional inequality, etc.

In the context of the random walk among random conductances, early results of this kind have been derived by, e.g., Mourrat~\cite{Mourrat1,Mourrat2} and Gloria and Mourrat~\cite{Gloria-Mourrat}. There is a large body of literature on quantitative estimates of various relevant quantities such as the corrector; e.g., Gloria and Otto~\cite{GO}, Duerinckx, Gloria and Otto~\cite{DGO} or a sequence of papers by Armstrong, Kuusi and Mourrat culminating in their book~\cite{AKM}. Its recent follow-up by Armstrong and Kuusi~\cite{AK} pushes this further by deriving quantitative homogenization results (for diffusions in the continuum) under a unifying mixing condition called the Concentration for Sums. While the derivations in these works are by and large deterministic, and they meticuously separate the deterministic steps from those where stochastic input is needed, the conclusions remain inherently stochastic (i.e., valid for~$\BbbP$-a.e.\ sample of the environment). 

\subsection{Key obstacles}
The effective outcome of the proofs of IIPs using any of the  above   methods is the assignment~$\nu\mapsto\Omega_\nu$ of a measurable subset of~$\Omega$ to each ergodic conductance law~$\nu$ such that an IIP holds for all~$\omega\in\Omega_\nu$. The set~$\Omega_\nu$ is of full $\nu$-measure but is otherwise inexplicit and its dependence on~$\nu$ is singular; in fact, $\Omega_\mu$ and~$\Omega_{\nu}$ are generically disjoint when~$\mu$ and~$\nu$ are ergodic but not equal. As singletons are $\nu$-null for non-periodic~$\nu$, we have no way to provably  conclude  an IIP for  any given $\omega\in\Omega$  unless it is periodic.

There are two good reasons to be unhappy about this situation. One of them is practical:  Disordered systems encountered in applications are not \emph{a priori} generic samples from a probability law but rather particular instances that determine their own statistics, very much in spirit of our Definition~\ref{def-1.2}. (An engineer will not be interested in the properties of a typical sample from that statistics but rather in that very system that sits on the lab desk.) The second reason is conceptual: As witnessed throughout the proofs of the QIPs,  most of the arguments are deterministic with stochastic tools being called upon only when deterministic arguments run out of steam. It is thus quite reasonable to ask whether a fully deterministic proof is possible.

Fortunately, as an inspection  of the existing proofs of the QIP  reveals, the need for stochastic setting (and lack of explicitness of the set~$\Omega_\nu$) arises from basically just two sources: first, the use of the Spatial Ergodic Theorem for space and time averages and, second, from the reliance on $L^2$-limits in the definition and control of the corrector. Both of these inherently produce an exceptional set which is null but otherwise impossible to control. All we need to do is to find a way to bypass these two steps.

\subsection{Proof outline}
\label{sec-2.2}\noindent
A ``reader's digest'' version of our proofs  opens by the statement of our main  guiding principle: The random walk by itself contains enough randomness, and thus ability to induce averaging, to make stochasticity of the conductance configuration redundant.  
A formal  version  of this comes in  Theorem~\ref{thm-2.1} which asserts that averaging does take place at large temporal scales for observables from the ``point of view of the particle.'' Explicitly, for each~$\omega\in\Omega^\star\cap\Omega_{p,q}$ (with~$p,q$ as in Theorem~\ref{thm-2.4}),
\begin{equation}
\label{E:2.12u}
\frac1n\sum_{k=0}^{n-1}f\circ\tau_{X_k}(\omega)\,\underset{n\to\infty}\longrightarrow\,E_{\Q_\omega}(f),
\qquad\text{\rm in }L^1(P_\omega^0),
\end{equation} 
holds for any~$f\in\Cb(\Omega)$ with~$\Q_\omega$ being related to~$\BbbP_\omega$ as~$\Q$ in \eqref{E:2.6u} is related to~$\BbbP$. 

In stochastic setting, the limit \eqref{E:2.12u} holds~$\BbbP$-a.s.\ thanks to Wiener's Ergodic Theorem and the fact that~$\Q$ is stationary ergodic for the environment Markov chain~$n\mapsto\tau_{X_n}(\omega)$ on~$\Omega$. (We also use that~$\BbbP$ and~$\Q$ are equivalent.) Since this avenue is not available for us, we instead prove \eqref{E:2.12u} by relying on heat-kernel upper bounds that Andres, Deuschel and Slowik~\cite{ADS16} proved to hold for each~$\omega\in\Omega_{p,q}$; see Proposition~\ref{prop-ADS}. (This is where the restriction on~$p$ and~$q$ comes from.) 

The heat-kernel bounds allow us to effectively dominate averages with respect to the time and position of the walk by block averages in space. Explicitly, for all~$\omega\in\Omega^\star\cap\Omega_{p,q}$ there exists  $c(\omega)\in(0,\infty)$ such that
\begin{equation}
\label{E:2.13}
\limsup_{n\to\infty}\frac1n\sum_{k=0}^{n-1}E_\omega^0\bigl(f(X_k)\bigr)\le c(\omega)\,\limsup_{n\to\infty}\frac1{|\Lambda_n|}\sum_{x\in\Lambda_n}f(x)\pi_\omega(x),
\end{equation}
holds for all~$f\colon\Z^d\to[0,\infty)$; see the ``Conversion Lemma'' (Lemma~\ref{lemma-3.3a}). Since expectation with respect to~$P_\omega^0$ is taken on the left, the best we can hope to get from this is convergence in $L^1(P_\omega^0)$  in~\eqref{E:2.12u}. 

For~$f$ of the form~$f(x):=h\circ\tau_x(\omega)$ for some~$h\in\Cb(\Omega)$, the limit on the right-hand side of \eqref{E:2.13} is computed using the fact that~$\omega$ is averaging,
\begin{equation}
%\label{}
\limsup_{n\to\infty}\frac1{|\Lambda_n|}\sum_{x\in\Lambda_n}h\circ\tau_x(\omega)\pi_\omega(x)=E_{\BbbP_\omega}\bigl(\pi(0)\bigr)\,E_{\Q_\omega}(h).
\end{equation}
(Technically,  since  $\pi(0)$ is not bounded, we  also  need to strengthen the convergence \eqref{E:1.1} beyond bounded~$f$'s.) This is the point where spatial averages  morph  into stochastic averages  thus giving us access to  standard calculations from the stochastic homogenization approach. Still, as \eqref{E:2.13} is just a bound, to get an actual limit in \eqref{E:2.12u} we have to  first ``remove the mean,'' which we achive by way of  a martingale approximation. It is this  step  where the assumption that~$\omega$ is ergodic  crucially  enters; see the final part of the proof of Theorem~\ref{thm-2.1}.

With the averaging \eqref{E:2.12u} settled, we  harvest  our first limit results about the random walk. Indeed, along a very standard argument, the averaging alone implies vanishing speed (a.k.a.\ Law of Large Numbers) in all spatial dimensions (see Corollary~\ref{cor-3.2}). With the convergence \eqref{E:1.1} extended to  a suitable class of  unbounded functions, it also yields an IIP in~$d=1$. The latter relies on the convenient fact that the one-dimensional corrector can be constructed explicitly and its sublinearity deduced from block averaging for deterministic conductance configurations alone.

Our proof of an IIP in dimensions $d\ge2$ parallels known arguments from homogenization theory but the specifics are quite different. We do start along a route familiar from Kipnis and Varadhan~\cite{KV86}: A martingale approximation is invoked to write
\begin{equation}
\label{E:2.15u}
X_n = X_0+\bigl[\chi_\epsilon(\omega)-\chi_\epsilon\circ\tau_{X_n}(\omega)\bigr]+
\biggl[\,\sum_{k=0}^{n-1}\epsilon\chi_\epsilon\circ\tau_{X_k}(\omega)\biggr]+M_n^{(\epsilon)},
\end{equation}
where~$M_n^{(\epsilon)}$ is a martingale and~$\chi_\epsilon$ is an approximation of the corrector obtained by solving a ``massive'' Poisson equation with ``mass''~$\epsilon$; see \twoeqref{E:5.2r}{E:3.7w}. 
The intuition we have from stochastic homogenization is that all the terms involving~$\chi_\epsilon$ are small for~$\epsilon$ small and the martingale thus captures the bulk of the fluctuations of the walk. 

This intuition notwithstanding, as we  are not  in a stochastic setting, we cannot employ the usual argument (based, for instance, on spectral theory in~\cite{KV86}) to take~$\epsilon\downarrow0$, which is what is  needed  to eliminate the additive term and replace~$\chi_\epsilon(\tau_x(\omega))-\chi_\epsilon(\omega)$ by~$\chi(\omega,x)$; see \eqref{E:5.9i}. Instead we invoke a martingale approximation \emph{one more time}, this time localized in the physical space, to write
\begin{equation}
\label{E:2.16u}
\sum_{k=\lfloor\delta n\rfloor}^{n\wedge\ttau_{\Lambda_r}}\epsilon\chi_\epsilon\circ\tau_{X_{k-1}}(\omega) =\theta_{r,\epsilon}(\omega,X_{\lfloor \delta n\rfloor\wedge\ttau_{\Lambda_r}})-\theta_{r,\epsilon}(\omega,X_{n\wedge\ttau_{\Lambda_r}})+\wt M_n^{(r,\epsilon)},
\end{equation}
where $\wt M_n^{(r,\epsilon)}$ is another martingale and~$H_{\Lambda_r}$ is the first exit time from~$\Lambda_r$, for $r:=\sqrt n/\delta$ with~$\delta>0$ small. Here~$\theta_{r,\epsilon}(\omega,\cdot)$ is a kind of ``second-order corrector'' being defined as a solution of a (massless) Poisson equation in~$\Lambda_r$ with ``charge density''~$x\mapsto\epsilon\chi_\epsilon\circ\tau_x(\omega)$ and Dirichlet boundary conditions outside~$\Lambda_r$; see \eqref{E:3.24}.

The rewrite \eqref{E:2.16u} does kill the pesky additive term in \eqref{E:2.15u} as soon as we can control the exit time~$H_{\Lambda_r}$.  Here we draw on  earlier work of  Chen, Kumagai, Wang and the present author~\cite{BCKW}  that supplies most of the analytic estimates  assuming $\omega\in\Omega_{p,q}$ for~$p$ and~$q$ as above.  (The needed  adaptations  are given in Section~\ref{sec7}.) It remains to control the objects on the right of \eqref{E:2.16u}.
As it turns out, this reduces to homogenization of the Dirichlet energy  of the potential  associated with the ``charge density'' $x\mapsto\epsilon\chi_\epsilon\circ\tau_x(\omega)$ in~$\Lambda_r$ (see Theorem~\ref{thm-3.1} and Proposition~\ref{prop-4.4}). 

The above steps effectively dominate the limiting variance of the~$\chi_\epsilon$-dependent terms in \eqref{E:2.15u} by the inner product
\begin{equation}
\label{E:2.17i}
\bigl\langle\epsilon\chi_\epsilon,(-\LL)^{-1}\epsilon\chi_\epsilon\bigr\rangle_{L^2(\Q_\omega)},
\end{equation}
where~$\LL$ is the generator of the environment Markov chain; see \eqref{E:4.1w}. A standard calculation based on the Spectral Theorem shows that  the inner product in \eqref{E:2.17i}  tends to zero as~$\epsilon\downarrow0$; see \eqref{E:5.34}. Controlling variances  of the terms in \eqref{E:2.15u}  suffices because the aforementioned exit time estimate also supplies tightness of~$B^{(n)}$ (see Proposition~\ref{prop-2.7}) and thus reduces an IIP to convergence of finite dimensional distributions.

Summarizing, while we still  rely on  ideas and tools from stochastic homogenization, the stochastic setting is not present from the outset but rather arises naturally inside proofs. We do not attempt to define the actual corrector for individual conductance configurations in~$d\ge2$ as we think that it is unlikely to exist and/or be regular for a sufficiently large subset of~$\Omega$. This is related to the fact that the need to evaluate functions at specific conductance configurations forces us to base our proofs on uniform convergence of bounded continuous functions rather than~$L^2$-convergence. The $L^2$-calculus typical for stochastic homogenization is still used, but only after the averages in the physical space have been converted to their stochastic counterpart. 

\subsection{ Organization}
%\smallskip
The remainder of this paper is  structured  as follows. Section~\ref{sec3} gives  proofs of Propositions~\ref{prop-2.2} and~\ref{prop-1.5} providing the connection between spatial averages and the associated stochastic framework. A number of preliminary results are also proved  that are useful throughout the rest of the  paper.  Section~\ref{sec4} proves averaging at large temporal scales for functions of the environment  ``observed''  by the random walk  thus making a connection to the ``point of view of the particle.''  Section~\ref{sec5} is devoted to homogenization of finite-volume Dirichlet energies  that drive our control of~$\chi_\epsilon$-dependent terms in \eqref{E:2.15u}.  The proof of the main theorem is completed in Section~\ref{sec6}. Section~\ref{sec7} contains the proof of the aforementioned exit time estimate which is largely drawn from~\cite{BCKW}. Section~\ref{sec8} gives conjectures and open problems along with a general discussion of other models where we think a similar approach could be attempted.

\section{Preliminaries}
\label{sec3}\noindent
Here we prove preliminary results that set up the stochastic counterpart of the physical space for averaging and ergodic conductance configurations. We then list facts from  heat-kernel technology that will be needed in our control of space and time averaging of the walk. The conclusions in this section are largely routine but they will be frequently called upon in later arguments.

\subsection{Stochastic framework}
We start by developing the stochastic-homogenization framework associated with admissible conductance configurations. In particular, we give proofs of Propositions~\ref{prop-2.2}  and~\ref{prop-1.5} as well as  the properties of~$\Omega^\star$ stated in Theorem~\ref{thm-2.4}.

\begin{proofsect}{Proof or Proposition~\ref{prop-2.2}}
The proof is quite simple if we limit ourselves to uniformly elliptic conductances. Indeed, this amounts to replacing~$\Omega$ by $\Omega\cap [a,b]^{E(\Z^d)}$ (for some~$a,b\in(0,\infty)$ with $a<b$) which is compact in the product topology. The set of local functions is then dense in the set of all continuous functions and the limit \eqref{E:1.1} extends readily to these. The claim then follows from the Riesz Representation Theorem  and the fact that $f\mapsto\ell_\omega(f)$ is a continuous positive linear functional on~$C(\Omega\cap [a,b]^{E(\Z^d)})$. 

Dealing with non-elliptic cases requires a number of additional steps  needed, mainly, thanks to the fact that local functions do not have compact support.  We will apply the Riesz Representation Theorem only to finite-dimensional projections and boost the conclusion to the whole space using the Kolmogorov Extension Theorem.

Pick~$\omega\in\Omega'$ averaging. Given a non-empty finite symmetric set $B\subseteq E(\Z^d)$,  denote  $\Omega_B:=\{\omega\in(0,\infty)^B\colon \text{\eqref{E:1.2w} holds in }B\}$ and let $\Pi_B\colon\Omega\to\Omega_B$ be the projection on the coordinates in~$B$. Write $C(\Omega_B)$ be the set of continuous functions~$\Omega_B\to\R$ and denote
\begin{equation}
%\label{}
C_0(\Omega_B) :=\Bigl\{f\in C(\Omega_B)\colon\bigl(\forall\epsilon>0\colon \{|f|\ge\epsilon\}\text{ compact in }\Omega_B\bigr)\Bigr\}
\end{equation}
and
\begin{equation}
%\label{}
C_\cc(\Omega_B) :=\Bigl\{f\in C(\Omega_B)\colon \text{ closure of }\{|f|>0\}\text{ compact in }\Omega_B\Bigr\}.
\end{equation}
Observe that then
\begin{equation}
%\label{}
\bigl\{f\circ\Pi_B\colon f\in C_\cc(\Omega_B)\bigr\}\subseteq\Cloc(\Omega).
\end{equation}
As~$\omega\in\Omega$ is averaging, \eqref{E:1.1} applies to $f\circ\Pi_B$ for all $f\in C_\cc(\Omega_B)$ and, using that $C_\cc(\Omega_B)$ is dense in $C_0(\Omega_B)$, to all $f\in C_0(\Omega_B)$ as well.

The map $f\mapsto\ell_\omega(f\circ\Pi_B)$ is a continuous positive linear functional on $C_0(\Omega_B)$, which is the space of continuous functions on a locally-compact Hausdorff space that vanish at infinity. By the Riesz Representation Theorem, there exists a finite Borel measure $\mu_B$ on~$\Omega_B$ such that
\begin{equation}
\label{E:3.1t}
\forall f\in C_0(\Omega_B)\colon\,\,\ell_\omega(f\circ\Pi_B) = \int_{\Omega_B} f\textd\mu_B.
\end{equation}
In light of $|\ell_\omega(f)|\le\sup|f|$ we have $\mu_B(\Omega_B)\le1$. 

Given $\epsilon\in(0,1/2)$, pick~$g_\epsilon\in C_0(\Omega_B)$ with $0\le g_\epsilon\le1$ and such that $g_\epsilon=1$ on $\bigcap_{e\in B}\{2\epsilon\le\cc(e)\le 1/\epsilon\}$ and $g_\epsilon=0$ on $\Omega\smallsetminus\bigcap_{e\in B}\{\epsilon\le\cc(e)\le 2/\epsilon\}$. Then 
\begin{equation}
\label{E:3.4i}
0\le 1-g_\epsilon\circ\Pi_B\le\sum_{e\in B}1_{\R\smallsetminus[\epsilon,2/\epsilon]}(\cc(e)).
\end{equation}
 Let~$k\ge1$ be  such that $B\subseteq E(\Lambda_k)$. Then 
\begin{equation}
\label{E:3.2}
0\le 1-\ell_\omega(g_\epsilon\circ\Pi_B)\le |B|\,\limsup_{n\to\infty}\frac1{|\Lambda_n|}\sum_{e\in E(\Lambda_{n+k})}1_{\R\smallsetminus [\epsilon,2/\epsilon]}(\cc_\omega(e))
\end{equation}
with the right-hand side vanishing as $\epsilon\downarrow0$ in light of $\omega\in\Omega'$. The Bounded Convergence Theorem applied to \eqref{E:3.1t} with $f:=g_\epsilon$ then gives $\mu_B(\Omega_B)=1$.

Let $\{B_n\colon n\ge1\}$ be  an increasing sequence of  symmetric sets with~$\bigcup_{n\ge1}B_n=E(\Z^d)$. Since \eqref{E:3.1t} shows that the  probability  measures $\{\mu_{B_n}\colon n\ge1\}$ form a consistent family, the Kolmogorov Extension Theorem  shows that these are restrictions of  a unique probability measure~$\BbbP_\omega$ on~$(\Omega,\FF)$ such that
\begin{equation}
\label{E:3.3}
\forall f\in \Cloc(\Omega)\colon\ell_\omega(f)=E_{\BbbP_\omega}(f),
\end{equation}
where the restriction to local functions permits us to identify~$f$ with $f'\circ\Pi_B$, for a suitable~$B$ and $f'\in C_\cc(\Omega_B)$, and then apply \eqref{E:3.1t}.

In order to extend \eqref{E:3.3} to all bounded continuous functions, let~$f\in\Cb(\Omega)$ and pick $\epsilon>0$. By continuity of~$f$ in the product topology, there exists a finite $B\subseteq E(\Z^d)$ and a bounded $f'\in C(\Omega_B)$ such that $\sup|f-f'\circ\Pi_B|<\epsilon$ and $\sup|f'|\le2\sup|f|$. Set
\begin{equation}
%\label{}
f_\epsilon:=(f' \cdot g_\epsilon)\circ\Pi_B
\end{equation}
for $g_\epsilon$ as in \eqref{E:3.4i} and note that $f_\epsilon\in\Cloc(\Omega)$. Using \eqref{E:3.4i} we get
\begin{equation}
\label{E:3.9a}
\sum_{x\in\Lambda_n}\bigl|f\circ\tau_x(\omega)-f_\epsilon\circ\tau_x(\omega)\bigr|<\epsilon|\Lambda_n|+2|B|(\sup|f|)\sum_{e\in E(\Lambda_{n+k})}1_{\R\smallsetminus [\epsilon,2/\epsilon]}(\cc_\omega(e)).
\end{equation}
Dividing by $|\Lambda_n|$, taking~$n\to\infty$ followed by $\epsilon\downarrow0$, the fact that $\ell_\omega(f_\epsilon)$ exists for all $\epsilon>0$ implies that so does~$\ell_\omega(f)$ and that $\ell_\omega(f_\epsilon)\to \ell_\omega(f)$ as~$\epsilon\downarrow0$. A similar argument applied on the right of \eqref{E:3.3} gives $E_{\BbbP_\omega}(f_\epsilon)\to E_{\BbbP_\omega}(f)$, and so we get \eqref{E:1.2} as desired. The translation invariance of~$\BbbP_\omega$ is then a consequence of $\ell_\omega(f\circ\tau_x)=\ell_\omega(f)$.
\end{proofsect}

 We are now also able to give:

\begin{proofsect}{Proof or Proposition~\ref{prop-1.5}}
Using Jensen's inequality to pass the averaging over~$x\in\Lambda_n$ inside the absolute value shows that \eqref{E:1.9i} forces~$\omega$ to be averaging. Under~$\omega\in\Omega'$ the property \eqref{E:1.9i} is thus equivalent to $\omega$ being averaging and
\begin{equation}
\label{E:1.10i}
\lim_{r\to\infty}E_{\BbbP_\omega}\biggl|\frac1{|\Lambda_r|}\!\sum_{y\in \Lambda_r}f\circ\tau_{y}-E_{\BbbP_\omega}(f)\biggr|=0
\end{equation}
being valid for all~$f\in\Cloc(\Omega)$.
Assuming $\BbbP_\omega$ ergodic, \eqref{E:1.10i} is implied by the Spatial Ergodic Theorem. Conversely, since~$\Cloc(\Omega)$ is dense in~$L^1(\BbbP_\omega)$, \eqref{E:1.10i} extends to all $f\in L^1(\BbbP_\omega)$ and, for~$f:=1_A$ with~$A$ invariant under the shifts, gives $1_A=\BbbP_\omega(A)$ $\BbbP_\omega$-a.s., proving~$\BbbP_\omega$ to be ergodic.
\end{proofsect}

The argument \eqref{E:3.9a} extends the convergence in \eqref{E:1.1} from local functions to all boun\-ded continuous functions using only the ``tightness'' assumption built into the definition of~$\Omega'$. The reader may wonder whether the convergence applies to a larger class of functions if a stricter condition is assumed on~$\omega$. Here is a criterion in this vain:

\begin{lemma}
\label{lemma-3.1}
Suppose that~$f\in C(\Omega)$ is such that, for some $h\colon\Omega\to[0,\infty)$,
\begin{equation}
\label{E:3.10i}
\forall\omega'\in\Omega\colon\,\, \bigl| f(\omega')\bigr|\le h(\omega').
\end{equation}
Then $f\in L^1(\BbbP_\omega)$ and \eqref{E:1.1} holds for all $\omega\in\Omega'$ such that
\begin{equation}
\label{E:3.11i}
\sup_{n\ge1}\,\frac1{n^d}\sum_{x\in\Lambda_n}[h\circ\tau_x(\omega)]^\alpha<\infty
\end{equation}
for some $\alpha>1$.
\end{lemma}

\begin{proofsect}{Proof}
Given  $f\in C(\Omega)$, assume that $h\colon\Omega\to[0,\infty)$ and $\omega\in\Omega'$ are such that \twoeqref{E:3.10i}{E:3.11i} hold. We will first prove that $f\in L^1(\BbbP_\omega)$. To this end, given $r>0$, define $f_r\colon\Omega\to\R$ by
\begin{equation}
\label{E:3.12i}
f_r(\omega'):=\max\bigl\{\min\{\,f(\omega'),r\},-r\bigr\}.
\end{equation}
 Then $f_r\in\Cb(\Omega)$ and so \eqref{E:1.1} applies. The assumption \eqref{E:3.10i} combined with Jensen's inequality give
\begin{equation}
\label{E:3.13i}
E_{\BbbP_\omega}\bigl(|f_r|\bigr)=\lim_{n\to\infty}\,\frac1{|\Lambda_n|}\sum_{x\in\Lambda_n}|f_r|\circ\tau_x(\omega)
\le\limsup_{n\to\infty}\biggl(\frac1{|\Lambda_n|}\sum_{x\in\Lambda_n}[h\circ\tau_x(\omega)]^\alpha\biggr)^{1/\alpha}
\end{equation}
 showing that, under  \eqref{E:3.11i}, the left-hand side is bounded uniformly in~$r\ge1$. Taking $r\to\infty$, the Monotone Convergence Theorem gives $E_{\BbbP_\omega}(|f|)<\infty$.

Next
note that \eqref{E:3.10i} gives
\begin{equation}
\label{E:3.14}
\bigl|f(\omega')-f_r(\omega')\bigr|\le h(\omega')g_r(\omega')\quad\text{with}\quad g_r(\omega'):=\frac{2|f(\omega')|}{r+|f(\omega')|}
\end{equation}
for all~$\omega'\in\Omega$. The H\"older inequality then shows
\begin{equation}
\label{E:3.15i}
\begin{aligned}
 \biggl|\,\frac1{|\Lambda_n|}\sum_{x\in\Lambda_n} f\circ\tau_x(&\omega)-\frac1{|\Lambda_n|}\sum_{x\in\Lambda_n} f_r\circ\tau_x(\omega)\biggr|\\
\\*[-5mm]&\le \biggl(\frac1{|\Lambda_n|}\sum_{x\in\Lambda_n} [h\circ\tau_x(\omega)]^\alpha\biggr)^{1/\alpha}
\biggl(\frac1{|\Lambda_n|}\sum_{x\in\Lambda_n} [g_r\circ\tau_x(\omega)]^{\alpha'}\biggr)^{1/\alpha'},
\end{aligned}
\end{equation}
where $\alpha'$ is the H\"older conjugate of~$\alpha$. 
Under \eqref{E:3.11i}, the first term on the right is bounded uniformly in~$n\ge1$. Noting that~$g_r\in\Cb(\Omega)$,  as~$n\to\infty$ the second term tends to
\begin{equation}
%\label{}
\epsilon_r:=[E_{\BbbP_\omega}(g_r^{\alpha'})]^{1/\alpha'}
\end{equation}
by Proposition~\ref{prop-2.2}. Since~$0\le g_r\le2$ and $g_r\to0$ pointwise as $r\to\infty$, the Bounded Convergence Theorem gives $\epsilon_r\to0$ as $r\to\infty$. Since also $f_r\in\Cb(\Omega)$, Proposition~\ref{prop-2.2} yields $|\Lambda_n|^{-1}
%n^{-d}\sum_{x\in\Delta_n(v)} 
\sum_{x\in\Lambda_n}
f_r\circ\tau_x(\omega)\to E_{\BbbP_\omega}(f_r)$. As~$r\to\infty$, the latter expectation tends to~$E_{\BbbP_\omega}(f)$ by the Dominated Convergence Theorem enabled by~$f\in L^1(\BbbP_\omega)$, $|f_r|\le |f|$ and the pointwise convergence~$f_r\to f$. This gives the claim.
\end{proofsect}

Next we will address the properties of the set~$\Omega^\star$ listed in Theorem~\ref{thm-2.4}.
Recall that~$\FF$ denotes the product $\sigma$-algebra on~$\Omega$. 

\begin{proposition}
\label{prop-3.1}
The set $\Omega^\star$ is measurable (i.e., a member of~$\FF$) and translation invariant. Moreover, $\nu(\Omega^\star)=1$ for each shift-invariant, ergodic probability measure $\nu$ on $(\Omega,\FF)$.
\end{proposition}

\begin{proofsect}{Proof}
Note that~$\Omega'$ is measurable and translation invariant. Abbreviate
\begin{equation}
%\label{}
\Omega'':=\{\omega\in\Omega'\colon\text{ averaging}\}.
\end{equation}
Approximation arguments show that $\omega$ is averaging if and only if the limit $\ell_\omega(f)$ exists for all~$f$ in a dense subset of~$\Cloc(\Omega)$.  Since the latter space  is separable and  the set  $\{\omega\in\Omega'\colon\ell_\omega(f)\text{ exists}\}$ is measurable and translation invariant for each~$f\in\Cloc(\Omega)$,  we conclude that  $\Omega''$ is measurable and translation invariant as well. 

In order to prove that $\Omega^\star$ is measurable,  let~$\CalS$ be  a countable dense subset of $\{f\in \Cloc(\Omega)\colon 0\le f\le1\}$. A translation-invariant measure $\mu$ on $(\Omega,\FF)$ then fails to be ergodic if and only if there is $\delta\in(0,1/2)$  such that,  for each $\epsilon>0$, there  exists  $f\in\CalS$  with
\begin{equation}
%\label{}
\delta<E_\mu(f)<1-\delta\,\text{ but }\, E_\mu\bigl(|f\circ\tau_x-f|\bigr)<\epsilon\,\text{ whenever }\, (0,x)\in E(\Z^d).
\end{equation}
Enumerate~$\CalS$ into a sequence $\{f_k\}_{k\ge1}$. Discretizing $\delta:=2^{-m}$ and $\epsilon:=2^{-n}$ yields
\begin{equation}
%\label{}
\{\omega\in\Omega''\colon \text{ NOT ergodic}\} 
= 
\bigcup_{m\ge1}\bigcap_{n\ge1}\bigcup_{k\ge1}\,\,\bigcap_{x\colon|x|=1}
A_{k,m,n}(x),
\end{equation}
where
\begin{equation}
%\label{}
A_{k,m,n}(x):=\Bigl\{\omega\in\Omega''\colon\,2^{-m}<E_{\BbbP_\omega}(f_k)<1-2^{-m}\,\wedge\,
E_{\BbbP_\omega}|f_k\circ\tau_x-f_k|<2^{-n}\Bigr\}.
\end{equation}
As $\omega\mapsto 1_{\Omega''}(\omega)E_{\BbbP_\omega}(f)=1_{\Omega''}(\omega)\ell_\omega(f)$ is measurable for each $f\in \Cloc(\Omega)$, we have $A_{k,m,n}(x)\in\FF$. It follows that $\Omega^\star=\{\omega\in\Omega''\colon\text{ ergodic}\}\in\FF$ as well.

To finish the proof, note that $\ell_\omega(f\circ\tau_x)=\ell_\omega(f)$ gives $E_{\BbbP_{\tau_x(\omega)}}(f)=E_{\BbbP_\omega}(f)$. This in turn shows $\BbbP_{\tau_x(\omega)}=\BbbP_{\omega}$ for all~$x\in\Z^d$ and thus implies that~$\omega\in\Omega^\star$ if and only if $\tau_x(\omega)\in\Omega^\star$ for all~$x\in\Z^d$. It follows that~$\Omega^\star$ is translation invariant. Given an ergodic probability measure~$\nu$ on~$(\Omega,\FF)$, Birkhoff's Ergodic Theorem and separability of~$\Cloc(\R^d)$ ensure that, for $\nu$-a.e.~$\omega'\in\Omega$, the double limit in the definition of~$\Omega'$ vanishes for~$\omega'$ and, moreover, $\ell_{\omega'}(f)$ exists and equals~$E_\nu(f)$ for all~$f\in\Cloc(\R^d)$. Each such~$\omega'$ thus belongs to~$\Omega''$ and obeys~$\BbbP_{\omega'}=\nu$, and so~$\omega'\in\Omega^\star$. It follows that~$\nu(\Omega^\star)=1$.
\end{proofsect}

We remark that, while ``large'' from the perspective of any ergodic measure, the set~$\Omega^\star$ is nowhere dense and is thus a ``small'' subset of~$\Omega$ in the sense of topology.

\subsection{Heat-kernel bounds and  the  Conversion lemma}
\noindent
We now return to the random walk. As discussed  earlier,  an important technical input for our derivations are heat-kernel estimates.  These  are typically Gaussian-type bounds~on
\begin{equation}
%\label{}
y\mapsto\frac{P^x_\omega(X_n=y)}{\pi_\omega(y)}.
\end{equation}
One is generally interested in both upper and lower bounds (see, e.g., the monographs by Kumagai~\cite{Kumagai} or Barlow~\cite{B17}) but we will only need upper bounds in this paper.

Heat-kernel estimates for random walks on weighted graphs are the subject of an early work of Delmotte~\cite{D99}  which   would be sufficient for dealing with  uniformly elliptic conductances. The more general setting of our result requires a non-trivial extension due to Andres, Deuschel and Slowik~\cite{ADS16}. Both studies are based on advanced techniques from elliptic-regularity theory; namely, the Moser iteration.

A minor technical hurdle is that the conclusions of~\cite{ADS16} are cast for the continuous time version~$Y$ of the (discrete-time) process~$X$, both constant and variable speed. We will only need the former which is defined as follows. Assume that~$X$ is coupled on the same probability space with an independent sequence of i.i.d.\ exponentials $T_1,T_2,\dots$ with parameter one. Given~$t\ge0$, set
\begin{equation}
\label{E:2.22}
Y_t:=X_{ N(t)}\quad\text{for}\quad  N(t):=\sup\{n\ge0\colon T_1+\dots+T_n\le t\}.
\end{equation}
We will continue writing~$P^x_\omega$ for the joint law of these processes subject to the initial condition $P^x_\omega(Y_0=X_0=x)=1$.
Theorem~1.6 of~\cite{ADS16} then gives:

\begin{proposition}
\label{prop-ADS}
Let $p,q>1$ obey \eqref{E:2.6}. Then for each~$\omega\in\Omega_{p,q}$ there exist quantities $c_1,c_2,c_3,c_4\in(0,\infty)$ such that for all $x,y\in\Z^d$ and all $t\ge0$ with $\sqrt t\ge c_3(1+\min\{|x|_1,|y|_1\})$:
%\settowidth{\leftmargini}{(11)}
\begin{enumerate}
\item[(1)] if $|x-y|_1\le c_4 t$, then
\begin{equation}
\label{E:3.13r}
P^x_\omega(Y_t=y)\le c_1 t^{-d/2}\exp\Bigl\{-c_2\tfrac{|x-y|_1^2}{t}\Bigr\}\pi_\omega(y),
\end{equation}
\item[(2)] while if $|x-y|_1> c_4 t$, then
\begin{equation}
\label{E:3.14r}
P^x_\omega(Y_t=y)\le c_1 t^{-d/2}\exp\biggl\{-c_2|x-y|_1\Bigl(1+\log\bigl(1+\tfrac{|x-y|_1}{\sqrt t}\bigr)\Bigr)\biggr\}\pi_\omega(y).
\end{equation}
\end{enumerate}
Here~$|\cdot|_1$ is the $\ell^1$-norm on~$\R^d$.
\end{proposition}

\begin{proofsect}{Proof}
We need to verify the conditions of Theorem~1.6 of~\cite{ADS16}. (We adhere to the notation of~\cite{ADS16} in this proof noting, in particular, that~$N(\cdot)$ means something very different there than in \eqref{E:2.22}.) First, Remark~1.2 of~\cite{ADS16} tells us that $(\Z^d,E(\Z^d))$ is included among admissible graphs and that the quantities $d'$ and~$N_1(x)$ that appear in \cite[Theorem~1.6]{ADS16} take values $d'=d$ and~$N_1(x)=1$. The containment $\omega\in\Omega_{p,q}$ and the regularity of~$\Z^d$ in turn show that the quantities $\bar\mu_p(x)$ and~$\bar\nu_q(x)$ right before \cite[Assumption 1.5]{ADS16} are finite and independent of~$x$. 

For Assumption~1.5 of~\cite{ADS16}, we first find $N(0)$ such that the inequalities in the last line thereof hold for $x=0$ with, say, $3/2$ instead of~$2$ on the right-hand side. To get the inequalities as stated for~$x\ne0$, it then suffices to take $N(x)\ge N(0)$ so that $|\Lambda_{n+|x|_1}|/|\Lambda_n|\le 4/3$ for all~$n\ge N(x)$. This works for $N(x):=N(0)+c|x|_1$ once~$c$ is large enough. Theorem~1.6 of~\cite{ADS16} yields the claim modulo redefinition of constants. The fact that we can write $\min\{|x|_1,|y|_1\}$ in the condition on~$t$ arises from the fact that, by reversibility, $P^x_\omega(Y_t=y)/\pi_\omega(y)$ is symmetric upon exchange of~$x$ and~$y$.
\end{proofsect}

The heat kernel estimates enter our proofs in the conversion of time averages (and those with respect to the law of~$X$) to averages in the physical space:

\begin{lemma}[Conversion lemma]
\label{lemma-3.3a}
Let $p,q>1$ obey \eqref{E:2.6}. For each~$\omega\in\Omega_{p,q}$ there exists $c(\omega)\in(0,\infty)$ such that for all bounded~$f\colon\Z^d\to[0,\infty)$,
\begin{equation}
%\label{}
\limsup_{n\to\infty}\frac1n\sum_{k=0}^{n-1}E_\omega^0\bigl(f(X_k)\bigr)\le c(\omega)\,\limsup_{n\to\infty}\frac1{|\Lambda_n|}\sum_{x\in\Lambda_n}f(x)\pi_\omega(x).
\end{equation}
\end{lemma}

\begin{proofsect}{Proof}
Fix~$\delta\in(0,1/2)$ and recall the definition of $t\mapsto N(t)$ from \eqref{E:2.22}. Since~$N(t)$ is Poisson with parameter~$t$ and $\int_0^\infty P(N(t)=k)\textd t=1$ for all~$k\ge0$, routine estimates show that there is~$a>0$ such that
\begin{equation}
%\label{}
\min_{\lfloor 2\delta n\rfloor\le k\le n}\int_{\delta n}^{2n} P\bigl(N(t)=k\bigr) \textd t \ge1-\texte^{-a n}
\end{equation}
holds for all~$n$ large. Once $\texte^{-a n}\le\frac12$, the construction of~$Y$ gives
\begin{equation}
%\label{}
\int_{\delta n}^{2n}P_\omega^0(Y_t = x) \textd t \ge \frac12\sum_{k=\lfloor 2\delta n\rfloor}^{n} P_\omega^0(X_k = x).
\end{equation}
Hence we get
\begin{equation}
\label{E:3.17s}
\sum_{k=0}^{n-1}E_\omega^0\bigl(f(X_k)\bigr)\le 2\delta n(\sup|f|)+2\int_{\delta n}^{2n}\biggl(\,\sum_{x\in\Z^d} f(x)P_\omega^0(Y_t = x) \biggr)\textd t
\end{equation}
which we will now estimate using the heat kernel bound from Proposition~\ref{prop-ADS}.

For integer $r\ge0$ and real~$t>1$, let~$A_{r,t}:=\{x\in\Z^d\colon r\sqrt t\le|x|_1<(r+1)\sqrt t\}$.  Invoking \twoeqref{E:3.13r}{E:3.14r} with the logarithmic term dropped in \eqref{E:3.14r} and noting that $\min\{|x|_1,|x|_1^2/t\}$ exceeds~$r$ on~$A_{r,t}$ once~$t\ge1$,  the integrand on the right of \eqref{E:3.17s} is at most
\begin{equation}
%\label{}
\sum_{r\ge0}c_1 \texte^{-c_2r} \frac1{t^{d/2}}\sum_{x\in A_{r,t}}f(x)\pi_\omega(x),
\end{equation}
provided that also $t\ge c_3$ (with all three $c_i$'s depending on~$\omega$). Noting that $A_{r,t}\subseteq \Lambda_{(r+1)\sqrt t}$, this is in turn no larger than
\begin{equation}
%\label{}
\Bigl(\sum_{r\ge0}c_1 (r+1)^d\texte^{-c_2r}\Bigr)\sup_{m\ge \lfloor\sqrt{t}\rfloor} \,\frac1{m^d}\sum_{x\in \Lambda_m}f(x)\pi_\omega(x).
\end{equation}
The integral in \eqref{E:3.17s} is thus at most~$2n$ times this quantity at~$t:=\delta n$. Dividing \eqref{E:3.17s} by~$n$ and noting that $|\Lambda_m|/m^d\le 3^d$, the claim follows by taking~$n\to\infty$ and~$\delta\downarrow0$.
\end{proofsect}

\subsection{Exit time estimate and tightness}
\noindent
Another place where input from heat-kernel technology is useful is tightness of the processes $t\mapsto B_t^{(n)}$, where $B_t^{(n)}$ is as in \eqref{E:1.3}. 
%By definition, this process has continuous sample paths and is thus a random variable taking values in~$C([0,\infty),\R^d)$. 
This comes through a bound on the lower tail of the first exit time of the walk from a domain. For~$\Lambda\subseteq\Z^d$, set
\begin{equation}
%\label{}
\ttau_{\Lambda}:=\inf\{k\ge0\colon X_k\not\in\Lambda\}.
\end{equation}
We then have:

\begin{proposition}
\label{prop-4.1}
Let $d\ge2$ and let $p,q>d/2$ obey \eqref{E:2.6}. For all~$\alpha>0$, $\sigma>0$ and~$\omega\in\Omega_{p,q}$ there exists $R_0\ge1$ and $\tilde c_1=\tilde c_1(\alpha,\sigma,\omega)\in[1,\infty)$ such that for all $R\ge R_0$ and all $t\ge0$,
\begin{equation}
\label{E:2.32}
\max_{x\in\Lambda_R} P^x_\omega\bigl(\ttau_{\Lambda_{\sigma R}(x)}\le t\bigr)\le \tilde c_1\Bigl(\frac t{R^2}\Bigr)^\alpha,
\end{equation}
where we set $\Lambda_r(x):=x+\Lambda_r$.
\end{proposition}

While this statement should in principle be obtainable from the heat-kernel bounds in Proposition~\ref{prop-ADS}, this does not seem to apply directly due to the restriction on~$t$ relative to~$x$ and~$y$. We thus employ an argument from Biskup, Chen, Kumagai and Wang~\cite{BCKW} that controls the exit time (and the heat kernel) using a different method than \cite{ADS16}. The above statement requires only minor adaptations of the proofs in \cite{BCKW} and so we relegate its proof to Section~\ref{sec7}. 

\begin{remark}
Note that Proposition~\ref{prop-4.1} makes no claim about~$d=1$. This is because its proof (just as the arguments in~\cite{BCKW}) is based on Sobolev inequalities that behave somewhat differently in dimension one. This is no loss because our proof of Theorem~\ref{thm-2.4} will need the exit time estimate only in~$d\ge2$. 
\end{remark}

\smallskip
We will now show how Proposition~\ref{prop-4.1} implies tightness. As in~\cite[Proposition~4.1]{BCKW}, we could proceed by invoking a ``stopping time'' criterion of Aldous~\cite[Theorem~1]{Aldous}. A drawback of this approach is that the criterion is cast in the Skorohod space and topology, while our setting naturally rests in the Wiener space~$C([0,\infty))$. Since a careful choice of various parameters necessitates that we present a full argument anyway, we will prove tightness directly in~$C([0,\infty))$.
  
Given a real~$T>0$ and a function $f\colon[0,T]\to\R^d$, recall that its $\ell^\infty$-oscillation over intervals of size~$\delta$ is defined by
\begin{equation}
%\label{}
\text{osc}_f\bigl([0,T],\delta\bigr):=\sup_{\begin{subarray}{c}
s,t\in[0,T]\\|t-s|<\delta
\end{subarray}}\,
\bigl|f(t)-f(s)\bigr|_\infty,
\end{equation}
where $|\cdot|_\infty$ is the $\ell^\infty$-norm on~$\R^d$.
The conclusion of Proposition~\ref{prop-4.1} then  bounds the oscillation of~$B^{(n)}$ from \eqref{E:1.3} as follows:

\begin{lemma}
\label{lemma-4.2}
Let $d\ge2$ and $p,q>d/2$ obey \eqref{E:2.6}. For all~$\epsilon>0$, $T\ge1$ and~$\omega\in\Omega_{p,q}$ there exists $\hat c_1=\hat c_1(\epsilon,T,\omega)\in(0,\infty)$ such that for all~$\delta>0$ and all~$n> 1/\delta$,
\begin{equation}
\label{E:2.34}
P^0_\omega\Bigl(\text{\rm osc}_{B^{(n)}}\bigl([0,T],\delta\bigr)>\epsilon\,\wedge\,\ttau_{\Lambda_{\sqrt{n}/\epsilon}}>3Tn\Bigr)
\le\hat c_1\sqrt\delta.
\end{equation}
\end{lemma}

\begin{proofsect}{Proof}
Fix~$T\ge1$ and,  assuming $\hat c_1\ge1$, pick $\delta\in(0,1)$ and $n\ge1$ such that  $n\delta>1$. For $j\ge1$ set $k_j:=\lfloor jn\delta\rfloor$ and note that then $k_{j+2}-k_{j-1}\ge 3\delta n-1>\delta n+2$. Given~$t,s\in[0,T]$ with $0<t-s<\delta$, we then have $\lfloor tn\rfloor-\lfloor sn\rfloor\le tn-sn+1<\delta n+1$ and so there exists~$j$ with~$1\le j\le T/\delta$ such that $\lfloor tn\rfloor$, $\lfloor tn\rfloor+1$ and~$\lfloor sn\rfloor$ all lie between~$k_{j-1}$ and~$k_{j+2}$. Examining the various subintervals of $[k_{j-1},k_{j+2}]$ that $\lfloor tn\rfloor$, $\lfloor tn\rfloor+1$ and~$\lfloor sn\rfloor$  may fall into, we get
\begin{equation}
\label{E:2.35}
\bigl\{\text{osc}_{B^{(n)}}\bigl([0,T],\delta\bigr)>\epsilon\bigr\}\subseteq
\bigcup_{1\le j\le 3T/\delta}\Bigl\{\,\,\max_{k_{j-1}\le\ell\le k_j}|X_\ell-X_{k_{j-1}}|_\infty>\tfrac13\epsilon\sqrt n\Bigr\},
\end{equation}
where we also noted that $T/\delta+2\le 3T/\delta$.

Next we note that the inequality corresponding to index~$j$ on the right of \eqref{E:2.35} implies that $\ttau_{\Lambda_{\epsilon\sqrt n/4}(x)}\le k_j-k_{j-1}$ occurs for the portion of the path started at~$X_{k_{j-1}}$. Given any~$R\ge1$, the Markov property along with the union bound therefore give
\begin{equation}
%\label{}
P^0_\omega\Bigl(\text{osc}_{B^{(n)}}\bigl([0,T],\delta\bigr)>\epsilon\,\wedge\ttau_{\Lambda_R}>3Tn\Bigr)
\le 3T\delta^{-1} \max_{x\in\Lambda_R} P^x_\omega\bigl(\ttau_{\Lambda_{\epsilon\sqrt n/3}(x)}
\le 2\delta n\bigr),
\end{equation}
where we also used that $k_j-k_{j-1}\le2\delta n$. Now set $R:=\sqrt{n}/\epsilon$ and apply Proposition~\ref{prop-4.1} with $\alpha:=3/2$,~$\sigma:=\epsilon^2/3$  and~$t:=2\delta n$  to get
\begin{equation}
%\label{}
\max_{x\in\Lambda_{\sqrt{n}/\epsilon}} P^x_\omega\bigl(\ttau_{\Lambda_{\epsilon\sqrt n/3}(x)}
\le 2\delta n\bigr)\le  2^{3/2}\tilde c_1 \epsilon^3 \delta^{3/2}
\end{equation}
 as soon as  $R\ge R_0$. Since $n>1/\delta$,~$\hat c_1\ge R_0\epsilon$ and $\sqrt n/\epsilon\le R_0$ altogether imply that $\hat c_1\sqrt\delta\ge1$, the bound \eqref{E:2.34} holds trivially for~$R\le R_0$ if $\hat c_1\ge R_0\epsilon$. The claim thus follows by setting $\hat c_1:=\max\{3\cdot 2^{3/2} T\epsilon^3\tilde c_1,R_0\epsilon,1\}$.
\end{proofsect}

With this in hand, we are ready to state:

\begin{proposition}
\label{prop-2.7}
Let $d\ge2$ and $p,q>d/2$ obey \eqref{E:2.6}. For~$n\ge1$ and~$\omega\in\Omega$, let~$P^{(n,\omega)}$ be the law of $B^{(n)}$ from \eqref{E:1.3} induced by~$P_\omega^0$ on $(C([0,\infty),\R^d),\BB(C([0,\infty),\R^d)))$. Then
\begin{equation}
%\label{}
\forall\omega\in\Omega_{p,q}\colon\,\, \{P^{(n,\omega)}\colon n\ge1\}\text{\rm\ is tight.}
\end{equation}
\end{proposition}

\begin{proofsect}{Proof}
Given~$\delta\in(0,1)$, $\epsilon>0$ and~$T\ge1$, let
\begin{equation}
%\label{}
\KK_{\delta,\epsilon,T}:=\Bigl\{f\in C\bigl([0,\infty)\bigr)\colon f(0)=0\,\wedge\,\text{osc}_f\bigl([0,T],\delta\bigr)\le\epsilon\Bigr\}.
\end{equation}
Proposition~\ref{prop-4.1} with~$\alpha:=1$, $\sigma:=1$, $t:=3Tn$ and $R:=\sqrt n/\epsilon$ gives
\begin{equation}
\label{E:2.40}
P^0_\omega\bigl(\ttau_{\Lambda_{\sqrt n/\epsilon}}\le 3Tn\bigr)\le \tilde c_1(1,1,\omega)3T\epsilon^2.
\end{equation}
The piece-wise linear nature of~$B^{(n)}$ and the fact that~$X$ jumps by a unit $\ell^\infty$-distance in each step imply $\text{\rm osc}_{B^{(n)}}([0,T],\delta)\le \delta\sqrt n\le\sqrt\delta$ when~$\delta\le1/n$. Once $\sqrt\delta<\epsilon$, the restriction $n>1/\delta$ in Lemma~\ref{lemma-4.2} is thus moot and \eqref{E:2.34} along with~\eqref{E:2.40} show
\begin{equation}
%\label{}
P^{(n,\omega)}(K_{\delta,\epsilon,T})\ge 1-\tilde c_1(1,1,\omega)3T\epsilon^2-\hat c_1(\epsilon,T,\omega)\sqrt\delta
\end{equation}
for all~$n\ge1$ and all~$\omega\in\Omega_{p,q}$. Now take $T_k\uparrow\infty$ and pick $\epsilon_k\downarrow0$ and~$\delta_k\downarrow0$ so that
\begin{equation}
%\label{}
\sqrt{\delta_k}<\epsilon_k\,\wedge\, \tilde c_1(1,1,\omega)3T_k\epsilon_k^2<2^{-k}\,\wedge\hat c_1(\epsilon_k,T_k,\omega)\sqrt{\delta_k}<2^{-k}
\end{equation}
for each~$k\ge1$. For the set
\begin{equation}
%\label{}
\KK_\ell:=\bigcap_{k\ge\ell}\KK_{\delta_k,\epsilon_k,T_k}
\end{equation}
the union bound shows
\begin{equation}
%\label{}
\forall n\ge1\,\forall\omega\in\Omega_{p,q}\colon\,\,P^{(n,\omega)}(\KK_\ell)\ge 1-2^{-\ell+2}.
\end{equation}
We finish by noting that, by the  Arzel\`a-Ascoli Theorem or a direct argument, $\KK_\ell$ is compact in~$C([0,\infty))$ relative to the topology of locally uniform convergence.
\end{proofsect}

We will use the above observations in two ways. First, we will continue using the exit-time estimate to effectively confine the walk to a fixed box on a diffusive scale. Second, the tightness reduces the proof of an IIP to convergence of finite-dimensional distributions for which it suffices, more or less, to prove just a Central Limit Theorem.

%\newpage
\section{Ergodicity of Markov chain on environments}
\label{sec4}\noindent
This section addresses the first important step in our proof of an IIP for deterministic conductance configurations; namely, averaging at large temporal scales of the  walk. We then use it to prove our main theorem in spatial dimension one.

\subsection{Averaging for the random walk}
A starting point of many approaches to random walks in random environments is the so called ``point of view of the particle.'' The phrase refers to an observer traveling with the random walk~$X$ who then sees environment $\tau_{X_n}(\omega)$ at time~$n$.
For~$X$ sampled from~$P_\omega^0$, the sequence $\{\tau_{X_n}(\omega)\colon n\ge0\}$ of environments is itself a Markov chain but this time on state space~$\Omega$. An important observation is that many properties of the physical chain  and, in particular,  the value~$X_n$ itself can be encoded as additive functionals of the environment chain.

The Markovian dynamics of the chain from the ``point of view of the particle '' is described either by its generator~$\LL$ that acts on $f\in \Cb(\Omega)$ as
\begin{equation}
\label{E:4.1w}
(\LL f)(\omega):=\frac1{\pi_\omega(0)}\sum_{\begin{subarray}{c}
x\in\Z^d\\(0,x)\in E(\Z^d)
\end{subarray}}
\cc_\omega(0,x)\bigl[f\circ\tau_x(\omega)-f(\omega)\bigr]
\end{equation}
or directly by its transition kernel~$\Pi$ that is obtained from~$\LL$ via
\begin{equation}
\label{E:3.41q}
\Pi f:=f+\LL f.
\end{equation}
As is easy to check, $\Pi f$ is continuous on~$\Omega$ if~$f$ is and
\begin{equation}
\label{E:4.3i}
\sup_{\omega\in\Omega}\,\bigl|(\Pi f)(\omega)\bigr|\le\sup_{\omega\in\Omega}\,\bigl|f(\omega)\bigr|.
\end{equation}
Consequently,~$\Pi$ is a bounded linear operator on~$\Cb(\Omega)$ with norm one.

A key benefit that the ``point of the view of the particle'' brings to the overall theory is a way to demonstrate averaging at large temporal scales of the walk. For environments sampled from an ergodic law, this follows from the Ergodic Theorem for Markov chains and the fact that, by reversibility, a measure closely related to the \emph{a priori} law, cf \eqref{E:3.5t} below, is invariant and ergodic for the Markovian dynamics (see, e.g., \cite[Proposition~2.3]{B11}). Our main observation here is that the same can be obtained while relying solely on block averages. The technical tool that makes this work is the Conversion Lemma (Lemma~\ref{lemma-3.3a}) enabled by the heat-kernel bounds (Proposition~\ref{prop-ADS}). 

\begin{theorem}
\label{thm-2.1}
Let $p,q>1$ obey \eqref{E:2.6}. Then for all~$\omega\in\Omega^\star\cap\Omega_{p,q}$ and all $f\in \Cb(\Omega)$,
\begin{equation}
\label{E:4.4i}
\frac1n\sum_{k=0}^{n-1}f\circ\tau_{X_k}(\omega)\,\underset{n\to\infty}\longrightarrow\,E_{\Q_\omega}(f)
\qquad\text{\rm in }L^1(P_\omega^0),
\end{equation}
where
\begin{equation}
\label{E:3.5t}
\Q_\omega(\textd\omega'):=\frac{\pi_{\omega'}(0)}{E_{\BbbP_\omega}\pi(0)}\BbbP_\omega(\textd\omega')
\end{equation}
for~$\BbbP_\omega$ related to~$\omega$ as in Proposition~\ref{prop-2.2}.
\end{theorem}

\begin{proofsect}{Proof}
Thanks to Lemma~\ref{lemma-3.1}, the assumed conditions on~$p$, $q$ and~$\omega$ ensure $\pi(0)\in L^1(\BbbP_\omega)$ and so~$\Q_\omega$ is well defined. We thus need to show \eqref{E:4.4i}.
Pick~$f\in \Cb(\Omega)$ and,  reducing~$f$ by its expectation under $\Q_\omega$,  assume
\begin{equation}
\label{E:2.4}
E_{\Q_\omega}(f)=0.
\end{equation}
The proof is based on martingale approximation. Let $\epsilon>0$ and, noting that $\sup|\Pi^n f|\le\sup|f|$ by \eqref{E:4.3i}, define $h_\epsilon\colon\Omega\to\R$ by
\begin{equation}
\label{E:3.7p}
h_\epsilon(\omega'):=\sum_{n\ge0}\frac{1}{(1+\epsilon)^{n+1}}(\Pi^n f)(\omega').
\end{equation}
As is readily checked, $h_\epsilon$ solves the ``massive'' (a.k.a.~screened) Poisson equation
\begin{equation}
\label{E:2.5}
(\epsilon-\LL)h_\epsilon(\omega')=f(\omega').
\end{equation}
The uniform boundedness of $n\mapsto\Pi^n f$ applied under \eqref{E:3.7p} gives~$h_\epsilon\in C(\Omega)$. We also get that $\epsilon h_\epsilon$ is bounded uniformly in~$\epsilon>0$. In particular, $h_\epsilon\in\Cb(\Omega)$ for all $\epsilon>0$.

We now use \eqref{E:2.5} to rewrite the sum in the statement as
\begin{equation}
\label{E:2.7}
\begin{aligned}
\sum_{k=0}^{n-1}f\circ\tau_{X_k}(\omega) 
&=\sum_{k=0}^{n-1}(\epsilon+1-\Pi)h_\epsilon\circ\tau_{X_k}(\omega)
\\
&=\epsilon\sum_{k=0}^{n-1}h_\epsilon\circ\tau_{X_k}(\omega)+\sum_{k=0}^{n-1}\bigl[h_\epsilon\circ\tau_{X_k}(\omega)-\Pi h_\epsilon\circ\tau_{X_k}(\omega)\bigr].
\end{aligned}
\end{equation}
Let $\FF_k:=\sigma(X_0,\dots,X_k)$. Noting that
\begin{equation}
%\label{}
\Pi h_\epsilon\circ\tau_{X_k}(\omega)=E^0_\omega\bigl(h_\epsilon\circ\tau_{X_{k+1}}(\omega)\big|\FF_k\bigr),
\end{equation}
the sum on the right of \eqref{E:2.7} can be shuffled to the form
\begin{equation}
\label{E:2.9}
\begin{aligned}
\sum_{k=0}^{n-1}\bigl[h_\epsilon\circ\tau_{X_k}(\omega)-\Pi &h_\epsilon\circ\tau_{X_k}(\omega)\bigr]
\\&=h_\epsilon(\omega)- h_\epsilon\circ\tau_{X_n}(\omega) 
\\&\qquad+\sum_{k=1}^n\bigl[h_\epsilon\circ\tau_{X_k}(\omega)-E\bigl(h_\epsilon\circ\tau_{X_k}(\omega)\big|\FF_{k-1}\bigr)\bigr].
\end{aligned}
\end{equation}
As $h_\epsilon$ is bounded, the first two terms on the right are bounded. The sum is in turn a martingale with uniformly bounded increments. By Doob's $L^2$-inequality (or Azuma's), it will thus typically be order $\sqrt{n}$. It follows that, upon division by~$n$ and taking~$n\to\infty$, the  quantity  \eqref{E:2.9} vanishes in $L^1(P^0_\omega)$.

It thus suffices to show that, under \eqref{E:2.4},
\begin{equation}
\label{E:2.10}
\frac1n\sum_{k=0}^{n-1}\epsilon h_\epsilon\circ\tau_{X_k}(\omega)\,\underset{\begin{subarray}{c}
n\to\infty\\\epsilon\downarrow0
\end{subarray}}
\longrightarrow\,0
\qquad\text{in }L^1(P_\omega^0).
\end{equation}
Here the Conversion Lemma (Lemma~\ref{lemma-3.3a}) enabled by the heat-kernel estimates in Proposition~\ref{prop-ADS} along with the fact that \eqref{E:3.11i} holds for $h:=\pi(0)$ with~$\alpha:=p$ gives
\begin{equation}
%\label{}
\begin{aligned}
\limsup_{n\to\infty}\,\frac1n\sum_{k=0}^{n-1}\bigl|\epsilon h_\epsilon\circ\tau_{X_k}(\omega)\bigr|
&\le c(\omega)\limsup_{n\to\infty}\frac1{|\Lambda_n|}\sum_{x\in\Lambda_n}\bigl|\epsilon h_\epsilon\circ\tau_{x}(\omega)\bigr|\pi_\omega(x)
\\
&= c(\omega) E_{\BbbP_\omega}\bigl(\pi(0)\epsilon |h_\epsilon|\bigr)
\le c(\omega)E_{\BbbP_\omega}\bigl(\pi(0)\bigr)\Vert \epsilon h_\epsilon\Vert_{L^2(\Q_\omega)},
\end{aligned}
\end{equation}
where we also used the definition of~$\Q_\omega$ and the Cauchy-Schwarz inequality.  To prove \eqref{E:2.10} it suffices  to show that $\epsilon h_\epsilon\to0$ in $L^2(\Q_\omega)$  as $\epsilon\downarrow0$,  which we achieve by way of  spectral calculus.

Write $\langle\cdot,\cdot\rangle$ for the inner product in~$L^2(\Q_\omega)$ and note that $-\LL$ is a bounded symmetric positive semi-definite operator on~$L^2(\Q_\omega)$. As $f\in L^2(\Q_\omega)$, the Spectral Theorem asserts the existence of a finite Borel measure~$\mu_f$ on~$\R$ such that $\langle f,\phi(-\LL)f\rangle=\int \phi(\lambda)\mu_f(\textd \lambda)$ holds for any $\phi\in L^1(\mu_f)$. In conjunction with \eqref{E:2.5}, this yields
\begin{equation}
\label{E:2.15}
\begin{aligned}
\Vert \epsilon h_\epsilon\Vert_{L^2(\Q_\omega)}^2 &= \langle\epsilon h_\epsilon,\epsilon h_\epsilon\rangle 
\\
&= \bigl\langle\epsilon(\epsilon-\LL)^{-1}f,\epsilon(\epsilon-\LL)^{-1}f\bigr\rangle=\int_{[0,2]} \Bigl(\frac{\epsilon}{\epsilon+\lambda}\Bigr)^2\mu_f(\textd\lambda),
\end{aligned}
\end{equation}
where we also noted that, since~$\Pi$ is a probability kernel and thus a contraction in~$L^2$, the spectrum of~$-\LL$ and so also the support of~$\mu_f$ are confined to~$[0,2]$. 

The function under  the  integral in \eqref{E:2.15} is bounded by one and tends to $1_{\{0\}}$ as $\epsilon\downarrow0$. The Bounded Convergence Theorem shows
\begin{equation}
%\label{}
\Vert \epsilon h_\epsilon\Vert_{L^2(\Q_\omega)}^2\,\underset{\epsilon\downarrow0}\longrightarrow\,\mu_f(\{0\}).
\end{equation}
In light of $-\LL\ge0$, the right-hand side is the norm-squared of the orthogonal projection of~$f$ on $\text{Ker}(\LL)$. Here we note that $g\in\text{Ker}(\LL)$ reads $g = \Pi g$ which  upon iteration gives $g(\omega')=E_{\omega'}^0(g\circ\tau_{X_k}(\omega'))$ for~$\Q_\omega$-a.e.~$\omega'\in\Omega$ and all~$k\ge0$.  Noting that ergodicity of~$\BbbP_\omega$ under spatial shifts induces ergodicity of~$\Q_\omega$ under this chain (see again~\cite[Proposition~2.3]{B11}),  the Birkhoff Ergodic  Theorem gives $g=E_{\Q_\omega}(g)$ a.s, showing that~$g$ is constant. But \eqref{E:2.4} means that~$f$ is orthogonal to constants, and so the projection of~$f$ on $\text{Ker}(\LL)$ vanishes. This gives
\begin{equation}
%\label{}
\mu_f(\{0\})=0
\end{equation}
and thus completes the proof of \eqref{E:2.10} and the whole theorem.
\end{proofsect}

We  remark  that the proof that $\epsilon h_\epsilon\to0$ in~$L^2(\Q_\omega)$ is already contained in Kipnis and Varadhan~\cite{KV86} albeit under the assumption that~$f\in\text{Dom}((-\LL)^{-1/2})$ which,  being stronger than \eqref{E:2.4}, gives that  $\epsilon^{1/2} h_\epsilon\to0$. We included the above argument to demonstrate the role of the ergodicity of~$\BbbP_\omega$ that becomes apparent only in the final step. 

\begin{comment}
Since the control of errors is based largely on a martingale argument, we in fact get even the following statement:

\begin{corollary}
\label{cor-2.2}
Let~$\omega\in\Omega^\star$ and let $\{\tau_n\}_{n=1}^\infty$ be a sequence of stopping times for the random walk~$X$. Then for all $f\in C(\Omega)$,
\begin{equation}
%\label{}
\lim_{n\to\infty}\,\frac1n\biggl|\,\sum_{k=1}^{n\wedge\tau_n}f\circ\tau_{X_{k-1}}(\omega)-(n\wedge\tau_n)E_{\Q_\omega}(f)\biggr|=0
\end{equation}

\end{corollary}

\end{comment}

\subsection{Zero speed and an IIP in dimension one}
\label{sec-3.2}\noindent
A direct (albeit well known) consequence of Theorem~\ref{thm-2.1} is the absence of speed of the  random walk. We recall its short proof in:

\begin{corollary}
\label{cor-3.2}
Let $p,q>1$ obey \eqref{E:2.6}. Then for all~$\omega\in\Omega^\star\cap\Omega_{p,q}$,
\begin{equation}
%\label{}
\frac{X_n}n\,\underset{n\to\infty}{\longrightarrow}\,0\quad\text{\rm in $P_\omega^0$-probability}.
\end{equation}
\end{corollary}

\begin{proofsect}{Proof}
We start with the standard rewrite
\begin{equation}
\label{E:3.2u}
\begin{aligned}
X_n &= X_0+\sum_{k=1}^n (X_k-X_{k-1})
\\
&=X_0+\sum_{k=0}^{n-1}V\circ\tau_{X_k}(\omega)
+\sum_{k=1}^n\Bigl[(X_k-X_{k-1})-E_\omega^0\big(X_k-X_{k-1}\,\big|\,\FF_{k-1}\bigl)\Bigr],
\end{aligned}
\end{equation}
where  $\FF_k:=\sigma(X_0,\dots,X_k)$ and 
\begin{equation}
\label{E:3.2r}
V(\omega):=E_\omega^0(X_1)
\end{equation}
is the local drift that the walk feels at the origin. 
Now observe that the second sum on the right of \eqref{E:3.2u} is a martingale with bounded increments which, by Doob's $L^2$-inequality, is at most order $\sqrt{n}$ in probability. Noting that  $V\in\Cb(\Omega)$, the first sum divided by~$n$ in turn converges to~$E_{\Q_\omega}(V)$ in~$P^0_\omega$-probability by  Theorem~\ref{thm-2.1}. A standard calculation shows $E_{\Q_\omega}(V)=0$ and so the claim follows.
\end{proofsect}

Theorem~\ref{thm-2.1} is also sufficient to give the proof of our main result in spatial dimension one. For this we first need to extend the observation from Lemma~\ref{lemma-3.1} to averages with respect to the random walk:

\begin{lemma}
\label{lemma-4.2w}
Let $p,q>1$ obey \eqref{E:2.6}.
Suppose~$f\in C(\Omega)$ and ~$h\colon\Omega\to[0,\infty)$ are such that \eqref{E:3.10i} holds. Assume~$\omega\in\Omega^\star\cap\Omega_{p,q}$ is such that, for some~$\alpha>1$,
\begin{equation}
\label{E:4.18i}
\sup_{ n\ge1}\frac1{n^d}\sum_{x\in\Lambda_n}[h\circ\tau_x(\omega)]^\alpha\pi_\omega(x)<\infty.
\end{equation}
Then $f\in L^1(\Q_\omega)$ and \eqref{E:4.4i} holds  for~$f$. 
\end{lemma}

\begin{proofsect}{Proof}
We proceed very much like in the proof of Lemma~\ref{lemma-3.1}. Define~$f_r$ by \eqref{E:3.12i} and~$g_r$ by \eqref{E:3.14}. Then \eqref{E:3.10i} and Jensen's inequality give
\begin{equation}
\label{E:4.18ii}
\frac1n\sum_{k=0}^{n-1}E_\omega^0\bigl(|f_r|\circ\tau_{X_k}(\omega)\bigr)
\le \biggl(\,\frac1n\sum_{k=0}^{n-1}E_\omega^0\bigl([h\circ\tau_{X_k}(\omega)]^\alpha\bigr)\biggr)^{1/\alpha}.
\end{equation}
The Conversion Lemma (Lemma~\ref{lemma-3.3a}) then bounds the $n\to\infty$ limit of the right-hand side in terms of the quantity in \eqref{E:4.18i}. As $f_r\in\Cb(\Omega)$, Theorem~\ref{thm-2.1} shows that the left-hand side tends to $E_{\Q_\omega}(|f_r|)$. The Monotone Convergence Theorem then allows us to take~$r\to\infty$ and conclude $f\in L^1(\Q_\omega)$.

For the convergence \eqref{E:4.4i} for~$f$, we invoke an analogue of \eqref{E:3.15i} to get
\begin{equation}
\begin{aligned}
E_\omega^0\biggl|\frac1n\sum_{k=0}^{n-1}&f\circ\tau_{X_k}(\omega)-\frac1n\sum_{k=0}^{n-1}f_r\circ\tau_{X_k}(\omega)\biggr|
\\
&\le\biggl(\,\frac1n\sum_{k=0}^{n-1}E_\omega^0\bigl([h\circ\tau_{X_k}(\omega)]^\alpha\bigr)\biggr)^{1/\alpha}\biggl(\,\frac1n\sum_{k=0}^{n-1}E_\omega^0\bigl([g_r\circ\tau_{X_k}(\omega)]^{\alpha'}\bigr)\biggr)^{1/\alpha'}.
\end{aligned}
\end{equation}
The \emph{limes superior} as~$n\to\infty$ of the first term on the right is finite by the argument after \eqref{E:4.18ii}. By the Conversion Lemma, the second term tends to $[E_{\Q_\omega}(g_r^{\alpha'})]^{1/\alpha'}$ which vanishes as~$r\to\infty$. Since Theorem~\ref{thm-2.1} applies to~$f_r$, the claim follows.
\end{proofsect}

With Lemma~\ref{lemma-4.2w} in hand we are ready to give:

\begin{proofsect}{Proof of IIP in Theorem~\ref{thm-2.4} in dimension one}
Suppose~$d=1$, pick $p,q>1$ and let~$\omega\in\Omega^\star\cap\Omega_{p,q}$. Abbreviate
\begin{equation}
\label{E:4.20i}
f(\omega):=\frac1{\cc_\omega(0,1)}
\end{equation}
and let~$f_r$ be as in \eqref{E:3.12i}. Then $f_r\in\Cb(\Omega)$ and so the argument in \eqref{E:3.13i} (albeit without the use of Jensen's inequality) bounds $E_{\BbbP_\omega}(|f_r|^q)$ by the second supremum in \eqref{E:2.6w}. Taking~$r\to\infty$ using the Monotone Convergence Theorem gives $f\in L^q(\BbbP_\omega)$.

We will now proceed by the corrector method while capitalizing on the fact that the corrector can be constructed in~$d=1$ explicitly. Indeed, define~$\psi_\omega\colon\Z\to\R$ recursively so that
\begin{equation}
\label{E:2.19}
\psi_\omega(0) =0 \,\,\wedge\,\,\forall x\in\Z\colon\,\, \psi_\omega(x+1)-\psi_\omega(x) = \frac a{\cc_\omega(x,x+1)}
\end{equation}
for
\begin{equation}
\label{E:4.25u}
a:=\biggl[E_{\BbbP_\omega}\Bigl(\frac1{\cc(0,1)}\Bigr)\biggr]^{-1},
\end{equation} 
where the expectation is finite thanks to $f\in L^q(\BbbP_\omega)$.
Then $\psi_\omega$ is harmonic with respect to~$\cmss P_\omega$ and thus defines a harmonic coordinate provided we can show the ``sublinearity of the corrector,''
\begin{equation}
\label{E:2.12}
\lim_{x\to\pm\infty}\frac{\psi_\omega(x)-x}{|x|}=0.
\end{equation}
We will for simplicity prove only the part~$x\to+\infty$. Note that, for each~$n\ge1$,
\begin{equation}
%\label{}
\psi_\omega(n)= a\sum_{x=0}^{n-1}f\circ\tau_x(\omega),
\end{equation}
where~$f$ is as in \eqref{E:4.20i}. Lemma~\ref{lemma-3.1} with~$h:=f$ and \eqref{E:3.11i} enabled by the second condition in \eqref{E:2.6w} then shows 
\begin{equation}
%\label{}
\lim_{n\to\infty}\frac{\psi_\omega(n)}n = a E_{\BbbP_\omega}(f).
\end{equation}
The right-hand side equals one and so \eqref{E:2.12} follows.

Let~$X$ be the path of the random walk and set, as before,~$\FF_k:=\sigma(X_0,\dots,X_k)$. The harmonicity of~$\psi_\omega$ translates into~$k\mapsto\psi_\omega(X_k)$ being a martingale with respect to $\{\FF_k\}_{k\ge0}$ and so we will aim to control its diffusive limit using the Martingale Functional Central Limit Theorem of Brown~\cite{Brown} (see also \cite[Theorem 2.11]{B11} and discussion thereafter). The key point is to verify the conditions of this theorem. 

The Markov property along with \eqref{E:2.19} give
\begin{equation}
%\label{}
E^0_\omega\bigl([\psi_\omega(X_k)-\psi_\omega(X_{k-1})]^2\big|\FF_k\bigr) = g\circ\tau_{X_{k-1}}(\omega)
\end{equation}
for
\begin{equation}
\label{E:4.30u}
g(\omega'):=E^0_{\omega'}\bigl([\psi_{\omega'}(X_1)]^2\bigr) 
= \frac{a^2}{\pi_{\omega'}(0)}\biggl(\frac{1}{\cc_{\omega'}(0,1)}+\frac {1}{\cc_{\omega'}(-1,0)}\biggr).
\end{equation}
 Note that $g=a^2 (f+f\circ\tau_{-1})\pi(0)^{-1}$ for~$f$ as above. Given any~$\alpha\in(1,q)$, the convexity of~$s\mapsto s^\alpha$ gives $g^\alpha\pi(0)\le 2^{\alpha-1}a^{2\alpha}(f^\alpha\pi(0)^{\alpha-1}+(f\circ\tau_{-1})^\alpha\pi(0)^{\alpha-1})$ and the H\"older inequality with parameters $q/\alpha$ and~$\frac q{q-\alpha}$ shows
\begin{equation}
\label{E:4.31i}
\begin{aligned}
\frac1{n^d}\sum_{x\in\Lambda_n}&\bigl[g\circ\tau_x(\omega)\bigr]^\alpha\pi_\omega(x)
\\
&\le 2^{\alpha}a^{2\alpha}\biggl(\frac1{n^d}\sum_{x\in\Lambda_{n+1}} \bigl[f\circ\tau_x(\omega)\bigr]^q\biggr)^{\alpha/q}\biggl(\frac1{n^d}\sum_{x\in\Lambda_n}\pi_\omega(x)^{\frac{q(\alpha-1)}{q-\alpha}}\biggr)^{\frac{q-\alpha}q}.
\end{aligned}
\end{equation}
The term in the first large parentheses on the right is bounded uniformly in~$n\ge1$ thanks to~$\omega\in\Omega_{p,q}$ and the same applies to the term in the second parentheses provided $\alpha-1$ is so small that $\frac{q(\alpha-1)}{q-\alpha}\le p$.
It follows that \eqref{E:4.18i} holds with~$h:=g$ and, since~$g\in C(\Omega)$, the  extension of Theorem~\ref{thm-2.1} in Lemma~\ref{lemma-4.2w} gives
\begin{equation}
\label{E:3.31y}
\frac1n\sum_{k=1}^{n}E^0_\omega\bigl([\psi_\omega(X_k)-\psi_\omega(X_{k-1})]^2\big|\FF_k\bigr)
\,\,\underset{n\to\infty}{\overset{L^1(P^0_\omega)}\longrightarrow}\,\,E_{\Q_\omega}(g).
\end{equation}
Looking instead at the function
\begin{equation}
\label{E:3.30u}
\tilde g_r(\omega'):=E^0_{\omega'}\Bigl(\bigl[\psi_{\omega'}(X_1)\bigr]^2 \vartheta\bigl(|\psi_{\omega'}(X_1)|/r\bigr)\Bigr)
\end{equation}
 for some continuous non-decreasing $\vartheta\colon[0,\infty)\to[0,1]$ such that $\vartheta=0$ on $[0,1]$ and $\vartheta=1$ on $[2,\infty)$, we similarly get that, for each~$\epsilon>0$,
\begin{equation}
\label{E:2.25}
\frac1n\sum_{k=1}^{n}E^0_\omega\Bigl(\bigl[\psi_\omega(X_k)-\psi_\omega(X_{k-1})\bigr]^2\,1_{\{|\psi_\omega(X_k)-\psi_\omega(X_{k-1})|>\epsilon\sqrt{n}\}}\,\Big|\,\FF_k\Bigr)
\,\,\underset{n\to\infty}{\overset{L^1(P^0_\omega)}\longrightarrow}\,\,0.
\end{equation}
Indeed, once~$n$ is so large that $2r\le\epsilon\sqrt n$, the sum is dominated by $\sum_{k\ge1}\tilde g_r\circ\tau_{X_k}(\omega)$ which upon normalization by~$n$ converges to $E_{\Q_\omega}(\tilde g_r)$ by Theorem~\ref{thm-2.1}  and the argument \eqref{E:4.31i} combined with~$\tilde g_r\le g$.  The Dominated Convergence Theorem along with~$\tilde g_r\le g$ then show $E_{\Q_\omega}(\tilde g_r)\to0$ as $r\to\infty$ proving \eqref{E:2.25}. 

Whenever conditions \eqref{E:3.31y} and \eqref{E:2.25} hold, the Martingale Functional Central Limit Theorem shows that the process
\begin{equation}
%\label{}
t\mapsto n^{-1/2}\bigl[\psi_\omega(X_{\lfloor nt\rfloor})+(tn-\lfloor tn\rfloor)\psi_\omega(X_{\lfloor nt\rfloor+1})\bigr]
\end{equation}
tends in law to Brownian motion with variance~$E_{\Q_\omega}(g)$. In light of \eqref{E:2.12}, then so does the process $t\mapsto B_t^{(n)}$ defined from~$X$ via \eqref{E:1.3}. 

It remains to check that the variance $E_{\Q_\omega}(g)$ is given by \eqref{E:1.5}.  A standard argument based on the Parallelogram Law shows that, for any minimizing sequence $\{\varphi_n\}_{n\ge1}$ in \eqref{E:1.5} and any~$z\in\Z$, the functions $\varphi_n\circ\tau_z-\varphi_n$ converge to some~$\chi(\cdot,z)$ in the sense
\begin{equation}
%\label{}
E_{\BbbP_\omega}\biggl(\,\sum_{x=\pm1} \cc(0,x)\bigl|\varphi_n\circ\tau_x-\varphi_n-\chi(\cdot,x)\bigr|^2\biggr)\,\,\underset{n\to\infty}\longrightarrow\,\,0.
\end{equation}
Using the Cauchy-Schwarz inequality and $\cc(0,x)^{-1}\in L^1(\BbbP_\omega)$, hereby we get that also
$\varphi_n\circ\tau_z-\varphi_n\to\chi(\cdot,z)$ in~$L^1(\BbbP_\omega)$ which implies
\begin{equation}
%\label{}
\chi(\cdot,z)\in L^1(\BbbP_\omega)\,\,\wedge\,\,E_{\BbbP_\omega}\bigl(\chi(\cdot,z)\bigr)=0
\end{equation}
for all~$z\in\Z$.
Since $\chi(\omega,n) = \sum_{k=0}^{n-1}\chi(\tau_k(\omega),1)$, the ergodicity of~$\BbbP_\omega$ and Birkhoff's Pointwise Ergodic Theorem imply that $\chi(\cdot,x)=o(|x|)$ as~$|x|\to\infty$.

Examining the Euler-Lagrange equations in turn shows that $\tilde\psi(\omega,x):= x+\chi(\omega',x)$ is~$\cmss P_{\omega'}$-harmonic for~$\BbbP_\omega$-a.e.~$\omega'$, which readily translates into
\begin{equation}
%\label{}
\tilde\psi(\omega',x+1)-\tilde\psi(\omega',x)= \frac{a(\omega')}{\cc_{\omega'}(x,x+1)}
\end{equation}
with~$a(\omega'):=\cc_{\omega'}(0,1)\tilde\psi(\omega',1)$. This gives
\begin{equation}
%\label{}
n+\chi(\omega',n)=\tilde\psi(\omega',n) = a(\omega')\sum_{k=0}^{n-1}\frac1{\cc_{\omega'}(k,k+1)}
\end{equation}
which upon division by~$n$ and taking~$n\to\infty$ with the help of Birkhoff's Ergodic Theorem, the ergodicity of~$\BbbP_\omega$ and~$\cc(0,x)^{-1}\in L^1(\BbbP_\omega)$ shows that~$a(\omega')$ is constant and equal to~$a$ in \eqref{E:4.25u} $\BbbP_\omega$-a.s. It follows that~$\tilde\psi(\omega',x)=\psi_{\omega'}(x)$ for $\BbbP_\omega$-a.e.~$\omega'$\ which in light of \eqref{E:4.30u} shows that~$E_{\Q_\omega}(g)$ equals the quantity in \eqref{E:1.5} as desired.
\end{proofsect}

%For this observe that the functional optimized in \eqref{E:1.5} takes the form
%\begin{equation}
%\label{E:3.35}
%\varphi\mapsto\int \Q_\omega(\textd \omega') E^0_{\omega'}\Bigl(\bigl[X_1+\varphi\circ\tau_{X_1}(\omega')-\varphi(\omega')\bigr]^2\Bigr).
%\end{equation}
%Since $c(0,1),c(0,1)^{-1}\in L^1(\BbbP_\omega)$, for any minimizing sequence $\{\varphi_n\colon n\ge1\}$, the standard domination of~$L^1$-norms by weighted $L^2$-norms (cf, e.g., proof of~\cite[Lemma~4.8]{B11}) shows that $\varphi_n\circ\tau_x-\varphi_n$ tends to some~$\chi(x)$ in~$L^1(\BbbP_\omega)$ for each~$x\in\Z^d$. Then $E_{\BbbP_\omega}(\chi(x))=0$ and, using the fact that this constructs a minimizer of \eqref{E:3.35}, $x\mapsto x+\chi_{\omega'}(x)$ is $\cmss P_{\omega'}$-harmonic for~$\BbbP_\omega$-a.e.~$\omega'\in\Omega$. Corollary~5.6 of Biskup and Spohn~\cite{BS11} (based on \cite[Theorem~5.4]{BS11} whose proof needs a correction albeit only in~$d\ge2$; cf Biskup and Rodriguez~\cite[Appendix~A2]{BR18}), asserts that these properties determine~$\chi$ uniquely and so $x+\chi(x)=\psi(x)$ $\BbbP_\omega$-a.s. We now readily identify the minimum of \eqref{E:3.35} with~$E_{\Q_\omega}(g)$.

%\newpage
\section{Homogenization of inverse generator}
\label{sec5}\noindent
In this section we prove a technical claim needed in our control of long-time behavior of random walks among deterministic conductances. This claim concerns homogenization of inverse generator of the physical walk killed upon exiting a finite set. Our strategy is to convert this to a  homogenization problem for Dirichlet energy.

\subsection{Notation and statement}
We start by introducing notation for various relevant objects and structures needed in statements and proofs in this section. The main objective is to link concepts from (elementary) differential geometry in the physical space to the corresponding concepts in the stochastic counterpart thereof.

The first item on the list are the notions of gradient and divergence.
Let $e_i$ denote the $i$-th unit coordinate vector in~$\R^d$. The gradient $\nabla$ maps a function~$f\colon\Omega\to\R$ to a $d$-tuple of functions written in vector notation as
\begin{equation}
%\label{}
(\nabla f)(\omega):=\bigl(f\circ\tau_{e_1}(\omega)-f(\omega),\dots,f\circ\tau_{e_d}(\omega)-f(\omega)\bigr).
\end{equation}
The divergence $\nabla^\star$ in turn acts on $d$-tuples $v=(v_1,\dots v_d)$ of functions via
\begin{equation}
%\label{}
(\nabla^\star v)(\omega) := \sum_{i=1}^d\bigl[v_i(\omega)- v_i\circ\tau_{-e_i}(\omega)\bigr]
\end{equation}
thus turning them back into a scalar function on~$\Omega$.  These operations naturally introduce vector-valued functions into play; we will use $C(\Omega,\R^d)$ to denote the set of continuous functions $\Omega\to \R^d$ (sometimes called vector fields below) and write $\Cb(\Omega,\R^d)$ for the subspace $C(\Omega,\R^d)$ consisting of bounded functions. 

Next we need a notation for two Hilbert spaces naturally associated with the Markov chains in the physical and stochastic context. Given any translation invariant law~$\mu$ on $(\Omega,\FF)$, write $L^2(\mu,\R^d)$ for the closure of~$\Cb(\Omega,\R^d)$ in the topology induced by the inner product
\begin{equation}
%\label{}
\langle f,g\rangle_{L^2(\mu,\R^d)}:=\int_{\Omega} f(\omega')\cdot g(\omega')\mu(\textd \omega'),
\end{equation}
 where the Euclidean dot product is used to contract the coordinates under the integral. 
Similarly, given a finite set~$\Lambda\subseteq\Z^d$ and a function~$\nu\colon\Lambda\to[0,\infty)$, write $\ell^2(\Lambda,\nu)$ be the linear vector space of functions $\Z^d\to\R^d$ that vanish outside~$\Lambda$; this space is endowed with the inner product 
\begin{equation}
%\label{}
\langle\, f,g\rangle_{\ell^2(\Lambda,\nu)}:=\sum_{x\in\Lambda} f(x)\cdot g(x)\nu(x).
\end{equation}
The natural choices here are $\mu:=\Q_\omega$ (defined in \eqref{E:3.5t}) and~$\nu:=\pi_\omega$.

Given a translation-invariant probability measure~$\mu$ on~$(\Omega,\FF)$ with~$\pi(0)\in L^1(\mu)$, the (negative) generator~$-\LL$  (defined in \eqref{E:4.1w})  of the chain from the ``point of view of the particle'' is bounded, self-adjoint and positive semi-definite on~$L^2(\pi\mu,\R)$, where~$\pi\mu$ is a shorthand for the measure on~$(\Omega,\FF)$ defined by $\pi\mu(A):=\int_A\pi_{\omega'}(0)\mu(\textd\omega')$.
Methods of spectral calculus permit us to define $(-\LL)^{-1/2}f$ for all~$f\in L^2(\pi\mu,\R)$ satisfying
\begin{equation}
\label{E:3.6v}
\sup_{\lambda>0}\,\bigl\langle\, f,(\lambda-\LL)^{-1}f\bigr\rangle_{L^2(\pi\mu,\R)}<\infty.
\end{equation}
The supremum then equals the norm-squared of the vector $(-\LL)^{-1/2}f$ in~$L^2(\pi\mu,\R)$. As usual, we will write \eqref{E:3.6v} shortly as
\begin{equation}
%\label{}
\langle\, f,(-\LL)^{-1}f\rangle_{L^2(\pi\mu,\R)}<\infty
\end{equation}
in the sequel.

The connection with the above notions of gradient and divergence arises from the fact that~$\nabla^\star$ is a formal adjoint of~$-\nabla$ in the sense that (continuing to assume~$\mu$ to be translation invariant)
\begin{equation}
\label{E:3.5}
\langle \nabla^\star v, f\rangle_{L^2(\mu,\R)} = -\langle v, \nabla f\rangle_{L^2(\mu,\R^d)}
\end{equation}
holds for each~$f\in\Cb(\Omega)$ and each~$v\in \Cb(\Omega,\R^d)$. 
Writing $\cc$ be the operator of ``multiplication by conductance value'' defined for $v=(v_1,\dots,v_d)$ by
\begin{equation}
%\label{}
(\cc v)(\omega):=\bigl(\omega(0,e_1)v_1(\omega),\dots,\omega(0,e_d)v_d(\omega)\bigr),
\end{equation}
the generator~$\LL$ then acts as
\begin{equation}
\label{E:3.11}
(\LL f)(\omega)= \pi_\omega(0)^{-1}\nabla^\star(\cc\nabla f)(\omega).
\end{equation}
giving it a ``divergence-form'' structure.

While $\cmss P_\omega$ denotes the transition probability of the random walk in the physical space under environment~$\omega$, the same notation is often used for the linear operator mapping $f\colon\Z^d\to\R$ to $(\cmss P_\omega f)(x):=\sum_{y\in\Z^d}\cmss P_\omega(x,y)f(y)$. The operator
\begin{equation}
\label{E:4.8w}
\cmss L_\omega:=\cmss P_\omega-1
\end{equation}
is the generator of the chain in the physical space. For each~$\omega\in\Omega$, the operator $-\cmss L_\omega$ restricted (by projections which we do not mark explicitly) to $\ell^2(\Lambda,\pi_\omega)$ is self-adjoint and positive definite and thus admits an inverse that we denote as $(-\cmss L_\omega)^{-1}_{\Lambda}$. Recalling the definition of~$\Lambda_r$ from~\eqref{E:2.2i} we now claim:

\begin{theorem}
\label{thm-3.1}
Assume that~$f\in \Cb(\Omega)$ takes the form $f(\omega):=\pi_\omega(0)^{-1}\nabla^\star (\cc v)(\omega)$ for some $v\in \Cb(\Omega,\R^d)$. For all $p,q>1$ and all~$\omega\in\Omega^\star\cap\Omega_{p,q}$, 
\begin{equation}
\label{E:3.8b}
\langle\, f,(-\LL)^{-1}f\rangle_{L^2(\Q_\omega,\R)}<\infty
\end{equation}
and,  if~$p,q$ also obey $1/p+1/q<2/d$ in~$d\ge3$, then 
\begin{equation}
\label{E:3.9}
\frac1{\pi_\omega(\Lambda_r)}\bigl\langle\,f\circ\tau(\omega),\,(-\cmss L_\omega)^{-1}_{\Lambda_r} f\circ\tau(\omega)\bigr\rangle_{\ell^2(\Lambda_r,\pi_\omega)}
\,\,\underset{r\to\infty}\longrightarrow\,\, \bigl\langle\, f,\,(-\LL)^{-1} f\bigr\rangle_{L^2(\Q_\omega,\R)}.
\end{equation}
Here~$f\circ\tau(\omega)$ abbreviates the function~$x\mapsto f\circ\tau_x(\omega)$.
\end{theorem}

Note that, while we are assuming only the first-plus ``moment conditions'' on~$\omega$ in dimensions $d=1,2$,  in~$d\ge3$ we have to resort back to the condition in \eqref{E:2.6}. This is because the heat-kernel technology is ultimately invoked in the proof, albeit only in a very final step; see Lemma~\ref{lemma-5.9}. We believe that this is not needed, but have not been able to come up with a better argument. See Remark~\ref{rem-5.10} for more discussion. 

The restriction to~$f$ of the  gradient  form $f=\pi(0)^{-1}\nabla^\star(\cc v)$ is made as this suffices for our purposes (and gives the claim at the same time). We do not know if the convergence \eqref{E:3.9} holds when this is relaxed, e.g., by assuming only~\eqref{E:3.8b} instead. 

Our application of Theorem~\ref{thm-3.1} only needs that the \emph{limes superior} in \eqref{E:3.9} is bounded by the right-hand side, but we prove both directions as the argument for the \emph{limes inferior} is easier and the claim is more informative.

\subsection{Reduction to Dirichlet energy}
The proof of \eqref{E:3.8b} uses the idea that initiates the derivations in Kipnis and Varadhan's~\cite{KV86} approach to random walks in reversible random environments. Roughly speaking, we will rewrite the inner product for a function~$f$ as an inner product for the vector field~$v$, thus  getting rid of   the inverse generator whose direct control may be difficult.  We illustrate this on an observation drawn more or less directly from~\cite{KV86}: 

\begin{lemma}
\label{lemma-3.2}
Given a translation invariant probability measure~$\mu$ on~$(\Omega,\FF)$ with $\pi(0)\in L^1(\mu)$, for each~$v\in \Cb(\Omega,\R^d)$ and with $f\in \Cb(\Omega)$ defined by $f(\omega):=\pi_\omega(0)^{-1}\nabla^\star (\cc v)(\omega)$,
\begin{equation}
\label{E:3.12a}
\langle\, f,(-\LL)^{-1}f\rangle_{L^2(\pi\mu,\R)}\le \langle v, \cc v\rangle_{L^2(\mu,\R^d)}.
\end{equation}
The right-hand side is finite for each~$v\in \Cb(\Omega,\R^d)$.
\end{lemma}

\begin{proofsect}{Proof}
The assumption $\pi(0)\in L^1(\mu)$ ensures that $\Cb(\Omega)\subseteq L^2(\pi\mu,\R)$. For each~$h\in \Cb(\Omega)$, \eqref{E:3.5} and the Cauchy-Schwarz inequality give
\begin{equation}
\label{E:3.29w}
\begin{aligned}
\bigl|\langle h,f\rangle_{L^2(\pi\mu,\R)}\bigr| &= \bigl|\langle h,\nabla^\star(\cc v)\rangle_{L^2(\mu,\R)}\bigr|
\\
&=\bigl|\langle\nabla h,\cc v\rangle_{L^2(\mu,\R^d)}\bigr|
\le \langle\nabla h,\cc \nabla h\rangle_{L^2(\mu,\R^d)}^{1/2}\langle v,\cc v\rangle_{L^2(\mu,\R^d)}^{1/2}
\\
&=  \langle h,-\LL h\rangle_{L^2(\pi\mu,\R)}^{1/2}\langle v,\cc v\rangle_{L^2(\mu,\R^d)}^{1/2}.
\end{aligned}
\end{equation}
Taking $h:=(\lambda-\LL)^{-1} f$ for $\lambda>0$, which is well defined by the above mentioned properties of~$\LL$, yields
\begin{equation}
%\label{}
\langle f,(\lambda-\LL)^{-1} f\rangle_{L^2(\pi\mu,\R)}\le \langle f,\,(-\LL)(\lambda-\LL)^{-2} f\rangle_{L^2(\pi\mu,\R)}^{1/2}\langle v,\cc v\rangle_{L^2(\mu,\R^d)}^{1/2}.
\end{equation}
Since $\Vert\LL(\lambda-\LL)^{-1}\Vert\le1$, the first inner product on the right is dominated by that on the left-hand side.  Since the second inner product is finite due to $\pi(0)\in L^1(\mu)$, this gives both stated claims. 
\end{proofsect}

While \eqref{E:3.12a} will imply \eqref{E:3.8b}, a problem with its further
use is that the inequality is far away from sharp. This is because the gap in the inequality in \eqref{E:3.29w} is small only if~$v$ is close to a gradient. As it turns out, the proof of \eqref{E:3.9} actually requires  us  to compute exactly the gap in the inequality in \eqref{E:3.12a}. This  boils down  to a computation of the distance of~$v$ from the linear subspace of vector fields that are gradients.

We will address this first in the physical space, for which further notation is needed. Recall that, for each~$\Lambda\subseteq\Z^d$, we write $E(\Lambda)$ for the set of directed edges of~$\Z^d$ with at least one endpoint in~$\Lambda$. Given edge weights $\kappa\colon E(\Lambda)\to[0,\infty)$  subject to the symmetry condition $\kappa(x,y)=\kappa(y,x)$,  let $\ell^2(E(\Lambda),\kappa)$ be  the linear vector space of functions $u\colon E(\Lambda)\to\R$ such that
\begin{equation}
\label{E:4.25p}
\forall (x,y)\in E(\Lambda)\colon\,\,u(x,y)=-u(y,x)
\end{equation}
endowed with the inner product 
\begin{equation}
\label{E:4.25w}
\langle u,v\rangle_{\ell^2(E(\Lambda),\kappa)}:=\frac12\sum_{(x,y)\in E(\Lambda)} u(x,y)v(x,y)\kappa(x,y).
\end{equation}
The factor $1/2$ compensates for overcounting caused by both orientations of each edge being included in~$E(\Lambda)$. 

\begin{remark}
\label{rem-4.3}
If~$u$ and~$v$ also~$\R^d$ valued, then we also contract the indices using the dot product; i.e., writing $u(x,y)\cdot v(x,y)$ instead of just $u(x,y)v(x,y)$ in \eqref{E:4.25w}. 
\end{remark}

We think if~$u$ satisfying \eqref{E:4.25p} as a ``flow'' in the physical space.
A representative example of such a ``flow'' in all of~$\Z^d$ is $u(x,y):=\alpha\cdot(y-x)$ for~$\alpha\in\R^d$ corresponding to constant vector field $\alpha\in C(\Omega,\R^d)$. Other examples arise from taking the \emph{discrete gradient}~$\nabla f$ of a function~$f\colon\Z^d\to\R$ which is defined for each edge $(x,y)$ by
\begin{equation}
\label{E:4.27p}
(\nabla f)(x,y):=f(y)-f(x).
\end{equation}
(The difference between the discrete gradient and the gradient on~$\Omega$ will be clear from context or will be noted explicitly.)
A natural choice of edge weights is $\kappa:=\cc_\omega$. 
With these notions in hand, we now state:

\begin{lemma}
\label{lemma-3.5}
Given~$v\in\Cb(\Omega,\R^d)$, let~$f\in\Cb(\Omega)$ be defined by $f(\omega):=\pi_\omega(0)^{-1}\nabla^\star (\cc v)(\omega)$. Then for all finite $\Lambda\subseteq\Z^d$ and all~$\omega\in\Omega$,
\begin{equation}
\label{E:3.24a}
\begin{aligned}
\bigl\langle\,f\circ\tau(\omega),&\,(-\cmss L_\omega)^{-1}_{\Lambda} f\circ\tau(\omega)\bigr\rangle_{\ell^2(\Lambda,\pi_\omega)}
\\&= \bigl\langle\,v\circ\tau(\omega),\,v\circ\tau(\omega)\bigr\rangle_{\ell^2(E(\Lambda),\,\cc_\omega)}
\\&\quad\qquad-\inf_{h\in\ell^2(\Lambda,\pi_\omega)}
\bigl\langle\,v\circ\tau(\omega)-\nabla h,\,v\circ\tau(\omega)-\nabla h\bigr\rangle_{\ell^2(E(\Lambda),\,\cc_\omega)},
\end{aligned}
\end{equation}
where $v\circ\tau(\omega)$ is the function $E(\Z^d)\to\R$ satisfying \eqref{E:4.25p} (with~$\Lambda:=\Z^d$) and assigning value $v_i\circ\tau_x(\omega)$ to edge $(x,x+e_i)$.
\end{lemma}

\begin{proofsect}{Proof}
Fix~$\omega\in\Omega$, abbreviate $f'(x):=f\circ\tau_x(\omega)$ and $v'(x,x+e_i)=-v'(x+e_i,x):= v_i\circ\tau_x(\omega)$ and use the shorthand $\langle f,g\rangle_{\Lambda} = \langle f,g\rangle_{\ell^2(\Lambda,\pi_\omega)}$. The invertibility and positive-definiteness of~$-\cmss L_\omega$ on~$\ell^2(\Lambda,\pi_\omega)$ then give
\begin{equation}
\label{E:4.15o}
\bigl\langle\,f',\,(-\cmss L_\omega)^{-1}_{\Lambda} f'\bigr\rangle_{\Lambda}=\sup_{h\in \ell^2(\Lambda,\pi_\omega)}\Bigl[ 2\,\langle h,f'\rangle_{\Lambda}-\bigl\langle h,-\cmss L_\omega h\rangle_{\Lambda}\Bigr].
\end{equation}
The assumption on~$f$ along with the definition of~$\nabla^\star$ shows that, for each~$x$,
\begin{equation}
%\label{}
f'(x) = \pi_\omega(x)^{-1}\sum_{y\colon(x,y)\in E(\Z^d)} \cc_\omega(x,y)v'(x,y).
\end{equation}
As~$h$ vanishes outside~$\Lambda$, this yields
\begin{equation}
%\label{}
\begin{aligned}
\langle h,f'\rangle_{\Lambda} &= \sum_{(x,y)\in E(\Lambda)} h(x)\cc_\omega(x,y)v'(x,y)
\\&=\frac12\sum_{(x,y)\in E(\Lambda)} \bigl[h(x)-h(y)\bigr]\cc_\omega(x,y)v'(x,y)
=-\langle \nabla h,v'\rangle_{E(\Lambda)},
\end{aligned}
\end{equation}
where we invoked the symmetries $\cc_\omega(x,y)=\cc_\omega(y,x)$ and $v'(x,y)=-v'(y,x)$, recalled the notation \eqref{E:4.27p} for the discrete gradient and wrote $\langle u,v\rangle_{E(\Lambda)}$ for $\langle u,v\rangle_{\ell^2(E(\Lambda),\,\cc_\omega)}$. A completely analogous reasoning shows
\begin{equation}
\label{E:4.25x}
\bigl\langle h,-\cmss L_\omega h\rangle_{\Lambda} = \langle \nabla h,\nabla h\rangle_{E(\Lambda)}.
\end{equation}
But then
\begin{equation}
\begin{aligned}
\label{E:3.29}
\qquad
2\,\langle h,f'\rangle_{\Lambda}-\bigl\langle h,-\cmss L_\omega h\rangle_{\Lambda}
&=-2\langle \nabla h,v'\rangle_{E(\Lambda)}-\langle \nabla h,\nabla h\rangle_{E(\Lambda)}
\\
&=\langle v',v'\rangle_{E(\Lambda)}-\langle v'+\nabla h,v'+\nabla h\rangle_{E(\Lambda)}.
\end{aligned}
\end{equation}
Plugging this in \eqref{E:4.15o} yields the claim.
\end{proofsect}

We will henceforth abbreviate the infimum  in \eqref{E:3.24a}  as
\begin{equation}
\label{E:4.26y}
\EE_\Lambda^\omega(v):=\inf_{h\in\ell^2(\Lambda,\pi_\omega)}
\bigl\langle\,v\circ\tau(\omega)-\nabla h,\,v\circ\tau(\omega)-\nabla h\bigr\rangle_{\ell^2(E(\Lambda),\,\cc_\omega)}.
\end{equation}
Since the inner product has a representation reminiscent of a Dirichlet form, we will refer to $\EE_\Lambda^\omega(v)$ as the \emph{Dirichlet energy}.

In light of the positivity of all terms and the fact that only nearest neighbor edges are involved, for the first inner product on the right of \eqref{E:3.24a} we get the bounds
\begin{equation}
\label{E:4.26}
\sum_{x\in\Lambda_r}g\circ\tau_x(\omega)\le \bigl\langle\,v\circ\tau(\omega),\,v\circ\tau(\omega)\bigr\rangle_{\ell^2(E(\Lambda),\,\cc_\omega)}
\le \sum_{x\in\Lambda_{r+1}}g\circ\tau_x(\omega),
\end{equation}
where
\begin{equation}
%\label{}
g(\omega'):=\sum_{i=1}^d \cc_{\omega'} v_i(\omega')^2.
\end{equation}
If $\omega\in\Omega^\star\cap\Omega_{p,0}$ for some~$p>1$, then Lemma~\ref{lemma-3.1} shows
\begin{equation}
\label{E:4.29}
\frac1{|\Lambda_r|}\bigl\langle\,v\circ\tau(\omega),\,v\circ\tau(\omega)\bigr\rangle_{\ell^2(E(\Lambda_r),\,\cc_\omega)}\,\underset{r\to\infty}\longrightarrow\,\langle v,\cc v\rangle_{L^2(\BbbP_\omega,\R^d)}.
\end{equation}
The corresponding limit of the second term on the right of \eqref{E:3.24a} is harder due to the infimum over~$h$  getting in the way.  In fact, the main conclusion to be proved is that the limit effectively commutes around the infimum:

\begin{proposition}
\label{prop-4.4}
Let $v\in\Cb(\Omega,\R^d)$. Then for all~$p,q>1$  such that~$1/p+1/q<2/d$ in~$d\ge3$  and all~$\omega\in\Omega^\star\cap\Omega_{p,q}$,
\begin{equation}
%\label{}
\qquad\frac1{|\Lambda_r|}
\EE_{\Lambda_r}^\omega(v)
\,\,\underset{r\to\infty}\longrightarrow\,\,
\inf_{h\in\Cb(\Omega)}\bigl\langle v-\nabla h,\cc(v-\nabla h)\bigr\rangle_{L^2(\BbbP_\omega,\R^d)}.
\qquad
\end{equation}
\end{proposition}

Let us first check that this is  indeed  all we need:

\begin{proofsect}{Proof of Theorem~\ref{thm-3.1} from Proposition~\ref{prop-4.4}}
The bound \eqref{E:3.8b} was proved in Lemma~\ref{lemma-3.2}. For \eqref{E:3.9}, Lemma~\ref{lemma-3.5} along with Proposition~\ref{prop-4.4} and \eqref{E:4.29} give
\begin{equation}
\begin{aligned}
\frac1{|\Lambda_r|}\bigl\langle\,f\circ\tau(&\omega),\,(-\cmss L_\omega)^{-1}_{\Lambda_r} f\circ\tau(\omega)\bigr\rangle_{\ell^2(\Lambda_r,\pi_\omega)}
\\&\underset{r\to\infty}\longrightarrow\,\,
\langle v,\cc v\rangle_{L^2(\BbbP_\omega,\R^d)}-
\inf_{h\in\Cb(\Omega)}\bigl\langle v-\nabla h,\cc(v-\nabla h)\bigr\rangle_{L^2(\BbbP_\omega,\R^d)}
\end{aligned}
\end{equation}
whenever~$\omega\in\Omega^\star\cap\Omega_{p,q}$ for some~$p,q>1$  as in Proposition~\ref{prop-4.4}.  Reversing the arguments  \twoeqref{E:4.25x}{E:3.29} equates the right-hand side with $E_{\BbbP_\omega}(\pi(0))$-multiple of
\begin{equation}
%\label{}
\sup_{h\in\Cb(\Omega)}\Bigl[-2\langle h,f\rangle_{L^2(\Q_\omega,\R)}-\langle h,(-\LL)h\rangle_{L^2(\Q_\omega,\R)}\Bigr]
\end{equation}
which is then identified with $\langle f,(-\LL)^{-1} f\rangle_{L^2(\Q_\omega,\R)}$ as in \eqref{E:4.15o}. Since the assumptions on~$\omega$ also give $\pi_\omega(\Lambda_r)/|\Lambda_r|\to E_{\BbbP_\omega}(\pi(0))$, this proves the desired claim.
\end{proofsect}

\subsection{Homogenization of Dirichlet energy}
The proof of Theorem~\ref{thm-3.1} has been reduced to the limit of Dirichlet energies in Proposition~\ref{prop-4.4} which we will show by proving separately upper and lower bounds. The proof of the upper bound is quite simple:

\begin{lemma}
\label{lemma-4.5}
Let $v\in\Cb(\Omega,\R^d)$. Then for all~$p,q>1$ and all~$\omega\in\Omega^\star\cap\Omega_{p,q}$,
\begin{equation}
%\label{}
\limsup_{r\to\infty}\frac1{|\Lambda_r|}\EE_{\Lambda_r}^\omega(v)
\le
\inf_{h\in\Cb(\Omega)}\bigl\langle v-\nabla h,\cc(v-\nabla h)\bigr\rangle_{L^2(\BbbP_\omega,\R^d)}.
\end{equation}
\end{lemma}

\begin{proofsect}{Proof}
Given~$h\in\Cb(\Omega)$ and~$\omega\in\Omega$, for each~$r\ge1$, set $h'_r(x):=h\circ\tau_x(\omega)$ for $x\in\Lambda_r$ and~$h'_r(x):=0$ otherwise. Writing $\langle\cdot,\cdot\rangle_{E(\Lambda_r)}$  as a shorthand  for $\langle\cdot,\cdot\rangle_{\ell^2(E(\Lambda_r),\cc_\omega)}$,  the Minkowski inequality gives 
\begin{equation}
\begin{aligned}
 \EE_{\Lambda_r}^\omega(v)\le\bigl\langle\,&v\circ\tau(\omega)-\nabla h'_r,\,v\circ\tau(\omega)-\nabla h'_r\bigr\rangle_{E(\Lambda_r)}^{1/2}
\\
&\le\bigl\langle\,(v-\nabla h)\circ\tau(\omega),\,(v-\nabla h)\circ\tau(\omega)\bigr\rangle_{E(\Lambda_r)}^{1/2}
\\
&\qquad+\bigl\langle\,(\nabla h)\circ\tau(\omega)-\nabla h'_r,\,(\nabla h)\circ\tau(\omega)-\nabla h'_r\bigr\rangle_{E(\Lambda_r)}^{1/2}.
\end{aligned}
\end{equation}
(Here $\nabla h$ is a $\R^d$-valued function on~$\Omega$ while $(\nabla h)\circ\tau(\omega)$ and $\nabla h'_r$ are $\R$-valued functions on~$E(\Z^d)$.)
Since $(\nabla h)\circ\tau(\omega)-\nabla h'_r\ne 0$ only on edges between~$\Lambda_r$ and its complement, the second inner product on the right is at most order~$r^{d-1}$.  Upon  normalization by~$|\Lambda_r|$,  the same argument as  \twoeqref{E:4.26}{E:4.29} shows that the first inner product, and thus also the  whole right-hand  side,  tends to $\langle v-\nabla h,\cc(v-\nabla h)\rangle_{L^2(\Q_\omega)}$. Optimizing over~$h$ we get the claim.
\end{proofsect}

For the  corresponding  lower bound we have to work harder.  We start by noting that, since  the conductances are strictly positive, the infimum in \eqref{E:4.26y} is  achieved by a unique function $h_r^{\omega,v}\in\ell^2(\Lambda_r,\pi_\omega)$  that vanishes outside~$\Lambda_r$.  Writing $\VV:=(\R^d)^{\Z^d}$ for the space of vector fields endowed with the product $\sigma$-algebra~$\GG$ of Borel sets,  the pair~$(\omega,h^{\omega,v}_r)$ induces a  probability measure~$\mu_r^{\omega,v}$ on~$(\Omega\times\VV,\FF\otimes\GG)$ via
\begin{equation}
\label{E:4.35}
\mu_r^{\omega,v}(A):=\frac1{|\Lambda_{r+1}|}\sum_{x\in\Lambda_{r+1}}1_A\bigl((\tau_{x}(\omega),\nabla h_r^{\omega,v}(x+\cdot))\bigr),
\end{equation}
 where in the second coordinate we violate our earlier convention by regarding $\nabla h_r^{\omega,v}(x)$ as a $d$-dimensional vector whose $i$-th component equals $h_r^{\omega,v}(x+e_i)- h_r^{\omega,v}(x)$. 
Using~$\tau_x$ also for the joint shift of the variables in $\Omega\times\VV$, we observe the following facts:

\begin{lemma}
\label{lemma-5.7i}
Let $v\in\Cb(\Omega)$ and~$\omega\in\Omega^\star\cap\bigcup_{p,q>1}\Omega_{p,q}$. Then
\begin{equation}
\label{E:4.36}
\{\mu_r^{\omega,v}\colon r\ge1\}\,\,\text{\rm is tight}.
\end{equation}
Moreover, if $r_k\to\infty$ is a subsequence such that~$\mu_{r_k}^{\omega,v}\overset{\text{\rm w}}\longrightarrow\mu$, then also:
\settowidth{\leftmargini}{(1111)}
\begin{enumerate}
\item[(1)] $\mu$ is a probability measure, i.e., $\mu(\Omega\times\VV)=1$,
\item[(2)] $\mu$ is translation invariant, i.e., $\mu\circ\tau_x^{-1}=\mu$ for all~$x\in\Z^d$,
\item[(3)] the first marginal of~$\mu$ is~$\BbbP_\omega$, i.e.,
$\mu(A\times\VV)=\BbbP_\omega(A)$ for all~$A\in\FF$,
\item[(4)] the second marginal is integrable and of zero mean, i.e.,
\begin{equation}
\label{E:4.37}
\int_{\Omega\times\VV}|u(0)|_1\,\mu(\textd\omega'\textd u)<\infty\,\,\,\wedge\,\,\,\int_{\Omega\times\VV}u(0)\,\mu(\textd\omega'\textd u)=0,
\end{equation}
\item[(5)] energy is not generated in the limit,
\begin{equation}
\label{E:4.38}
\int_{\Omega\times\VV}\sum_{i=1}^d\cc_{\omega'}(0,e_i)\bigl[v_i(\omega')-u_i(0)\bigr]^2\,\mu(\textd\omega'\textd u)
\le
\liminf_{k\to\infty}\frac1{|\Lambda_{r_k}|}\EE_{\Lambda_{r_k}}^\omega(v).
\end{equation}
\end{enumerate}
Here, in~(4), $|\cdot|_1$ is the $\ell^1$-norm on~$\R^d$.
\end{lemma}

\begin{proofsect}{Proof}
Fix~$p,q>1$  and let~$\omega\in\Omega_{p,q}$.  We start with some estimates.   Fix $v\in\Cb(\Omega)$ and $r\ge1$ and abbreviate
\begin{equation}
%\label{}
h_i(\omega,x):=\cc_\omega(x,x+e_i)\Bigl[v_i\circ\tau_x(\omega)-\bigl(h_r^{\omega,v}(x+e_i)-h_r^{\omega,v}(x)\bigr)\Bigr]^2.
\end{equation}
 Given any~$z\in\Z^d$, the  definition of~$\mu_r^{\omega,v}$ gives
\begin{equation}
\label{E:4.39}
\begin{aligned}
\int_{\Omega\times\VV}\sum_{i=1}^d&\,\cc_{\omega'}( z,z+ e_i)\bigl[v_i\circ\tau_z(\omega')-u_i( z)\bigr]^2\,\mu_r^{\omega,v}(\textd\omega'\textd u)
\\
&=
\frac1{|\Lambda_{r+1}|}\,\sum_{x\in\Lambda_{r+1}}
\,\sum_{i=1}^d h_i(\omega,x+z)
\\
&\le
\frac1{|\Lambda_{r+1}|}\EE_{\Lambda_{r}}^\omega(v)+\frac1{|\Lambda_{r+1}|}\sum_{x\in\Lambda_{r+|z|_1+1}\,\smallsetminus\,\Lambda_{r}} 
\,\sum_{i=1}^d\cc_{\tau_x(\omega)}(0,e_i)\bigl|v_i\circ\tau_x(\omega)\bigr|^2.
\end{aligned}
\end{equation}
The H\"older inequality with parameters $\frac{1+q}q$ and~$1+q$ turns this into
\begin{equation}
\begin{aligned}
\label{E:4.40}
\int_{\Omega\times\VV}\sum_{i=1}^d\,\bigl|&v_i\circ\tau_z(\omega')-u_i( z)\bigr|^{\frac{2q}{1+q}}\,\mu_r^{\omega,v}(\textd\omega'\textd u)
\\
&\le
\biggl(\frac1{|\Lambda_{r+1}|}\sum_{(x,y)\in E(\Lambda_{r+|z|_1+1})}\cc_\omega(x,y)^{-q}\biggr)^{\frac{1}{1+q}}\Bigl(\text{r.h.s.\ of\ \eqref{E:4.39}}\Bigr)^{\frac{q}{1+q}},\end{aligned}
\end{equation}
where the sum in the first term is a bound on $\int_{\Omega\times\VV}\sum_{i=1}^d\cc_{\omega'}( z,z+ e_i)^{-q}\mu_r^{\omega,v}(\textd\omega'\textd u)$.

 Note that~$\omega\in\Omega_{p,q}$ shows that the first term on the right of \eqref{E:4.40} is bounded uniformly in~$r\ge1$. Along with boundedness of~$v$, Jensen's inequality shows that  the second term on the right of \eqref{E:4.39} vanishes as~$r\to\infty$. 
Since $\EE_\Lambda^\omega(v)\le\langle v\circ\tau,v\circ\tau\rangle_{\ell^2(\Lambda,\cc_\omega)}$, the argument underlying \eqref{E:4.29} shows that the right-hand side of \eqref{E:4.39}, and then also \eqref{E:4.40}, are bounded uniformly in~$r\ge1$. The triangle inequality and boundedness of~$v$ then give 
\begin{equation}
\label{E:4.41}
\sup_{r\ge1}\int_{\Omega\times\VV}\bigl| u( z)\bigr|_{\frac{2q}{1+q}}\, \mu_r^{\omega,v}(\textd\omega'\textd u)<\infty,
\end{equation}
where $|\cdot|_\alpha$ is the $\ell^\alpha$-norm on~$\R^d$.
Jointly with the tightness of the first marginals, which follows from Proposition~\ref{prop-2.2}, this implies the tightness claim in \eqref{E:4.36}.

Assume now that~$\mu$ is a weak limit of~$\mu_{r_k}^{\omega,v}$ for some sequence~$r_k\to\infty$. Since each~$\mu_{r_k}^{\omega,v}$ is a probability  measure,  the tightness extends this to~$\mu$. Translation invariance follows from  the averaging over shifts that show that $\mu_r^{\omega,v}(\tau_z^{-1}(A))-\mu_r^{\omega,v}(A)$ is order $r^{-1}$.  The first marginal of~$\mu$ coincides with~$\BbbP_\omega$ by Proposition~\ref{prop-2.2}. An elementary truncation (or Fatou's lemma for weak convergence) shows that the bound \eqref{E:4.41} survives the limit. As $\frac{2q}{1+q}\ge1$, we get the first half of \eqref{E:4.37}.  Since (as observed before) the second term on the right of \eqref{E:4.39} vanishes as~$r\to\infty$, the bound \eqref{E:4.38} follows from \eqref{E:4.39}. 

It remains to prove the second half of \eqref{E:4.37}. Since $q>1$ and \eqref{E:4.41} applies with $\frac{2q}{1+q}>1$,  the laws of the random variables $\{u_i(0)\colon i=1,\dots,d\}$ under  the measures $\{\mu_r^{\omega,v}\colon r\ge1\}$ are uniformly integrable. It thus suffices to show convergence in the mean,
\begin{equation}
\label{E:4.42}
\int_{\Omega\times\VV}u(0)\,\mu_r^{v,\omega}(\textd\omega'\textd u)\,\underset{r\to\infty}\longrightarrow\,0.
\end{equation}
 As it turns out,  the integral actually vanishes.  Indeed, letting let~$i\in\{1,\dots,d\}$, the~$i$-th  component of the integral in \eqref{E:4.42} equals
\begin{equation}
%\label{}
 \frac1{|\Lambda_{r+1}|}\sum_{x\in\Lambda_{r+1}}\bigl[h_r^{\omega,v}(x+e_i)-h_r^{\omega,v}(x)\bigr]
\end{equation}
 which vanishes thanks to the fact that $\Lambda_r$ is a square box and that~$h_r^{\omega,v}$ is constant outside~$\Lambda_r$.  This shows \eqref{E:4.42} and thus proves the claim.
\end{proofsect}

One additional property of measures~$\mu_k^{\omega,v}$ has not been used so far. Namely, the second variable in the indicator in \eqref{E:4.35} is a discrete gradient. %The second one is the observation that $h_r^{\omega,v}$ is a solution to $\cmss L_\omega h^{\omega,v}=\pi(0)^{-1}\nabla^\star\cc v$ in~$\Lambda_r$ with zero values outside~$\Lambda_r$. We collect the consequences of this in:
This gives:

\begin{corollary}
\label{cor-4.7}
Let $v\in\Cb(\Omega)$ and~$\omega\in\Omega^\star\cap\bigcup_{p,q>1}\Omega_{p,q}$ and let~$\mu$ be a subsequential weak limit of measures in \eqref{E:4.36}. Then the second marginal of~$\mu$ is supported on curl-free fields,
\begin{equation}
\label{E:4.44u}
\mu\biggl(\,\Omega\times\!\bigcap_{1\le i<j\le d}\Bigl\{u\in\VV\colon u_i(0)+u_j(e_i) = u_j(0)+u_i(e_j)\Bigr\}\biggr)=1.
\end{equation}
\end{corollary}

\begin{proofsect}{Proof}
As was just noted, the second marginal of~$\mu_r^{\omega,v}$ comes from a discrete gradient and so the equality holds for~$\mu_r^{\omega,v}$ in place of~$\mu$ regardless of~$r\ge1$. 
With the event closed in product topology, the claim follows by the Portmanteau Theorem.
\end{proofsect}

The fact that the second coordinate of~$\mu$ above is curl-free (a.k.a.\ a cocycle or a shift-covariant function) will permit us to invoke the following observation:

\begin{lemma}
\label{lemma-5.9}
Let~$\nu$ be a translation-invariant, ergodic probability measure on~$(\Omega,\FF)$ such that
\begin{equation}
\label{E:5.45i}
\cc(e)\in L^p(\nu)\,\,\wedge\,\, \cc(e)^{-1}\in L^q(\nu)
\end{equation}
with~$p=q=1$ in~$d=1,2$ and~$p,q\ge1$ satisfying $1/p+1/q<2/d$ in~$d\ge3$, for all~$e\in E(\Z^d)$. Then each measurable~$u\colon\Omega\to\R^d$ obeying the cocycle-condition
\begin{equation}
\label{E:5.46i}
u_i+u_j\circ\tau_{e_i} = u_j+u_i\circ\tau_{e_j},\quad\nu\text{\rm-a.s.}
\end{equation}
for all~$i,j=1,\dots d$ and the moment conditions
\begin{equation}
\label{E:5.47i}
\langle u,\cc u\rangle_{L^2(\nu,\R^d)}<\infty
\end{equation}
and
\begin{equation}
\label{E:5.48i}
u\in L^1(\nu)\,\,\wedge\,\, E_\nu(u)=0
\end{equation}
is a generalized gradient in the sense
\begin{equation}
\label{E:5.49i}
\inf_{f\in \Cb(\R)}\bigl\langle u-\nabla f,\cc(u-\nabla f)\bigr\rangle_{L^2(\nu,\R^d)}=0.
\end{equation}
\end{lemma}

\begin{proofsect}{Proof}
A standard argument based on the Parallelogram Law shows that any minimizing sequence $\{f_n\}_{n\ge1}$ of \eqref{E:5.49i} is Cauchy in weighted-$L^2$ sense,
\begin{equation}
%\label{}
\lim_{n,m\to\infty}\,\bigl\langle \nabla f_n-\nabla f_m,\cc(\nabla f_n-\nabla f_m)\bigr\rangle_{L^2(\nu,\R^d)}=0.
\end{equation}
Denoting by~$v$ the limit of $u-\nabla f_n$, the vector field~$v$ is readily checked to obey the analogues of \twoeqref{E:5.46i}{E:5.47i}. Using that~$\cc(e)^{-1}\in L^1(\nu)$, the Cauchy-Schwarz inequality then shows that the convergence also takes place in~$L^1(\nu)$ and so~$v$ obeys \eqref{E:5.48i}.

From~$v$ being a cocycle we get existence of~$h\colon\Omega\times\Z^d\to\R$ such that $h(\cdot,0)=0$ and
\begin{equation}
\label{E:5.51i}
h(\omega,x+e_i)-h(\omega,x) = v_i\circ\tau_x(\omega)
\end{equation}
 for all~$x\in\Z^d$, all~$i=1,\dots,d$ and~$\nu$-a.e.\ $\omega\in\Omega$. The fact that~$v$ arises as a minimizer in \eqref{E:5.49i} in turn implies that~$x\mapsto h(\omega,x)$ is a.s.\ $\cmss P_\omega$-harmonic. We claim that \eqref{E:5.48i} along with our assumptions on~$\nu$ then give
 \begin{equation}
\label{E:5.52i}
\lim_{|x|\to\infty}\frac{h(\cdot,x)}{|x|}=0,\quad\nu\text{\rm-a.s.}
\end{equation}
Indeed, $\cc(e)\in L^1(\nu)$ along with Birkhoff's Ergodic Theorem gives sublinearity of~$h$ along coordinate directions; see~\cite[Lemma~4.8]{B11}. This shows \eqref{E:5.52i} in~$d=1$; for~$d=2$ we invoke the argument from~\cite{BB07} (or the proof of~\cite[Theorem~4.7]{B11}) that relies on nearest-neighbor nature of~$\cmss P_\omega$. In~$d\ge3$ we instead plug in the elliptic-regularity theory developed in Andres, Deuschel and Slowik~\cite{ADS15}  and finessed for this particular purpose in Nguyen~\cite{Nguyen}. 

We now conclude by showing that~$h$ is constant $\nu$-a.s. Arguments for absence of sublinear harmonic functions for walks subject to heat-kernel bounds have been given before; e.g., in Benjamini, Duminil-Copin, Kozma and Yadin~\cite{BDCKY}, albeit under assumptions that do not generally apply here. We proceed by an argument underlying the proof of the Cacciopoli inequality. Let~$\eta\colon\Z^d\to[0,1]$ have finite support. The fact that~$h(\omega,\cdot)$ is~$\cmss P_\omega$-harmonic then gives
\begin{equation}
%\label{}
\sum_{(x,y)\in E(\Z^d)}\cc_\omega(x,y)\bigl[h(\omega,y)-h(\omega,x)\bigr]\bigl[h(\omega,y)\eta(y)-h(\omega,x)\eta(x)\bigr]=0.
\end{equation}
Writing the last term as $[h(\omega,y)-h(\omega,x)]\eta(y)+h(\omega,x)[\eta(y)-\eta(x)]$ then shows
\begin{equation}
\label{E:5.54i}
\begin{aligned}
\sum_{(x,y)\in E(\Z^d)}&\cc_\omega(x,y)\bigl[h(\omega,y)-h(\omega,x)\bigr]^2\eta(y)
\\&\le \sum_{(x,y)\in E(\Z^d)}\cc_\omega(x,y)\bigl|h(\omega,y)-h(\omega,x)\bigr|\bigl|h(\omega,x)\bigr|\bigl|\eta(y)-\eta(x)\bigr|.
\end{aligned}
\end{equation}
Now take~$\eta$ so that~$\eta=1$ on~$\Lambda_r$, $\eta=0$ on $\Z^d\smallsetminus\Lambda_{2r}$ and~$\eta$ decreases linearly in~$\Lambda_{2r}\smallsetminus\Lambda_r$. Then $|\eta(y)-\eta(x)|\le 1/r$ for~$(x,y)\in\Lambda_{2r}$ and, given~$\epsilon>0$, the fact that $|h(\omega,x)|\le \epsilon r$ on~$\Lambda_{2r}$ once~$r$ is large enough thanks to \eqref{E:5.52i} gives
\begin{equation}
%\label{}
\sum_{(x,y)\in E(\Lambda_r)}\cc_\omega(x,y)\bigl[h(\omega,y)-h(\omega,x)\bigr]^2
\le \epsilon \sum_{(x,y)\in E(\Lambda_{2r})}\cc_\omega(x,y)\bigl|h(\omega,y)-h(\omega,x)\bigr|.
\end{equation}
Dividing by~$|\Lambda_r|$ and taking $r\to\infty$ with the help of the Spatial Ergodic Theorem and \eqref{E:5.51i}, and also the Cauchy-Schwarz inequality on the right-hand side, yields
\begin{equation}
%\label{}
\langle v,\cc v\rangle_{L^2(\nu,\R^d)}
\le 2^d\epsilon \sqrt{\langle v,\cc v\rangle_{L^2(\nu,\R^d)}\sum_{i=1}^d E_\nu(\cc(e_i))}.
\end{equation}
As~$\cc$ is positive~$\nu$-a.s., taking~$\epsilon\downarrow0$ gives~$v=0$ a.s.\ proving \eqref{E:5.49i}.
\end{proofsect}

\begin{remark}
\label{rem-5.10}
The fact that square-integrable mean-zero cocycles can be approximated by gradients was shown under the assumption of uniform ellipticity in Biskup and Spohn~\cite[Theorem~5.4]{BS11}. The construction of an approximate gradient in Biskup and Rodriguez~\cite[Lemma 4.4]{BR18} (which also corrects for a tacit assumption of separate ergodicity in \cite{BS11}) extends this to convergence in~$L^1$ under the minimal moment conditions $\cc(e),\cc(e)^{-1}\in L^1(\nu)$. However, that construction does not seem to help for weighted-$L^2$-convergence needed here.

Our proof of Lemma~\ref{lemma-5.9} links the desired claim to absence of non-constant sublinear harmonic functions which, in $d\ge3$, require elliptic regularity theory and stricter moment assumptions than just the minimal ones. This approach gives sublinearity in the maximal sense which may be too much to ask as \eqref{E:5.54i} only requires sublinearity under averaging. It is an interesting question whether the claim holds under the minimal moment conditions $\cc(e),\cc(e)^{-1}\in L^1(\nu)$ in all~$d\ge1$.
\end{remark}

We are  now  ready to give:

\begin{proofsect}{Proof of Proposition~\ref{prop-4.4}}
Suppose that $\liminf_{r\to\infty}\frac1{|\Lambda_r|}\EE_{\Lambda_r}^\omega(v)$ is achieved along a sequence $r_k\to\infty$. Reducing to a subsequence if necessary, we may assume that~$\mu_{r_k}$ tends weakly to a probability measure~$\mu$. Consider the $\R^d$-valued random variable~$U$ defined formally as the map $U\colon\Omega\times\VV\to\R^d$ by $U(\omega',u'):=u'(0)$. Then the first half of \eqref{E:4.37} shows that $U\in L^1(\mu)$ and so we may set
\begin{equation}
%\label{}
u(\omega):=E_\mu(U|\FF)(\omega).
\end{equation}
 Part~(3) of Lemma~\ref{lemma-5.7i} along with the  second half of \eqref{E:4.37} then gives
\begin{equation}
%\label{}
E_{\BbbP_\omega}(u)=0.
\end{equation}
By conditional Jensen's inequality applied to \eqref{E:4.38}, we also have
\begin{equation}
\label{E:4.48}
\bigl\langle (v-u),\cc (v-u)\bigr\rangle_{L^2(\BbbP_\omega,\R^d)}\le\liminf_{r\to\infty}\frac1{|\Lambda_r|}\EE_{\Lambda_r}^\omega(v).
\end{equation}
 This in particular gives 
\begin{equation}
%\label{}
\langle u,\cc u\rangle_{L^2(\BbbP_\omega,\R^d)}<\infty
\end{equation}
using Lemma~\ref{lemma-4.5} along with the boundedness of~$v$ and $\cc\in L^1(\BbbP_\omega)$. 

The identity in Corollary~\ref{cor-4.7} now implies that~$u$ is (a.s.) curl-free (a.k.a~{cocycle}), meaning that for all
$i,j=1,\dots,d$,
\begin{equation}
%\label{}
u_i+u_j\circ\tau_{e_i} = u_j+u_i\circ\tau_{e_j},\quad\BbbP_\omega\text{-a.s.}
\end{equation}
Since~$\BbbP_\omega$ is ergodic,  \eqref{E:5.49i} along with \eqref{E:4.48} give 
\begin{equation}
\label{E:4.50}
\inf_{f\in\Cb(\Omega)}\bigl\langle (v-\nabla f),\cc (v-\nabla f)\bigr\rangle_{L^2(\BbbP_\omega,\R^d)}\le\liminf_{r\to\infty}\frac1{|\Lambda_r|}\EE_{\Lambda_r}^\omega(v).
\end{equation}
In conjunction with the upper bound from Lemma~\ref{lemma-4.5}, this yields the claim.
\end{proofsect}

\begin{remark}
Since equality holds in \eqref{E:4.50}, it also does in \eqref{E:4.48}. Using the equality part of Jensen's inequality, we thus have $U=u$ $\mu$-a.s. By the aforementioned uniqueness of~$u$ also~$\mu$ is unique and so~$\mu_n\overset{\text{\rm w}}\longrightarrow\mu$.
\end{remark}

%\newpage
\section{Proof of the CLT}
\label{sec6}\noindent
We are now ready to move to the proof of Theorem~\ref{thm-2.4}. Thanks to the tightness established in Proposition~\ref{prop-2.7}, an IIP has been reduced to convergence of finite-dimensional distributions; i.e., a multivariate Central Limit Theorem. Before delving into the proof of the latter, let us review the salient steps of Kipnis and Varadhan's theory. 

\subsection{Review of Kipnis-Varadhan's theory}
Fix~$\omega\in\Omega$ and let~$X$ be a path of the random walk sampled from the measure $P_\omega^x$. Recall the rewrite \eqref{E:3.2u} of~$X_n-X_0$ into a martingale plus the sum
\begin{equation}
\label{E:5.1o}
\sum_{k=0}^{n-1}V\circ\tau_{X_k}(\omega),
\end{equation}
where~$V$ denotes the local drift from \eqref{E:3.2r}. The martingale has bounded increments and so its asymptotic diffusive scaling can presumably be handled by the Martingale Functional CLT that we invoked already in the one-dimensional case. The key challenge is to extract a similar conclusion for the sum \eqref{E:5.1o}.

Kipnis and Varadhan~\cite{KV86} approached this by developing a fluctuation theory of additive functionals of Markov chains. (In the case at hand, the Markov chain is that of ``point of view of the particle'' with state space~$\Omega$, transition kernel~$\Pi$ from \eqref{E:3.41q} and generator~$\LL$ in \eqref{E:4.1w}.)
An early work on this (likely unknown to Kipnis and Varadhan) is due to Gordin~\cite{G69} whose approach would apply if we had~$V=-\LL\chi$ for some bounded continuous~$\chi$. Unfortunately, this is likely false and definitely too hard to establish directly even for uniformly elliptic conductances. Kipnis and Varadhan therefore resorted to approximation. 

Define $\chi_\epsilon\colon\Omega\to\R^d$ by
\begin{equation}
\label{E:5.2r}
\chi_\epsilon(\omega):=\sum_{n\ge0}\frac1{(1+\epsilon)^n}(\Pi^n V)(\omega).
\end{equation}
Since~$V\in\Cb(\Omega,\R^d)$ and~$\Pi$ is an operator of norm one, the sum converges  uniformly  and so~$\chi_\epsilon\in \Cb(\Omega,\R^d)$. Moreover, $\chi_\epsilon$  obeys the ``massive'' Poisson equation
\begin{equation}
\label{E:3.7w}
(\epsilon-\LL)\chi_\epsilon = V
\end{equation}
making it an approximate inverse of~$V$ under~$-\LL$. Writing, as before, $\FF_k:=\sigma(X_0,\dots,X_k)$ and substituting \eqref{E:3.7w} into the sum in \eqref{E:5.1o}, we get
\begin{equation}
\begin{aligned}
\label{E:3.8r}
\quad
\sum_{k=0}^{n-1}V\circ\tau_{X_k}(\omega&) = \chi_\epsilon(\omega)-\chi_\epsilon\circ\tau_{X_n}(\omega)
\\
&+\sum_{k=0}^{n-1}\epsilon\chi_\epsilon\circ\tau_{X_k}(\omega)
+\sum_{k=1}^{n}\Bigl[\chi_\epsilon\circ\tau_{X_k}(\omega)
-E_\omega^0\bigl(\chi_\epsilon\circ\tau_{X_k}(\omega)
\,\big|\,\FF_{k-1}\bigr)\Bigr]
\end{aligned}
\end{equation}
using the same manipulations that led to \eqref{E:3.2u}.

In order to plug \eqref{E:3.8r} in \eqref{E:3.2u}, it is convenient to work with the shorthand 
\begin{equation}
\label{E:5.5i}
\psi_\epsilon(\omega,x):=x+\chi_\epsilon\circ\tau_x(\omega)-\chi_\epsilon(\omega).
\end{equation}
The definition implies the so called cocycle property
\begin{equation}
\label{E:5.6i}
\forall x,y\in\Z^d\colon\,\,\psi_\epsilon\bigl(\tau_x(\omega),y\bigr)=\psi_\epsilon(\omega,x+y)-\psi_\epsilon(\omega,x)
\end{equation}
which permits us to summarize the above as
\begin{equation}
\label{E:4.8u}
X_n = X_0+\bigl[\chi_\epsilon(\omega)-\chi_\epsilon\circ\tau_{X_n}(\omega)\bigr]+
\biggl[\,\sum_{k=0}^{n-1}\epsilon\chi_\epsilon\circ\tau_{X_k}(\omega)\biggr]+M_n^{(\epsilon)},
\end{equation}
where
\begin{equation}
\label{E:4.9u}
M_n^{(\epsilon)}:=\sum_{k=1}^{n}\biggl[\psi_\epsilon\bigl(\tau_{X_{k-1}}(\omega),X_k-X_{k-1}\bigr)-E_\omega^0\Bigl(\psi_\epsilon\bigl(\tau_{X_{k-1}}(\omega),X_k-X_{k-1}\bigr)\,\Big|\,\FF_{k-1}\Bigr)\biggr].
\end{equation}
As before,~$M_n^{(\epsilon)}$ is a martingale with bounded increments.

At the first sight, \eqref{E:4.8u}  may  seem to be no better than \eqref{E:3.2u} because we still have an additive functional of the Markov chain on the right-hand side. However, unlike \eqref{E:3.2u}, we now have a parameter to play with. Kipnis and Varadhan~\cite{KV86} utilized this by taking the limit~$\epsilon\downarrow0$ and showing that
\begin{equation}
\label{E:5.9i}
\epsilon\chi_\epsilon\,\,\underset{\epsilon\downarrow0}{\overset{L^2}\longrightarrow}\,\,0
\quad\wedge\quad\chi_\epsilon\circ\tau_z-\chi_\epsilon\,\,\underset{\epsilon\downarrow0}{\overset{L^2}\longrightarrow}\,\,\chi(\cdot,z),
\end{equation}
where the~$L^2$ limits are relative to the underlying probability space on environments weighted by the conductances in the latter case and~$\chi(\cdot, z)$ is a function on~$\Omega$ defined hereby. The  additive term  on the right of \eqref{E:4.8u} then disappears and we get
\begin{equation}
%\label{}
X_n = X_0-\chi(\omega,X_n)+M_n,
\end{equation}
where~$M_n$ is obtained by replacing~$\psi_\epsilon(\omega,x)$ by~$\psi(\omega,x):=x+\chi(\omega,x)$ in \eqref{E:4.9u}. Unfortunately,  these manipulations come with a price tag:  $n\mapsto\chi(\omega,X_n)$ is no  longer  bounded and  the  proof of a CLT-scaling requires control of its growth. See Section~\ref{sec-2.1}  for a discussion of contexts where this has been achieved. 

%%%%
\begin{comment}
Unfortunately, the reduction offered by the limit comes with a cost: $n\mapsto\chi(\omega,X_n)$ is no nonger bounded and a proof of a CLT-scaling requires showing that
\begin{equation}
%\label{}
\frac1{\sqrt n}\,\chi(\omega,X_n)\,\underset{n\to\infty}\longrightarrow\,0\quad\text{in probability}.
\end{equation}
In~\cite{KV86}, Kipnis and Varadhan managed to prove this under the joint law of the walk and the environment capitalizing on the fact that, under this law, the Markov chain from the ``point of view of the particle'' is reversible. Subsequent work by many authors (see the Introduction) pushed the same by showing instead that $\chi(\omega,x)=o(|x|)$ as~$|x|\to\infty$ for almost every environment~$\omega$. However, these proofs rely heavily on measure-theoretic arguments and there seems to be little hope that they could apply for a fixed~$\omega$.
\end{comment}
%%%%%

\subsection{The new idea}
The novelty of our approach is that we keep $\epsilon$ fixed and, instead of worrying about the growth of the corrector, we argue that,  for~$\epsilon$ small, the martingale $M_n^{(\epsilon)}$ already contains most of the variance of the walk. This boils down to showing that,  under the diffusive scaling of space and time, $\sum_{k=0}^{n-1}\epsilon\chi_\epsilon\circ\tau_{X_k}(\omega)$ is small with high probability once~$\epsilon$ is small.  We will achieve this  by another martingale approximation but this time localized to a finite box of size~$r$ that will be large on a diffusive scale. 

Recall the notation~$\cmss L_\omega$ from \eqref{E:4.8w} for the generator of the random walk in the physical space. Given~$r\ge1$,~$\epsilon>0$,~$\omega\in\Omega$ and~$\chi_\epsilon$ as in \eqref{E:5.2r}, let $\theta_{r,\epsilon}\colon\Z^d\to\R^d$ be the solution of
\begin{equation}
\label{E:3.24}
 \begin{alignedat}{2}
-\cmss L_\omega\theta_{r,\epsilon}(\omega,x) &= \epsilon\chi_\epsilon\circ\tau_x(\omega)\qquad & \text{if }x\in\Lambda_r,
\\
\theta_{r,\epsilon}(\omega,x)&=0\qquad & \text{if } x\not\in\Lambda_r.
\end{alignedat}
\end{equation}
This is uniquely solvable because, thanks to all the conductances being positive,~$\cmss L_\omega$ is invertible on the space of functions that vanish outside~$\Lambda_r$. Recalling our notation~$\ttau_{\Lambda}$ for the first exit time of~$X$ from~$\Lambda$, the same calculations as before show
\begin{equation}
\label{E:5.12}
\sum_{k=1}^{n\wedge\ttau_{\Lambda_r}}\epsilon\chi_\epsilon\circ\tau_{X_{k-1}}(\omega) = \theta_{r,\epsilon}(\omega,X_0)-\theta_{r,\epsilon}(\omega,X_{n\wedge\ttau_{\Lambda_r}})+\wt M_n^{(r,\epsilon)},
\end{equation}
where
\begin{equation}
\label{E:5.14}
\wt M_n^{(r,\epsilon)}:=\sum_{k=1}^{n\wedge\ttau_{\Lambda_r}}\Bigl(\theta_{r,\epsilon}(\omega,X_k)%-\theta_{r,\epsilon}(\omega,X_{k-1})
-E_\omega^0\bigl(\theta_{r,\epsilon}(\omega,X_k)
%-\theta_{r,\epsilon}(\omega,X_{k-1})
\big|\,\FF_{k-1}\bigr)\Bigr).
\end{equation}
Note that, since we are not using the ``massive'' Poisson equation to define~$\theta_{r,\epsilon}$, we are not just repeating the step  leading  from \eqref{E:3.2u} to \eqref{E:4.8u}. In particular, no additive functional of the Markov chain arises on the right of \eqref{E:5.12}.

%%%%
\begin{comment}
\begin{remark}
Incidentally, the function
\begin{equation}
%\label{}
\psi_{r,\epsilon}(\omega,x):=x+\chi_\epsilon\circ\tau_x(\omega) +\theta_{r,\epsilon}(\omega,x)
\end{equation}
then obeys
\begin{equation}
%\label{}
L_\omega\psi_{r,\epsilon}(\omega,x) = V\circ\tau_x(\omega)+\bigl(\epsilon\chi_\epsilon-V)\circ\tau_x(\omega)-\epsilon\chi_\epsilon\circ\tau_x(\omega)=0
\end{equation}
at each $x\in\Lambda_r$ meaning that~$x\mapsto\psi_{r,\epsilon}(\omega,x)$ is $L_\omega$-harmonic and can thus be used as a harmonic coordinate inside~$\Lambda_r$. Writing
\begin{equation}
%\label{}
\tau_r:=\inf\bigl\{k\ge0\colon X_k\not\in\Lambda_r\bigr\}
\end{equation}
for the first exit time from~$\Lambda_r$, the harmonicity translates into $\{\psi_{r,\epsilon}(\omega,X_{n\wedge\tau_r})\colon n\ge0\}$ being a martingale with respect to the filtration~$\FF_k:=\sigma(X_0,\dots,X_k)$.
\end{remark}
\end{comment}
%%%%% 
%

Besides the invertibility of~$\cmss L_\omega$, another reason for using zero boundary conditions is that the terms on the right of \eqref{E:5.12} can be controlled using the Dirichlet form associated with~$\theta_{r,\epsilon}$. The technical tool needed for the first term is the Sobolev inequality (which works best under zero boundary conditions) combined with heat-kernel estimates. As these are formulated for continuous time, we will have to work with that as well.

\begin{lemma}
\label{lemma-5.1}
Let $p,q>1$ obey \eqref{E:2.6} and let~$t\mapsto N(t)$ be the rate-1 Poisson process independent of~$X$. For all~$\omega\in\Omega_{p,q}$ there is~$c_1'(\omega)\in(0,\infty)$ such that for all~$\delta>0$,~$r\ge1$ and~$t\ge1$,
\begin{equation}
%\label{}
\begin{aligned}
P_\omega^0\Bigl(\bigl|\theta_{r,\epsilon}(\omega,X_{N(t)})&\bigr|_2>\delta r\,\wedge\,\ttau_{\Lambda_r}>N(t)\Bigr)
\\&\le c_1'(\omega)\frac{r^d}{t^{d/2}}\left[\frac1\delta\sqrt{\frac{\langle\nabla\theta_{r,\epsilon},\,\nabla\theta_{r,\epsilon}\rangle_{\ell^2(E(\Lambda_r),\cc)}}{r^d}}\,\right]^{1-1/p},
\end{aligned}
\end{equation}
where $|\cdot|_2$ is the Euclidean norm in~$\R^d$ and where we have suppressed the argument~$\omega$ of~$\theta_{r,\epsilon}$ in the inner product (which is defined in \eqref{E:4.25w}; see also Remark~\ref{rem-4.3}).
\end{lemma}

\begin{proofsect}{Proof}
Fix $p,q>1$ so that \eqref{E:2.6} holds  and let~$\omega\in\Omega_{p,q}$. 
Recalling our notation~$Y_t$ for the process~$X_{N(t)}$ and writing $\ttau_{\Lambda}(Y)$ for the exit time thereof, for each~$\lambda>0$ the probability in the statement is  at most 
\begin{equation}
\begin{aligned}
P_\omega^0\Bigl(\pi_\omega(Y_t)>&\lambda\,\wedge\,\ttau_{\Lambda_r}(Y)>t\Bigr)
\\&+
P_\omega^0\Bigl(\,\bigl|\theta_{r,\epsilon}(\omega,Y_t)\bigr|_2>\delta r\,\wedge\,\pi_\omega(Y_t)\le\lambda\,\wedge\,\ttau_{\Lambda_r}(Y)>t\Bigr).
\end{aligned}
\end{equation}
The Markov inequality  dominates  this further by
\begin{equation}
\label{E:5.17}
\frac1{\lambda^{p-1}}\sum_{x\in\Lambda_r}P_\omega^0(Y_t=x)\pi_\omega(x)^{p-1}+\frac1\delta\sum_{x\in\Lambda_r}P_\omega^0(Y_t=x)1_{\{\pi_\omega(x)\le\lambda\}}
\frac{|\theta_{r,\epsilon}(\omega,x)|_2}r.
\end{equation}
Proposition~\ref{prop-ADS} gives $P_\omega^0(Y_t=x)\le c_1(\omega) t^{-d/2}\pi_\omega(x)$ for some~$c_1$ depending on~$\omega$ and parameters~$p,q$.  Invoking also Jensen's inequality,  the containment~$\omega\in\Omega_{p,q}$  allows us bound  the first sum in \eqref{E:5.17} by~$c'(\omega)r^d/t^{d/2}$ uniformly in~$r\ge1$ and~$t\ge1$.

For the second sum the heat-kernel bound combined with the restriction on~$\pi_\omega(x)$ implies $P_\omega^0(Y_t=x)\le c_1(\omega) \lambda t^{-d/2}$ for each term under the sum. The $\ell^1$-Sobolev inequality (see, e.g., \cite[Lemma 2.1]{BR18}) in turn gives 
\begin{equation}
%\label{}
\sum_{x\in\Lambda_r}\frac{|\theta_{r,\epsilon}(\omega,x)|_2}r\le c(d) \sum_{(x,y)\in E(\Lambda_r)}\bigl|\,\theta_{r,\epsilon}(x,\omega)-\theta_{r,\epsilon}(y,\omega)\bigr|_2,
\end{equation}
where~$c(d)$ is a $d$-dependent constant and where the triangle inequality was used to  write the result using the Euclidean norm on~$\R^d$.  Invoking the Cauchy-Schwarz inequality, the sum  on the right  is further bounded by
\begin{equation}
%\label{}
\biggl(\,\sum_{(x,y)\in E(\Lambda_r)}\cc_\omega(x,y)^{-1}\biggr)^{1/2}\langle\nabla\theta_{r,\epsilon},\,\nabla\theta_{r,\epsilon}\rangle_{\ell^2(E(\Lambda_r),\cc)}^{1/2}.
\end{equation}
Since~$\omega\in\Omega_{p,q}$ with~$q\ge1$, the first term is at most an $\omega$-dependent constant times~$r^{d/2}$, uniformly in~$r\ge1$.

Combining the above estimates, the quantity in \eqref{E:5.17} is at most
\begin{equation}
%\label{}
c''(\omega)\frac{r^d}{t^{d/2}}\left[\,\frac1{\lambda^{p-1}}+\frac{\lambda}{\delta}\sqrt{\frac{\langle\nabla\theta_{r,\epsilon},\,\nabla\theta_{r,\epsilon}\rangle_{\ell^2(E(\Lambda_r),\cc)}}{r^{d}}}\,\right],
\end{equation}
where~$c''(\omega)$ dominates various  $\omega$-dependent  constants mentioned above. The claim follows by noting that $\inf_{\lambda>0}[\lambda^{1-p}+\lambda a]\le 2 a^{1-1/p}$.
\end{proofsect}

\begin{lemma}
\label{lemma-5.2}
Let $p,q>1$ obey \eqref{E:2.6} and let~$t\mapsto N(t)$ be a rate-$1$ Poisson process independent of~$X$. For all~$\omega\in\Omega_{p,q}$ there exists~$c_2'(\omega)\in(0,\infty)$ such that for all $r\ge1$,  $\epsilon>0$ and  $0< s\le t$, 
\begin{equation}
%\label{}
E_\omega^0\Bigl(\bigl|\wt M_{N(t)}^{(r,\epsilon)}-\wt M_{N(s)}^{(r,\epsilon)}\bigr|_2^2\Bigr) \le c_2'(\omega)\frac{t}{s^{d/2}}\langle\nabla\theta_{r,\epsilon},\,\nabla\theta_{r,\epsilon}\rangle_{\ell^2(E(\Lambda_r),\cc)}.
\end{equation}
 Here $\wt M_{n}^{(r,\epsilon)}$ is as in \eqref{E:5.14}. 
\end{lemma}

\begin{proofsect}{Proof}
Let us drop the superscript $(r,\epsilon)$ from~$\wt M$ for the duration of this proof. Note that $t\mapsto \wt M_{N(t)}$ is a  $\R^d$-valued  continuous-time martingale with quadratic variation
\begin{equation}
%\label{}
\langle \wt M_{N(\cdot)}\rangle_t = \int_0^t D(X_{N(u)})1_{\{\ttau_{\Lambda_r}>N(u)\}}\textd u,
\end{equation}
where
\begin{equation}
\label{E:5.23}
D(x):= E_\omega^x\biggl(\,\Bigl|\theta_{r,\epsilon}(\omega,X_1)-E_\omega^0\bigl(\theta_{r,\epsilon}(\omega,X_1)\bigr)
\Bigr|^2_2\biggr).
\end{equation}
Hence we get
\begin{equation}
\label{E:5.24}
\begin{aligned}
E_\omega^0\Bigl(\bigl|\wt M_{N(t)}-\wt M_{N(s)}\bigr|_2^2\Bigr)
&=\int_s^tE_\omega^0\bigl(D(Y_u)1_{\{\ttau_{\Lambda_r}(Y)>u\}}\bigr)\textd u
\\
&\le\int_s^t\sum_{x\in\Lambda_r}D(x)P^0_\omega(Y_u=x)\textd u,
\end{aligned}
\end{equation}
where we dropped the event $\ttau_{\Lambda_r}(Y)>u$ after restricting the sum to~$x\in\Lambda_r$.

We will now estimate the integrand in \eqref{E:5.24}. First, subtracting $\theta_{r,\epsilon}(\omega,x)$ from both occurrences of $\theta_{r,\epsilon}(\omega,X_1)$ in \eqref{E:5.23}, the inequality $(a+b)^2\le 2a^2+2b^2$ shows
\begin{equation}
%\label{}
D(x)\le 4\!\sum_{y\colon (x,y)\in E(\Lambda_r)}\!\!\!\cmss P_\omega(x,y)\bigl|\theta_{r,\epsilon}(\omega,y)-\theta_{r,\epsilon}(\omega,x)\bigr|_2^2.
\end{equation}
The heat-kernel bound $P_\omega^0(Y_u=x)\le c_1(\omega) u^{-d/2}\pi_\omega(x)$ from Proposition~\ref{prop-ADS} along with the fact that~$u\ge s$ throughout the integration domain then dominates \eqref{E:5.24} by
\begin{equation}
%\label{}
4c_1(\omega)\frac{t}{s^{d/2}}
\sum_{x\in\Lambda_r}\pi_\omega(x)\!\!\!\sum_{y\colon (x,y)\in E(\Lambda_r)}\!\!\!\cmss P_\omega(x,y)\bigl|\theta_{r,\epsilon}(\omega,y)-\theta_{r,\epsilon}(\omega,x)\bigr|_2^2.
\end{equation}
The claim now follows  by noting  that~$\pi_\omega(x)\cmss P_\omega(x,y)=\cc_{\omega}(x,y)$.
\end{proofsect}

We can now combine the above bounds into:

\begin{lemma}
\label{lemma-5.3}
Let $p,q>1$ obey \eqref{E:2.6} and let~$t\mapsto N(t)$ be a rate-$1$ Poisson process independent of~$X$. Then for all~$\omega\in\Omega^\star\cap\Omega_{p,q}$,
\begin{equation}
\label{E:5.20}
\lim_{\delta\downarrow0}\,\,\limsup_{\epsilon\downarrow0}\,\limsup_{t\to\infty}\,
P_\omega^0\Biggl(\,\biggl|\,\sum_{k=N(\delta t)+1}^{N(t)}\epsilon\chi_\epsilon\circ\tau_{X_{k-1}}\biggr|_2>\delta\sqrt t\Biggr)
=0.
\end{equation}
\end{lemma}

\begin{proofsect}{Proof}
Fix~$\delta\in(0,1)$ small and set $r:=\sqrt{t/\delta}$. (We will however continue omitting the superscript $(r,\epsilon)$ from~$\wt M$.) In light of the representation \eqref{E:5.12}, the probability in the statement is bounded by
\begin{equation}
\begin{aligned}
\label{E:5.28}
\qquad
P_\omega^0\bigl(\ttau_{\Lambda_r}\le N(t)\bigr)
+P_\omega^0\Bigl(&|\theta_{r,\epsilon}(X_{N(\delta t)})|_2>\tfrac13\delta\sqrt t\,\wedge\,\ttau_{\Lambda_r}>N(\delta t)\Bigr)\qquad\qquad
\\
&+P_\omega^0\Bigl(|\theta_{r,\epsilon}(X_{N(t)})|_2>\tfrac13\delta\sqrt t\,\wedge\,\ttau_{\Lambda_r}>N(t)\Bigr)
\\
&\qquad\qquad+P_\omega^0\Bigl(|\wt M_{N(t)}-\wt M_{N(\delta t)}|_2>\tfrac13\delta\sqrt t\Bigr).
\end{aligned}
\end{equation}
Abbreviating
\begin{equation}
%\label{}
Q_{r,\epsilon}:=\frac{1}{r^d}\langle\nabla\theta_{r,\epsilon},\,\nabla\theta_{r,\epsilon}\rangle_{\ell^2(E(\Lambda_r),\cc)},
\end{equation}
Proposition~\ref{prop-4.1} (with~$\alpha:=1$ and $\sigma:=1$) and Lemmas~\ref{lemma-5.1}--\ref{lemma-5.2} along with the Chebyshev inequality bound \eqref{E:5.28} by  $P(N(t)>2t)$, which decays exponentially as~$t\to\infty$, plus the quantity 
\begin{equation}
%\label{}
2\tilde c_1\frac t{r^2}
+2c_1'\frac{r^d}{t^{d/2}}\frac1{\delta^{d/2}}\bigl[3\delta^{-1/2}\sqrt{Q_{r,\epsilon}}\,\bigr]^{1-1/p}
+c_2'\frac 9{\delta^2 t}\frac{t}{(\delta t)^{d/2}}r^d Q_{r,\epsilon},
\end{equation}
provided that~$r\ge R_0$. (We have suppressed the $\omega$-dependence of~$\tilde c_1,c_1'$ and~$c_2'$.)
Plugging in for~$r$  in the prefactors,  this further simplifies into
\begin{equation}
\label{E:5.31}
2\tilde c_1\delta
+6c_1'\delta^{-(d+\frac12(1-1/p))}Q_{r,\epsilon}^{\frac12(1-1/p)}
+9c_2'\delta^{-(d+2)}Q_{r,\epsilon}.
\end{equation}
Our strategy is to prove that~$Q_{r,\epsilon}\to0$ as~$r\to\infty$ followed by~$\epsilon\downarrow0$. 

Note that the definition of~$\theta_{r,\epsilon}$ gives
\begin{equation}
%\label{}
\begin{aligned}
\langle\nabla\theta_{r,\epsilon},\,\nabla\theta_{r,\epsilon}\rangle_{\ell^2(E(\Lambda_r),\cc)}
&=\langle\theta_{r,\epsilon},\,(-\cmss L_\omega)\theta_{r,\epsilon}\rangle_{\ell^2(\Lambda_r,\pi_\omega)}
\\
&=\langle\theta_{r,\epsilon},\,\epsilon\chi_\epsilon\circ\tau\rangle_{\ell^2(\Lambda_r,\pi_\omega)}
\\
&=\bigl\langle\epsilon\chi_\epsilon\circ\tau,\,(-\cmss L_\omega)_{\Lambda_r}^{-1}\epsilon\chi_\epsilon\circ\tau\bigr\rangle_{\ell^2(\Lambda_r,\pi_\omega)}.
\end{aligned}
\end{equation}
Since $e_i\cdot(\epsilon\chi_\epsilon) = \nabla^\star\cc v$ for $v(\omega'):=e_i+\nabla(e_i\cdot\chi_\epsilon)(\omega')$,  applying Theorem~\ref{thm-3.1} coordinate by coordinate  shows
\begin{equation}
%\label{}
Q_{r,\epsilon}\,\underset{r\to\infty}\longrightarrow\,\bigl\langle\epsilon\chi_\epsilon,(-\LL)^{-1}\epsilon\chi_\epsilon\bigr\rangle_{L^2(\Q_\omega)},
\end{equation}
where the inner product  includes  contraction of the indices (via dot product) of the vector valued functions on both sides.
 With the limit object expressed as integral over~$\Omega$, we are now in a position to invoke standard calculations from stochastic homogenization:  Writing $\mu_V$ for the spectral measure of~$-\LL$ associated with (vector valued) function~$V$, the definition of~$\chi_\epsilon$ in \eqref{E:3.7w} yields
\begin{equation}
\label{E:5.34}
\bigl\langle\epsilon\chi_\epsilon,(-\LL)^{-1}\epsilon\chi_\epsilon\bigr\rangle_{L^2(\Q_\omega)} = 
\int_{[0,2]}\Bigl(\frac{\epsilon}{\epsilon+\lambda}\Bigr)^2\frac1\lambda\,\mu_V(\textd\lambda).
\end{equation}
Lemma~\ref{lemma-3.2} and the fact that~$e_i\cdot V=\pi(0)^{-1}\nabla^\star(\cc \tilde v)$ for $\tilde v(\cdot):=e_i$ along with the Spectral Theorem show $\int\lambda^{-1}\mu_V(\textd\lambda)<\infty$. The integral in \eqref{E:5.34} then vanishes in the limit~$\epsilon\downarrow0$ by the Dominated Convergence Theorem. Using this in \eqref{E:5.31}, the \emph{limes superior} as~$\epsilon\downarrow0$ in \eqref{E:5.20} is at most $2\tilde c_1\delta$. Taking~$\delta\downarrow0$ gives the claim.
\end{proofsect}

We are now finally ready to give:

\begin{proofsect}{Proof of Theorem~\ref{thm-2.4}}
The stated properties of~$\Omega^\star$ have already been shown in Proposition~\ref{prop-3.1}, so all we have to do is to prove an IIP. The one-dimensional case has been settled in Section~\ref{sec-3.2} so we may assume~$d\ge2$. Here Proposition~\ref{prop-2.7} gives tightness and so it suffices to prove convergence of finite-dimensional distributions. This amounts to proving that for all~$m\ge1$ natural, all~$0< t_0<t_1<\dots<t_m$ real and all~$\omega\in\Omega^\star\cap\Omega_{p,q}$ for some $p,q>1$ satisfying \eqref{E:2.6}, 
\begin{equation}
\label{E:5.35}
\text{\rm the law of }\Bigl(B^{(n)}_{t_1}-B^{(n)}_{t_{0}},\dots,B^{(n)}_{t_m}-B^{(n)}_{t_{m-1}}\Bigr)\text{\rm\ on~$\R^d$ under }P_\omega^0
\end{equation}
tends, as $n\to\infty$, weakly to the law of independent $d$-dimensional normals
\begin{equation}
\label{E:5.36}
\Bigl(\NN\bigl(0,(t_1-t_0)\Sigma\bigr),\dots\NN\bigl(0,(t_{m}-t_{m-1})\Sigma\bigr)\Bigr)
\end{equation}
for some~($\omega$-dependent) covariance~$\Sigma$. Here $t_0>0$ may be assumed thanks to the established tightness of~$B^{(n)}$ and the fact that~$B^{(n)}_0=0$ $P^0_\omega$-a.s.

Let~$Y_t:=X_{N(t)}$ denote, for~$N$ being a rate-$1$ Poisson process independent of~$X$, the continuous time version of~$X$. We start by noting that the tightness of processes $\{B^{(n)}\colon n\ge1\}$ along with $N(tn)/n\to t$ in probability shows
\begin{equation}
%\label{}
B_t^{(n)}-n^{-1/2}Y_{nt}\,\underset{n\to\infty}\longrightarrow\,0\quad\text{in $P_\omega^0$-probability}
\end{equation}
and so it suffices to prove \twoeqref{E:5.35}{E:5.36} with~$t\mapsto B^{(n)}_t$ replaced by $t\mapsto n^{-1/2}Y_{nt}$. Next recall the representation \eqref{E:4.8u} written, for $k\le n$, in the form 
\begin{equation}
%\label{E:4.8u}
X_n = X_k+\bigl[\chi_\epsilon\circ\tau_{X_k}(\omega)-\chi_\epsilon\circ\tau_{X_n}(\omega)\bigr]+
\biggl[\,\sum_{j=k}^{n-1}\epsilon\chi_\epsilon\circ\tau_{X_j}(\omega)\biggr]+(M_n^{(\epsilon)}-M_k^{(\epsilon)}).
\end{equation}
Given any $\eta>0$ and any $0<t_0<\dots<t_m$, Lemma~\ref{lemma-5.3} and the fact that $\chi_\epsilon$ is bounded (for~$\epsilon>0$ fixed) then show that
\begin{equation}
%\label{}
\,P_\omega^0\biggl(\max_{k=0,\dots,m}\Bigl| (Y_{nt_k}-Y_{n\delta})-\bigl(M_{N(n t_k)}^{(\epsilon)}-M_{N(n\delta)}^{(\epsilon)}\bigr)\Bigr|_2>\eta\sqrt n\biggr)\,\longrightarrow\,0
\end{equation}
in the limit as $n\to\infty$ followed by $\epsilon\downarrow0$ and by~$\delta\downarrow0$. In particular, it suffices to prove \twoeqref{E:5.35}{E:5.36} with~$t\mapsto B^{(n)}_t$  replaced by $t\mapsto n^{-1/2}M_{N(t)}^{(\epsilon)}$ in the limit~$n\to\infty$ and~$\epsilon\downarrow0$. (The fact that the law of~$X$, and thus also~$B^{(n)}$, has no dependence on~$\epsilon$ means that the $\epsilon\downarrow0$ limit must exist once $n\to\infty$ does. But we will prove this anyway.)

Returning to the discrete-time process, recall that $n\mapsto M_n^{(\epsilon)}$ is an $\R^d$-valued martingale with bounded increments. For each $a\in\R^d$, the Markov property of~$X$ along with the cocycle condition~\eqref{E:5.6i} give
\begin{equation}
\label{E:5.40}
E_\omega^0\bigl([a\cdot(M_k^{(\epsilon)}-M_{k-1}^{(\epsilon)})]^2\,\big|\,\FF_{k-1}\bigr) 
= f_\epsilon\circ\tau_{X_{k-1}}(\omega)
\end{equation}
for
\begin{equation}
%\label{}
 f_\epsilon(\omega'):=E_{\omega'}^0\bigl([a\cdot\psi_\epsilon(\omega',X_1)]^2\bigr),
\end{equation}
where~$\psi_\epsilon$ is as in \eqref{E:5.5i}. Since $ f_\epsilon\in\Cb(\Omega)$ and~$\omega\in\Omega^\star\cap\Omega_{p,q}$ for $p,q>1$ obeying \eqref{E:2.6}, Theorem~\ref{thm-2.1} implies
\begin{equation}
\label{E:5.27}
\frac1n\sum_{k=1}^{n} E_\omega^0\bigl([a\cdot(M_k^{(\epsilon)}-M_{k-1}^{(\epsilon)})]^2\,\big|\,\FF_{k-1}\bigr)\,\,\underset{n\to\infty}{\overset{P_\omega^0}\longrightarrow} a\cdot\Sigma^{(\epsilon)} a
\end{equation}
with~$\Sigma^{(\epsilon)}$ defined by
\begin{equation}
%\label{}
\forall a\in\R^d\colon\quad a\cdot\Sigma^{(\epsilon)} a =  
E_{\Q_\omega}( f_\epsilon)= \int 
E_{\omega'}^0\bigl([a\cdot\psi_\epsilon(\omega',X_1)]^2\bigr)\BbbP_\omega(\textd\omega').
\end{equation}
As $M^{(\epsilon)}$ has bounded increments, this verifies the conditions of the Martingale Functional CLT which then shows that $t\mapsto n^{-1/2}M_{nt}^{(\epsilon)}$ tends in law to  Brownian motion with covariance~$\Sigma^{(\epsilon)}$. The fact that $N(nt)/n\to t$ in probability as~$n\to\infty$ then extends the same conclusion to $t\mapsto n^{-1/2}M_{N(t)}^{(\epsilon)}$.

 In order to take the limit $\epsilon\downarrow0$ and identify the covariance via \eqref{E:1.5}, we observe
\begin{equation}
%\label{}
\begin{aligned}
 E_{\Q_\omega}( f_\epsilon)
&=\frac1{E_{\BbbP_\omega}\pi(0)}\,E_{\BbbP_\omega}\Biggl(\,\sum_{\begin{subarray}{c}x\in\Z^d\\(0,x)\in E(\Z^d)\end{subarray}}
\cc(0,x) |a\cdot x+\chi_\epsilon\circ\tau_x-\chi_\epsilon|^2\Biggr)
\\
&=\frac2{E_{\BbbP_\omega}\pi(0)}\bigl\langle a+\nabla (a\cdot\chi_\epsilon),\cc (a+\nabla (a\cdot\chi_\epsilon))\bigr\rangle_{L^2(\BbbP_\omega,\R^d)}.
\end{aligned}
\end{equation}
In particular, $E_{\Q_\omega}( f_\epsilon)$ is no smaller than the quantity in \eqref{E:1.5} for any~$\epsilon>0$.
Moreover, in light of \eqref{E:5.9i}, the inner product converges to
\begin{equation}
%\label{}
\bigl\langle a+\nabla (a\cdot\chi),\cc (a+\nabla (a\cdot\chi))\bigr\rangle_{L^2(\BbbP_\omega,\R^d)}
\end{equation}
as~$\epsilon\downarrow0$ with $x\mapsto x+\chi(\omega',x)$ $\cmss P_{\omega'}$-harmonic. This shows that
\begin{equation}
%\label{}
\bigl\langle a+\nabla (a\cdot\chi),\cc \nabla h\bigr\rangle_{L^2(\BbbP_\omega,\R^d)}=0
\end{equation}
for all~$h$ with~$\langle\nabla h,\cc\nabla h\rangle_{L^2(\BbbP_\omega,\R^d)}<\infty$ and so, for all such~$h$,
\begin{equation}
%\label{}
\begin{aligned}
\bigl\langle a+\nabla h,&\,\cc (a+\nabla h)\bigr\rangle_{L^2(\BbbP_\omega,\R^d)}
\\
&=\bigl\langle a+\nabla (a\cdot\chi),\cc (a+\nabla (a\cdot\chi))\bigr\rangle_{L^2(\BbbP_\omega,\R^d)}\\
&\qquad\quad+\bigl\langle \nabla h-\nabla (a\cdot\chi),\cc (\nabla h-\nabla (a\cdot\chi))\bigr\rangle_{L^2(\BbbP_\omega,\R^d)}
\\
&\ge\bigl\langle a+\nabla (a\cdot\chi),\cc (a+\nabla (a\cdot\chi))\bigr\rangle_{L^2(\BbbP_\omega,\R^d)}
\end{aligned}
\end{equation} 
proving that \eqref{E:1.5} is no smaller than $\lim_{\epsilon\downarrow0}E_{\Q_\omega}( f_\epsilon)$. It follows that $\lim_{\epsilon\downarrow0}E_{\Q_\omega}( f_\epsilon)$ coincides with the quantity in \eqref{E:1.5} for each~$a\in\R^d$. By polarization we get the pointwise convergence $\Sigma^{(\epsilon)}\to\Sigma$ and that of finite dimensional distributions as well. This completes the proof of our main theorem.
\end{proofsect}

%\newpage

\section{Exit time estimate}
\label{sec7}\noindent
Here we give the proof of the exit time estimate from Proposition~\ref{prop-4.1}. We follow the approach of Biskup, Chen, Kumagai and Wang~\cite{BCKW} whose main idea is to localize the effect of disordered environment to a finite box and thus embed the problem into a uniform setting to which general theory applies. 

Given $R\ge1$, set
\begin{equation}
\cc_\omega^R(e):=\begin{cases}
\cc_\omega(e),\qquad&\text{if }e\in E(\Lambda_{2R}),
\\
1,\qquad&\text{else},
\end{cases}
\end{equation}
and let
\begin{equation}
\pi_\omega^R(x):=\begin{cases}
\pi_\omega(x),\qquad&\text{if }x\in\Lambda_{2R},
\\
1+\pi_\omega(x),\qquad&\text{if }x\in\Lambda_{4R}\smallsetminus\Lambda_{2R},
\\
1,\qquad&\text{else}.
\end{cases}
\end{equation}
For $f\colon\Z^d\to\R$ with finite support, let %\bf CHECK IF NORMALIZATION NEEDED!! \rm
\begin{equation}
%\label{}
D_\omega^R(f,f):=\sum_{(x,y)\in E(\Z^d)}\cc_\omega^R(x,y)\bigl[f(y)-f(x)\bigr]^2
\end{equation}
be the Dirichlet form of~$f$ associated with the localized conductances. Write $\Vert\cdot\Vert_{\ell^p(\mu)}$ for the~$\ell^p$-norm on~$\Z^d$ relative to measure~$\mu$. The starting point is the following Nash-type inequality:

\begin{proposition}
\label{prop-Sob}
Let~$d\ge2$ and fix~$p,q\in(\frac d2,\infty)$ satisfying \eqref{E:2.6}. Then
\begin{equation}\label{r:remark}
\epsilon := 2\,\Bigl(\frac2d-\frac1p-\frac1q\Bigr)\Bigl(\frac{d-2}d+\frac1q\Bigr)^{-1}\in\bigl(0,\tfrac4{d-2}\bigr).
\end{equation}
 Moreover,  for each~$\omega\in\Omega_{p,q}$ there is $c_1(\omega)\in(0,\infty)$ such that
\begin{equation}
\label{e:Sob}
\|f\|^2_{\ell^{2+\epsilon}(\pi_\omega^R)}\le c_1(\omega) \left(R^{2-\frac{d\epsilon}{2+\epsilon}}D^R_\omega(f,f)+R^{-\frac{d\epsilon}{2+\epsilon}}\|f\|^2_{\ell^2(\pi_\omega^R)}\right)
\end{equation}
holds for all $f\in \ell^2(\pi_\omega^R)$.
\end{proposition}

This is, more or less, a restatement of~\cite[Proposition~3.3]{BCKW}. The proof of this proceeds by combining (via a mollifier technique) the conclusions of two Sobolev inequalities, one for functions supported in~$\Lambda_{8R}$ and the other for those supported outside~$\Lambda_{R}$. These can found in Lemmas~3.4 and~3.6 of~\cite{BCKW}, respectively. The constant $c_1(\omega)$ depends only on the spatial dimension~$d$, the parameters~$p$ and~$q$ and the suprema in the definition of~$\Omega_{p,q}$. The suprema enter the argument via \cite[Lemma~3.4]{BCKW}.

We will now derive some desired probabilistic consequences of \eqref{e:Sob}. For this, let~$Y^R$ be the continuous-time Markov chain on~$\Z^d$ with conductances $c^R_\omega$ and speed measure~$\pi^R_\omega$. 
A simple use of H\"older's inequality  turns \eqref{e:Sob} into
\begin{equation}\label{e:Sob1}
\|f\|^{4\frac{1+\epsilon}{2+\epsilon}}_{\ell^{2}(\pi_\omega^R)}\le c_1(\omega) \left(R^{2-\frac{d\epsilon}{2+\epsilon}}D^R_\omega(f,f)+R^{-\frac{d\epsilon}{2+\epsilon}}\|f\|^2_{\ell^2(\pi_\omega^R)}\right)\|f\|^{\frac{2\epsilon}{2+\epsilon}}_{\ell^{1}(\pi_\omega^R)}
\end{equation}
which is now a proper Nash inequality.
By general equivalence between heat-kernel bounds and Nash inequalities proved in Carlen, Kusuoka and Stroock~\cite[Theorem~(2.1)]{CKS}, the Nash inequality of the form (for $d'>0$ real)
\begin{equation}
\Vert f\Vert_2^{2+4/d'}\le A\Bigl(D(f,f)+\delta\Vert f\Vert_2^2\Bigr)\Vert f\Vert_1^{4/d'}
\end{equation}
yields a uniform bound on the heat kernel by $(d' A/t)^{d'/2}\texte^{\frac12\delta t}$. (Our Dirichlet form is normalized differently than in \cite{CKS}, so we adjusted the statement from \cite[Theorem~(2.1)]{CKS} accordingly.) 
Using this in \eqref{e:Sob} yields
\begin{equation}
\label{E:8.8}
P^x_\omega(Y^R_t = y)\le c_2(\omega) R^{-d}\Bigl(\frac t{R^2}\Bigr)^{-\frac{2+\epsilon}\epsilon}\texte^{\frac12t R^{-2}}\pi_\omega^R(y),
\end{equation}
where $c_2(\omega)$ is simple function of $c_1(\omega)$; see \cite[Lemma~3.7]{BCKW}.

Another argument from \cite{CKS} based on the so called Davies method~\cite{Davies} upgrades \eqref{E:8.8} to the following off-diagonal bound:

\begin{proposition}
\label{prop-8.2}
Let~$d\ge2$ and fix~$p,q\in(\frac d2,\infty)$ satisfying \eqref{E:2.6} and let~$\epsilon\in(0,\frac4{d-2})$ be as in~\eqref{r:remark}. For each~$\kappa\in(0,1]$ and $\omega\in\Omega_{p,q}$ there exists~$c_3(\omega)\in(0,\infty)$ such that
\begin{equation}
\label{E:8.9}
P^x_\omega(Y^R_t = y)\le c_3(\omega) R^{-d}\Bigl(\frac{t}{R^2}\Bigr)^{-\frac{2+\epsilon}{\epsilon}} \exp\left(-\frac{|x-y|_1}{5\kappa R}\log \left(\frac{R^2}{t}\right)\right)\pi_\omega^R(y)
\end{equation}
holds for all~$R\ge1$ with~$\kappa R\ge1$, all $t\in(0,R^2)$ and all $x,y\in\Z^d$.
\end{proposition}

\noindent
This is a restatement of \cite[Proposition~3.8]{BCKW} except that we keep the~$\omega$ dependence in the prefactor~$c_3(\omega)$ which is a function of the previous $\omega$-dependent constants and~$\kappa$, but not of~$R$,~$t$,~$x$ or~$y$. The reason for requiring $\kappa R\ge1$ is that \cite{BCKW} deals primarily with walks that can take arbitrarily long jumps; the value~$\kappa R$ is then a cutoff for jumps of ``small' size. With Proposition~\ref{prop-8.2} in hand, we are ready to give:

\begin{proofsect}{Proof of Proposition~\ref{prop-4.1}}
We proceed as in the proof of \cite[Proposition 3.9]{BCKW} with the main computation actually drawn from \cite[Lemma 3.10]{BCKW}. Let $d\ge2$ and assume $p,q>d/2$ obey \eqref{E:2.6}. Noting that the probability in question is non-decreasing in~$\sigma$, assume $\sigma\in(0,1/2)$. Fix~$\alpha>0$ and~$\omega\in\Omega_{p,q}$. Most of the proof is carried out for the process~$Y^R$; we will return to~$X$ only in the very last step. 

For each~$R\ge1$ and $n\ge1$, denote $A_{n,R}(x):=\{y\in\Z^d\colon n R\le |y|_1\le (n+1) R\}$. Then any~$\kappa\in(0,1]$, the bound \eqref{E:8.9} gives
\begin{equation}
%\label{}
P^x_\omega\bigl(Y^{R}_t\not\in\Lambda_{\sigma R}(x)\bigr)\le c_3(\omega) R^{-d}\Bigl(\frac{t}{R^2}\Bigr)^{-\frac{2+\epsilon}{\epsilon}}
\sum_{n\ge1}\texte^{-\frac\sigma{5\kappa}\log(R^2/t)n}\sum_{y\in A_{n,\sigma R}(x)}\pi_\omega(y)
\end{equation}
for all~$R\ge1$ with~$\kappa R\ge1$, all $t\in(0,R^2)$ and all $x\in\Z^d$.
Assuming $t\le R^2/3$ to ensure that $\log(R^2/t)\ge1$, note that for $x\in\Lambda_R$, we have $A_{n,\sigma R}(x)\subseteq \Lambda_{R+(n+1)\sigma R}\subseteq\Lambda_{(n+2)R}$ and $|\Lambda_{(n+2)R}|\le 5^d (n+2)^d R^d$. The Cauchy-Schwarz inequality then gives
\begin{equation}
\label{E:6.11}
P^x_\omega\bigl(Y^{R}_t\not\in\Lambda_{\sigma R}(x)\bigr)\le
c_4(\omega) \Bigl(\frac{t}{R^2}\Bigr)^{-\frac{2+\epsilon}{\epsilon}}
\biggl(\,\sum_{n\ge1}\texte^{-\frac\sigma{5\kappa}\log(R^2/t)n}\biggr)^{1/2},
\end{equation}
where
\begin{equation}
%\label{}
c_4(\omega):=5^d c_3(\omega)\biggl(\,\sum_{r\ge1}(r+2)^{2d}\texte^{-\frac\sigma{5\kappa}\log(R^2/t)r}\biggr)^{1/2}\sup_{n\ge1}\frac1{|\Lambda_n|}\sum_{y\in\Lambda_n}\pi_\omega(y).
\end{equation}
Bounding the geometric series in \eqref{E:6.11} under the assumption that $\kappa/\sigma$ is so small that $\texte^{-\frac\sigma{5\kappa}\log(R^2/t)}\le1/2$ shows
\begin{equation}
\label{E:8.13}
P^x_\omega\bigl(Y^{R}_t\not\in\Lambda_{\sigma R}(x)\bigr)\le
 c_5(\omega)\Bigl(3\frac{t}{R^2}\Bigr)^{\frac\sigma{10\kappa}-\frac{2+\epsilon}{\epsilon}},
\end{equation}
where $c_5(\omega)=\min\{1,\sqrt2\,c_4(\omega)\}$ and where we assumed the exponent to be positive. The factor~$3$ ensures that the bound holds with no restrictions on~$t$ relative to~$R$. In the rest of the proof, we also assume that
\begin{equation}
%\label{}
\frac\sigma{10\kappa}-\frac{2+\epsilon}{\epsilon}>\alpha
\end{equation}
which is achieved by choosing~$\kappa$ small enough. The restriction $\kappa R\ge1$ then forces us to assume that $R\ge R_0:=\kappa^{-1}$.

Denote $ H_\Lambda(Y):=\inf\{t\ge0\colon Y_t\not\in\Lambda\}$ and $ H_\Lambda(Y^R):=\inf\{t\ge0\colon Y_t^R\not\in\Lambda\}$  and pick any $x\in\Lambda_{R/2}$.
Since~$Y$ and~$Y^R$  started from~$x$  have the same law until the first exit time from~$\Lambda_R$, from $\sigma<1/2$ we get
\begin{equation}
\begin{aligned}
 P^x_\omega\bigl( H_{\Lambda_{\sigma R}(x)}&(Y)\le 2t\bigr)
=P^x_\omega\bigl( H_{\Lambda_{\sigma R}(x)}(Y^R)\le 2t\bigr)
\\
&\le P^x_\omega\bigl(Y_{3t}^R\not\in\Lambda_{\sigma R/2}(x)\bigr)
+P^x_\omega\Bigl(Y_{3t}^R\in\Lambda_{\sigma R/2}(x)\,\wedge\,H_{\Lambda_{\sigma R}(x)}(Y^R)\le 2t\Bigr).
\end{aligned}
\end{equation}
The first probability is estimated directly from \eqref{E:8.13} while for the second  probability  we invoke the strong Markov property at time $ H_{\Lambda_{\sigma R}(x)}(Y^R)$ to bound it by
\begin{equation}
%\label{}
 \max_{z\in\partial\Lambda_{\sigma R}(x)}\,\sup_{0<s\le 3t}P^z_\omega\bigl(Y_{s}^R\not\in\Lambda_{\sigma R/2}(z)\bigr),
\end{equation}
where $\partial\Lambda$ is the set of vertices in~$\Lambda^\cc$ that have an edge to~$\Lambda$ and where we noticed that $\Lambda_{\sigma R/2}(x)\cap\Lambda_{\sigma R/2}(z)=\emptyset$ for all~$z\in\partial\Lambda_{\sigma R}(x)$.  This  is now estimated via \eqref{E:8.13} as well, uniformly in~$z$.  Using the definition of~$Y$ via  the discrete-time walk~$X$  and an independent Poisson process~$t\mapsto N(t)$ we get
\begin{equation}
 P^x_\omega\bigl(H_{\Lambda_{\sigma R}(x)}(Y)\le 2t\bigr)
\ge P\bigl(N(2t)>t\bigr) P^x_\omega\bigl(\ttau_{\Lambda_{\sigma R}(x)}\le t\bigr)
\end{equation}
 uniformly in~$x\in\Lambda_{R/2}$. The claim \eqref{E:2.32} now  follows by noting that the first probability on the right is uniformly positive for~$t\ge1$. (There is nothing to prove for~$t<1$.) 
\end{proofsect}

%\newpage
\section{Extensions and open problems}
\label{sec8}\noindent
We will close this paper by a discussion of open problems and extensions that naturally arise from this work. Let us first address the limitations of the present setting.

\subsection{Beyond moment conditions and ergodicity}
We believe that, just as in the results based on stochastic homogenization (particularly, references \cite{ADS15} and \cite{BCKW}), the moment conditions expressed by \eqref{E:2.6} are likely artifacts of the method of proof. In fact, we will go as far as to pose:

\begin{conjecture}
%\label{conj}
In all $d\ge1$, an IIP holds for all $\omega\in\Omega^\star\cap\bigcup_{p,q>1}\Omega_{p,q}$.
\end{conjecture}

\noindent
Note that, in spite of some temptation to the contrary, we do not not include~$\Omega_{1,1}$. This is because the $L^1(P_\omega^0)$-convergence in Theorem~\ref{thm-2.1} is not automatic under $\Q_\omega$-integrability of~$f$ alone. (Note that, even in the stochastic seting, this calls for the use of the Dominated Ergodic Theorem that needs~$f$ to be in at least $L\log L$ instead of just~$L^1$. This matters at times; see~\cite[Remark~3.3]{BLRV} or \cite[Lemma~2.6]{B11}.)

Requiring the containment in~$\bigcup_{p,q>1}\Omega_{p,q}$ is natural as without the first moment conditions we may lose tightness and/or non-degeneracy of the limiting process (see  Barlow, Burdzy and Tim\'ar~\cite{BBT} or the examples in Biskup~\cite{B19}). Unfortunately, except in~$d=1,2$ the status of the conjecture is not clear even for samples from ergodic measures, let alone generic members of~$\Omega^\star$.

Since we are  discussing borderline situations  where an IIP can be expected, the reader may wonder whether the assumption of ergodicity of~$\omega$ (defined by ergodicity of~$\BbbP_\omega$) is necessary for an IIP to hold. It seems reasonable to ask:

\begin{problem}
Characterize the averaging non-ergodic $\omega\in\Omega$, or at least a ``robust'' subset thereof, for which an IIP provably holds.
\end{problem}

The lack of ergodicity for averaging~$\omega$ is tantamount to~$\omega$ being dominated by large regions where averaging does take place but with limit values that depend on the region.
An approximate IIP likely holds in each of these regions  and further averaging will presumably take place as soon as the regions are of diameters that grow sublinearly with distance. The challenge is to find a way to build a proof based on this intuition.

\subsection{Random walk in subdomains of~$\Z^d$}
As was pointed out  to us  by Jean-Dominique~Deuschel, the proofs of our main result seem to generalize to the situation when the underlying environment  changes  with the time scale of the walk. This should upgrade the main result to the form:

\begin{conjecture} 
\label{cor-2.5}
Let $p,q>1$ obey \eqref{E:2.6} and let~$\{\omega_n\}_{n\ge1}$ be a sequence of environments from~$\Omega$ such that the following holds:
%\settowidth{\leftmargini}{(11)}
\begin{enumerate}
\item[(1)] $\sup_{L\ge1}\sup_{n\ge1}(Ln)^{-d}\sum_{e\in E(\Lambda_{ Ln})}\cc_{\omega_n}(e)^p<\infty$,
\item[(2)] $\sup_{L\ge1}\sup_{n\ge1}(Ln)^{-d}\sum_{e\in E(\Lambda_{ Ln})}\cc_{\omega_n}(e)^{-q}<\infty$,
\item[(3)] $\forall f\in\Cloc(\Omega)\,\exists\ell(f)\in\R\,\forall L>0\colon$
\begin{equation}
%\label{}
\lim_{n\to\infty}\frac1{|\Lambda_{L\sqrt n}|}\sum_{x\in\Lambda_{L\sqrt n}}f\circ\tau_x(\omega_n)=\ell(f),
\end{equation}
\item[(4)] the probability measure~$\BbbP$ on~$(\Omega,\FF)$ defined implicitly by $\forall f\in\Cloc(\Omega)\colon E_\BbbP(f)=\ell(f)$ is ergodic w.r.t.\ translations.
\end{enumerate}
Then, as $n\to\infty$, the law of $t\mapsto B_t^{(n)}$ under $P^0_{\omega_n}$ 
tends to that of a centered Brownian motion with a non-degenerate covariance structure.
\end{conjecture}

\noindent
Note that the environments $\{\omega_n\}_{n\ge1}$ are \emph{static}. All that we allow is that the walk at time scale~$n$ sees a different environment than that on time scale~$n+1$.

A situation where such an $n$-dependent setting is useful is the random walk confined to  subdomains of~$\Z^d$. Examples of these are half-spaces, quarter spaces, etc, which are domains of the form
\begin{equation}
\label{E:1.17}
\V:=\bigl\{(x_1,\dots,x_d)\in\Z^d\colon x_1,\dots,x_{d_1}\ge0\bigr\}
\end{equation}
for  natural~$d_1\ge1$ such that $d-d_1\ge1$.  The Markov chain then moves as before inside~$\V$ except that the moves outside~$\V$ are suppressed. 

For stochastic environments, including supercritical bond percolation, an IIP has been proved for these geometries by Chen, Croydon and Kumagai~\cite{CCK15}. The conclusion is that the resulting Brownian motion (still with covariance structure determined by full-space homogenization) now reflects on the boundaries of
\begin{equation}
\label{E:7.3a}
\bigl\{(x_1,\dots,x_d)\in\R^d\colon x_1,\dots,x_{d_1}\ge0\bigr\}.
\end{equation}
A  technical challenge  is that, while the environment~$\omega$ is kept fixed relative to~$\V$, the initial state of the walk is allowed to ``slide off'' to infinity inside~$\V$ at the diffusive scale; i.e., along points $z_n:=\lfloor z\sqrt n\rfloor$ for some~$z$ in the set \eqref{E:7.3a}. Thus, even the control of the walk until its first hit of the boundary of~$\V$ requires proving a so called \emph{arbitrary starting point} version of a Quenched Invariance Principle, i.e., one for environments $\omega_n:=\tau_{z_n}(\omega)$. The relevance of this problem was first pointed out by Rhodes~\cite{Rhodes}. 

 The statement of  Conjecture~\ref{cor-2.5} accommodates these situations because the reflection can be imitated by extending the environment in~$E(\V)$ to all of~$E(\Z^d)$ via reflections in planes of the form $-\frac12e_i+\{(x_1,\dots,x_d)\colon x_i=0\}$ where~$i=1,\dots,d_1$. The corresponding full-space random walk in the reflected environment then coincides with the random walk in~$\V$. As soon as  we know  that the former scales to Brownian motion~$B$, the walk in~$\V$ scales to the reflected version of~$B$. (For the walk started from the origin, this in fact follows already from Theorem~\ref{thm-2.4}.)

The reflection tricks are slick and informative but we feel that the current approach has the potential to solve the problem in more general domains as well:

\begin{problem}
Construct a ``large'' set~$\Omega_1\subseteq\Omega$ such that, for each~$\omega\in\Omega_1$, each bounded Lipschitz domain~$D\subseteq\R^d$ and each~$z\in D$, the law of~$B^n$ under environment~$\omega$ of the walk started from~$z_n:=\lfloor z\sqrt n\rfloor$ and reflected on the boundary of~$\D_n:=\{x\in\Z^d\colon x/\sqrt n\in D\}$ scales to a reflected Brownian motion in~$D$ started from~$z$.
\end{problem}

\noindent
Here the attribute ``large'' again refers to the property that~$\Omega_1$ is translation invariant and of full mass under all ergodic laws on the conductances. We expect that all that will ultimately be needed is that the condition in Definition~\ref{def-1.2} be upgraded to sequences of volumes of the form $\lfloor zn\rfloor+\Lambda_n$, for $z\in\R^d$. We may in fact even consider the situation that the limit measure then depends on~$z$. This would allow for inclusion of spatially modulated environments.

\subsection{Long-range and time-dependent walks}
Throughout this work we have confined ourselves to models with nearest-neighbor jumps. However, the idea to take on this problem was partially inspired by recent joint work of the author with X.~Chen, T.~Kumagai and J.~Wang~\cite{BCKW} that studied walks with jumps of arbitrary size. Based on this work, we believe in:

\begin{conjecture}
Let~$d\ge2$ and $p,q>d/2$ be such that $1/p+1/q<2/d$. Then an IIP holds for random walk among long-range conductance configurations $\{\cc(x,y)=\cc(y,x)\colon x,y\in\Z^d\}$ that are averaging and ergodic (in spirit of Definitions~\ref{def-1.2} and~\ref{def-1.4}) and obey
\begin{equation}
%\label{}
\sup_{n\ge1}\frac1{|\Lambda_n|}\sum_{x\in\Lambda_n}\Bigl(\sum_{y\in\Z^d}\cc(x,y)|y-x|^2\Bigr)^p<\infty
\end{equation}
and
\begin{equation}
%\label{}
\sup_{n\ge1}\frac1{|\Lambda_n|}\sum_{(x,y)\in E(\Lambda_n)}\cc(x,y)^{-q}<\infty.
\end{equation}
Here $E(\Lambda)$ is the set of nearest-neighbor edges incident with~$\Lambda$. (We require nearest-neighbor conductances to be strictly positive.)
\end{conjecture}

\noindent
We find this worthwhile because the setting includes models of long-range percolation (superimposed on already connected~$\Z^d$). However, unlike the stochastic setting, the above allows for zero-density modifications of the environment and so, in particular, for spatially inhomogeneous truncation of long jumps.  This matters in light of  the obstructions posed by the lack of everywhere sublinearity of the corrector (cf~\cite[Theorem 2.5]{BCKW})  which prevent us from addressing the full range of long-range percolation models in which an IIP is expected to hold. 

Another area where our approach should be applicable are continuous time random walks in time dependent environments. This particularly concerns the variable-speed model with the generator
\begin{equation}
\label{E:7.3}
\cmss L_t f(x) = \sum_{y\colon (x,y)\in E(\Z^d)}a_t(x,y)\bigl[f(y)-f(x)\bigr],
\end{equation}
where $t\mapsto a_t(e)$ is, for each edge~$e$, a positive function on~$\R$. Here the work of Andres~\cite{A14} and later Andres, Chiarini, Deuschel and Slowik~\cite{ACDS18} shows that an IIP holds for a.e.\ sample from a space-time ergodic random environment subject to a variant of the $p$,$q$-condition. We take a bold step and pose:

\begin{conjecture}
Let $d\ge2$ and assume that  a locally Lebesgue integrable function $t\mapsto a_t(e)\in(0,\infty)$ is given for each~$e\in E(\Z^d)$.  Writing $\tau_{t,x}$ for the space-time shift acting as $[\tau_{t,x}(a)]_s(y,z) = a_{t+s}(y+x,z+x)$, suppose that $t,e\mapsto a_t(e)$ is averaging in the sense that for each bounded and continuous function~$f$ on the space of environments,
\begin{equation}
%\label{}
\lim_{n\to\infty}
\,\frac1{n^2}\int_0^{n^2}\frac1{|\Lambda_n|}\sum_{x\in\Lambda_n}f\circ\tau_{t,x}(a)\,\textd t
\quad\text{\rm\ exists},
\end{equation}
 the law~$\BbbP$ defined from these limits as in Proposition~\ref{prop-2.2} is ergodic under space-time translations and, in addition,
\begin{equation}
%\label{}
\sup_{n\ge1}\frac1{n^{d+2}}\int_0^{n^2} \sum_{e\in E(\Lambda_n)}a_t(e)^p\,\textd t<\infty
\end{equation}
and
\begin{equation}
%\label{}
\sup_{n\ge1}\frac1{n^{d+2}}\int_0^{n^2} \sum_{e\in E(\Lambda_n)}a_t(e)^{-q}\,\textd t<\infty
\end{equation}
hold for for some $p,q>1$ satisfying
\begin{equation}
\label{E:7.9i}
\frac1{p-1}+\frac1{q(p-1)}+\frac1q<\frac2d.
\end{equation}
Then an IIP holds for the continuous-time Markov chain on~$\Z^d$ with generator \eqref{E:7.3}.
\end{conjecture}

The reason why we find condition \eqref{E:7.9i}, which is exactly that of \cite{ACDS18}, reasonable is because this is also the condition under which a local CLT holds in stochastic context (Andres, Chiarini and Slowik~\cite{ACS21}). An important difference from our problem is that the random walk with generator \eqref{E:7.3} is generally not reversible. Due to the time dependence, the question is non-trivial even in~$d=1$ (see Deuschel and Slowik~\cite{DS16}, Biskup~\cite{B19} and Biskup and Pan~\cite{B-Pan}).

\subsection{Quantitative theory}
As discussed in Section~\ref{sec-2.1}, a great deal of progress has been made in recent years on quantitative versions of stochastic homogenization. While the present work aims at an extreme end of qualitative theory, there is potentially a new perspective on quantitative  homogenization as well, this time focused on deterministic environments.  Motivated by conclusions obtained by Mourrat~\cite{Mourrat1,Mourrat2}, we pose: 

\begin{problem}
Given a specific polynomial rate of convergence in \eqref{E:1.1},  resp.,~\eqref{E:1.9i},   derive a rate of convergence in \eqref{E:4.4i} or even in the IIP. Work under uniform ellipticity if convenient.
\end{problem}

\noindent
We note that quantitative stochastic homogenization results exist even in non-elliptic settings such as the random walk on the supercritical percolation cluster (e.g., Armstrong and Dario~\cite{AD} and Dario~\cite{Dario}). A quantitative theory might in fact be the right way to address degenerate situations represented by~$\Omega:=\{0,1\}^{E(\Z^d)}$ which includes the random walk on percolation cluster.

An inherent part of both qualitative and quantitative homogenization is control of the corrector~$\chi$. For i.i.d.\ environments,~$\chi$ oscillates quite heavily with the overall fluctuation structure close to that of the Gaussian Free Field (see Mourrat and Otto~\cite{MO16b}, Mourrat and Nolen~\cite{MN17}, Gu and Mourrat~\cite{GM16}). This indicates that the dependence of~$\chi$ on the conductances is likely too irregular to make it defined for all configurations for which an IIP holds. The following problem has been on the author's bucket list for a while:

\begin{problem}
Find a ``large'' subset~$\Omega^{\star\star}$ of~$\Omega^\star$ such that the corrector is  well defined and (reasonably) ``well behaved'' for all~$\omega\in\Omega^{\star\star}$.
\end{problem}

\noindent
Here the word ``large'' definitely includes the attribute that~$\Omega^{\star\star}$ is $\FF$-measurable and of full measure under each ergodic law of the conductances. The attribute ``well behaved'' in turn refers to the properties relevant for the use of the corrector; e.g., sublinearity at large separations.

\subsection{Other transition mechanisms}
The final direction we find worthy of exploring concerns different transition mechanisms of the Markov chain. While the question ``what deterministic choices of the transition probabilities lead to the usual behavior'' applies in full generality of random walks in random environment, we will focus mainly on one recently studied case.

A \emph{balanced random walk} is a Markov chain on~$\Z^d$ whose transition probability satisfies the symmetry requirement
\begin{equation}
\label{E:7.1}
\forall x,z\in \Z^d\colon\, \cmss P(x,x+z)=\cmss P(x,x-z).
\end{equation}
Typically, we also assume that there exists a finite set~$\Lambda\subseteq\Z^d$ of ``full dimensionality'' such that $\cmss P(x,x+z)=0$ unless~$z\in\Lambda$. Under this condition, the Markov chain~$X$ is a martingale with bounded increments and so an IIP will hold as soon as we verify the conditions of the Martingale Functional CLT. This boils down to a proof of averaging for the ``point of view of the particle'' (i.e., an analogue of Theorem~\ref{thm-2.1}).

The stochastic approach proceeds by constructing an invariant measure~$\Q$ for the ``point of view of the particle'' and proving that~$\Q$ is equivalent to the \emph{a priori} law of the environment; the desired averaging then follows from the Ergodic Theorem (in time). This approach was pioneered by Lawler~\cite{L82} in uniformly elliptic cases and extended beyond uniform ellipticity by Guo and Zeitouni~\cite{GZ12} under suitable moment conditions and by Berger and Deuschel~\cite{BD14}, who even allow degenerate~$\cmss P$.

We believe that progress can made for deterministic uniformly elliptic environments (where~$\omega$ represents the collection $\{\cmss P(x,x+z)\colon x\in\Z^d, z\in\Lambda\}$). Here Mustapha~\cite{Mustapha} proved the existence and uniqueness of a function $\frakm_\omega\colon\Z^d\to[0,\infty)$ solving $\cmss L^\star \frakm_\omega=0$, where~$\cmss L:=\cmss P-1$, subject to the normalization~$\frakm_\omega(0)=1$. In addition, he also claims a uniform heat-kernel bound of the form
\begin{equation}
%\label{}
P^0(X_n=x)\le \frac{c_1}{\frakm_\omega(\Lambda_{\sqrt n})}\texte^{-c_2|x|^2/n} \,\frakm_\omega(x),
\end{equation}
where $\frakm_\omega(A):=\sum_{x\in A}\frakm_\omega(x)$. With this in hand, the proof of an IIP boils down to proving that, for~$\omega$ uniformly elliptic, averaging and ergodic (again, in spirit of Definitions~\ref{def-1.2} and~\ref{def-1.4}), there exists a probability measure~$\Q_\omega$ on~$(\Omega,\FF)$ such that 
\begin{equation}
\label{E:7.13}
f\in \Cb(\Omega)\colon\quad
\frac1{\frakm_\omega(\Lambda_r)}\sum_{x\in\Lambda_n}\frakm_\omega(x) f\circ\tau_x(\omega)\,\underset{n\to\infty}\longrightarrow\,E_{\Q_\omega}(f)
%\quad\text{\rm\ in C\'esaro sense}
\end{equation}
and
\begin{equation}
%\label{}
\forall f\in \Cb(\Omega)\colon\, E_{\Q_\omega}(f)=0\,\,\Rightarrow\,\,\liminf_{\epsilon\downarrow 0} E_{\Q_\omega}\Bigl(\epsilon\bigl|(\epsilon-\LL)^{-1}f\bigr|\Bigr)=0,
\end{equation}
where~$\LL$ is the generator of the chain from the ``point of view of the particle.'' Indeed, these are exactly the conditions that make the conclusion of Theorem~\ref{thm-2.1}  true. 

Another recently studied transition mechanism is that of \emph{doubly stochastic} or \textit{div-free} environments. These are exactly those environments for which the \emph{a priori} translation-invariant law remains stationary for the ``point of view of the particle.'' (The conductance models are thus special cases of this.) An early work on an IIP the case of such environments that admit a bounded cycle representation is that of Deuschel and K\"osters~\cite{DK08}. A more general approach based on a technical $\HH_{-1}$-condition has been developed by T\'oth~\cite{T18} (drawing on Kozma and T\'oth~\cite{KT17}). 

Unfortunately, the $\HH_{-1}$-condition seems difficult to verify so, while extending these results to a ``deterministic'' setting seems a valuable task, the likely first cases to try are those admitting a bounded cycle decomposition.

 As a final remark, inserted per suggestion from a referee, we note the recent work of Gwynne, Miller and Sheffield~\cite{GMS} on random walk in scale-free random environments. The proof of an invariance principle there is based on a novel way of averaging that leverages (properly formulated) assumptions of scale invariance. It would be of interest to see if the methods of the present work can be translated into that context as well.

%\newpage
\section*{Acknowledgments}
\noindent
This research has been partially supported by NSF award DMS-1954343. Discussions with Jean-Dominique Deuschel are greatly appreciated and so are the comments received from Scott Armstrong, Claude Le Bris, Mitia Duerinckx, Jean-Christophe Mourrat and Felix Otto.

\bibliographystyle{abbrv}

\end{document}